\documentclass[11pt,oneside]{amsart}
\usepackage{epsfig,psfrag,subfigure}
\usepackage{amsmath,amssymb,amsfonts,latexsym, mathrsfs}
\usepackage{graphicx,graphics,cite}
\usepackage{wrapfig}
\usepackage[usenames]{color}
\usepackage{geometry}                
\usepackage{colortbl,slashbox,booktabs,ctable,multirow}
\usepackage[table]{xcolor}
\usepackage{epstopdf}

\topmargin by -\headsep \textheight 8.9in \oddsidemargin 0.15in
\makeatletter \@addtoreset{figure}{section}
\def\thefigure{\thesection.\@arabic\c@figure} \def\fps@figure{h, t}
\@addtoreset{table}{section}
\def\thetable{\thesection.\@arabic\c@table}\def\fps@table{h, t}
\@addtoreset{equation}{section}

\makeatother
\newtheorem{theorem}{Theorem}[section]
\newtheorem{lemma}{Lemma}[section]
\newtheorem{proposition}{Proposition}[section]
\newtheorem{corollary}{Corollary}[section]
\newtheorem{algorithm}{Algorithm}[section]

\newcommand{\beq}{\begin{equation}}
\newcommand{\eeq}{\end{equation}}
\newcommand{\beqa}{\begin{eqnarray}}

\newcommand{\eeqa}{\end{eqnarray}}

\newcounter{nfig}

\begin{document}

\title{ML($n$)BiCGStab: Reformulation, Analysis and Implementation}
\author{Man-Chung Yeung\\ \\
Department of Mathematics, University of Wyoming, Laramie, 82071, USA}
\date{}
\maketitle

\begin{center}
{\it Dedicated to the Memory of Gene Golub}\footnote{This paper was presented in Gene Golub Memorial Conference, Feb. 29-Mar. 1, 2008,  University of Massachusetts, Dartmouth, U.S.A..}
\end{center}

\begin{abstract}
With the aid of index functions, we re-derive the ML($n$)BiCGStab
algorithm in \cite{yeungchan} in a more systematic way. It turns out that there are
$n$ ways to define the ML($n$)BiCGStab residual vector. Each
definition will lead to a different ML($n$)BiCGStab
algorithm.
We demonstrate this by presenting a second algorithm which
requires less storage.
In theory, this second algorithm serves as a bridge
that connects the Lanczos-based BiCGStab and the Arnoldi-based FOM while ML($n$)BiCG a bridge connecting BiCG
and FOM.
We also
analyze the breakdown situations from the probabilistic point of
view and summarize some useful properties of ML($n$)BiCGStab.
Implementation issues are also addressed.
\end{abstract}




\section{Introduction}
\label{Introduction}
If we express the BiCG\cite{fletcher, lanczos} residual as ${\bf
r}_k^{BiCG} = p_k({\bf A}) {\bf r}_0$ in terms of a polynomial,
the residual vector ${\bf r}_k$ of a Lanczos-type product
method\footnote{For this type of Krylov subspace methods, one can
consult \cite{gut10}. They are called hybrid BiCG methods in
\cite{SVF}.} based on BiCG is defined to be ${\bf r}_k
= \phi_k({\bf A}) p_k({\bf A}) {\bf r}_0$ where
$\phi_k(\lambda)$ is some polynomial of degree $k$ with $\phi_k(0) =
1$. In CGS\cite{sonn}, $\phi_k = p_k$.
Since, in every iteration, CGS searches for an approximate solution
in a larger Krylov subspace, it often converges much faster than
BiCG. However, CGS usually behaves irregularly due to a lack of a
smoothing mechanism. In BiCGStab\cite{van}, 
the $\phi_k$ is
\begin{equation}\label{equ:12-23}
\phi_k(\lambda) = \left\{ \begin{array}{lcl} 1& & \mbox{if } k = 0\\
(\rho_k \lambda + 1) \phi_{k-1}(\lambda)& & \mbox{if } k > 0.
\end{array} \right.
\end{equation}
Here $\rho_k$ is a free parameter selected to minimize the
$2$-norm of ${\bf r}_k^{BiCGStab}$ in the $k$-iteration.
As a
result, BiCGStab is generally more stable and robust than CGS.
BiCGStab has been
extended to BiCGStab2\cite{gut} and BiCGStab($l$)\cite{SF, SVF}
through the use of higher degree minimizing polynomials. In
BiCGStab2, the $\phi_k$ is defined by the recursion
$$
\phi_k(\lambda) = \left\{
\begin{array}{lll}
1 & & \mbox{if } k = 0
\\
(\rho_k \lambda + 1) \phi_{k-1}(\lambda) & & \mbox{if } k \mbox{ is
odd}\\
((\alpha_k \lambda + \beta_k) (\rho_{k-1} \lambda + 1) + 1- \beta_k)
\phi_{k-2}(\lambda) & &\mbox{if } k \mbox{ is even}.
\end{array}\right.
$$
The parameters are again chosen to minimize BiCGStab2 residuals.
Likewise, BiCGStab($l$) defines its $\phi_k$
as
$$\phi_k(\lambda) = \left\{ \begin{array}{lcl}
1& &\mbox{if } k = 0\\
 (1+\sum_{j = 1}^l \alpha_j
\lambda^j)\phi_{k-l}(\lambda)& & \mbox{if } k \mbox{ is a multiple
of } l \end{array}\right.
$$
where the parameters in the factor $1+\sum_{j = 1}^l \alpha_j
\lambda^j$ yields an $l$-dimensional minimization in every $l$th step.
 BiCGStab2 and BiCGStab($l$) usually converge
faster than BiCGStab because of smaller residuals in magnitude while
avoiding near-breakdowns caused by a possibly too small $\rho_k$. CGS,
BiCGStab and BiCGStab2 have been summarized and generalized by
GPBi-CG\cite{zhang}. Here $\phi_k$ is
$$\phi_k(\lambda) = \left\{ \begin{array}{lcl}
1& &\mbox{if } k = 0\\
 \rho_1 \lambda + 1& & \mbox{if } k = 1 \\
 (\alpha_k \lambda + 1 + \beta_k) \phi_{k-1}(\lambda) - \beta_k \phi_{k-2}(\lambda)& &\mbox{if } k >
 1.
 \end{array}\right.
$$
GPBi-CG will become CGS, BiCGStab or BiCGStab2 when the $\alpha,
\beta, \rho$ are appropriately chosen. For detailed descriptions of
these and other product-type methods, one is referred to \cite{gut30, gut31, saad,
savan, van1} and the references therein. Moreover, a history of product-type methods can be found in
\cite{gut1}. The history starts three decades ago with
IDR\cite{ws} method which can be
considered as the predecessor of CGS and
BiCGStab\cite{ssg}. Recently, IDR has been generalized to
IDR($s$) with a shadow space of higher dimension, see \cite{
ssg, sonn1, gs10}. IDR($s$) has close relations with ML($s$)BiCGStab.

Generalizations of BiCGStab
to methods
based on generalizations of BiCG have been made. For example,
BL-BiCGStab\cite{gue} is a BiCGStab variant built on the
BL-BiCG\cite{ol} for the solution of systems with multiple
right-hand sides. ML($n$)BiCGStab\cite{yeungchan} is another
BiCGStab variant built on ML($n$)BiCG, a BiCG-like method derived
from a variant of the band Lanczos process described in
\cite{abfh96} with $n$ left-starting vectors and a single
right-starting vector.

The derivation of the ML($n$)BiCGStab algorithm in \cite{yeungchan}
was complicated. In this paper, we exploit the concept of
index functions to re-derive the algorithm in a more systematic way, step by step.
Index functions were introduced in \cite{yeungboley} by Boley for
the purpose of simplifying the development of the transpose-free
multiple starting Lanczos process, and they proved to be very
helpful.

It turns out that the definition of the ML($n$)BiCGStab residual
vector ${\bf r}_k$ in \cite{yeungchan} is not unique. There are at
least $n$ different ways to define ${\bf r}_k$. Let $\widehat{\bf
r}_k$ be the residual of ML($n$)BiCG and $\phi_k(\lambda)$ as in
(\ref{equ:12-23}). Then, the ML($n$)BiCGStab residual ${\bf r}_k$
in \cite{yeungchan} is defined as
\begin{equation} \label{equ:607}
{\bf r}_{k} = \phi_j ({\bf A}) \,\widehat{\bf r}_{k}
\end{equation}
where $k = jn + i,\,\, 1 \leq i \leq n,\,\, j = 0, 1, 2, \cdots$.
Starting from $k = 1$, let us call every $n$ consecutive iterations
an iteration ``cycle''. For example, iterations $k = 1, 2, \cdots, n$
form the first cycle, iterations $k = n+1, n+2, \cdots, 2n$ the
second cycle and so on. Then definition (\ref{equ:607}) increases
the degree of $\phi$ by $1$ at the beginning of a cycle. One actually can define ${\bf r}_k$ by increasing the degree of $\phi$
by $1$ anywhere within an iteration cycle. Each definition will lead to a
different ML($n$)BiCGStab algorithm. As an illustration, we derive a second ML($n$)BiCGStab algorithm associated with the
definition
\begin{equation} \label{equ:607a}
{\bf r}_{jn+i} = \left\{ \begin{array}{lccl} \phi_j({\bf A})
\,\widehat{\bf r}_{jn+i} & & \mbox{if} &1 \leq i \leq n-1 \\
 \phi_{j+1}
({\bf A}) \,\widehat{\bf r}_{jn+i} & &\mbox{if} & i = n.
\end{array} \right.
\end{equation}
(\ref{equ:607a}) increases the degree of $\phi$ by $1$ at the end of
a cycle. The resulting algorithm requires about $25\%$ less storage
(not counting the storage of the coefficient matrix and the preconditioner) than
the algorithm associated with definition (\ref{equ:607}). However,
one drawback with this storage-saving algorithm is that, in some experiments, its computed
residual ${\bf r}_k$ can easily diverge from the corresponding exact
residual when $n$ is moderately large.

Both ML($n$)BiCG and ML($n$)BiCGStab possess a set of left starting vectors (or, shadow vectors) ${\bf q}_1, \cdots, {\bf q}_n$
that can be chosen freely. This freedom appears to be an advantage of the methods. It not only helps stabilize the performance of
the algorithms, but also allows to see a connection between the Lanczos-based BiCG/BiCGStab and the Arnoldi-based FOM.

Just like BiCGStab, ML($n$)BiCGStab can suffer from three types of
breakdown, caused respectively by the failure of the underlying
Lanczos process, the non-existence of the $LU$ factorization during the
construction of ML($n$)BiCG and the parameters $\rho_k$.
We prove that the breakdown probability is zero when the shadow vectors are selected randomly.


The outline of the paper is as follows. In \S\ref{sec:8-1}, we
introduce index functions. In \S\ref{sec:12-24-1}, we present the
ML($n$)BiCG algorithm introduced in \cite{yeungchan}, from which
ML($n$)BiCGStab algorithms are derived.
 In \S\ref{sec:12-24-2}, we rederive the
ML($n$)BiCGStab algorithm in \cite{yeungchan} by index functions.
In \S\ref{sec:12-24-3}, we derive a storage-saving ML($n$)BiCGStab
algorithm from a different definition of the residual vector. In
\S\ref{sec:12-24-4}, we discuss relationships of ML($n$)BiCGStab
with some other methods. In \S\ref{sec:12-19}, implementation
issues are addressed. Concluding remarks
are made in \S\ref{con}.

\section{Index Functions}\label{sec:8-1} Let be given a positive integer $n$.
For all integers $k$, we define
$$
\begin{array}{lcl}
g_n(k) = \lfloor ( k - 1 ) /n \rfloor & \mbox{and}& r_n(k) = k - n
g_n(k)
\end{array}
$$
where $\lfloor \,
\cdot \,\rfloor$ rounds its argument to the nearest integer towards
minus infinity. We call $g_n$ and $r_n$ index functions; they are
defined on ${\mathcal Z}$, the set of all integers, with ranges ${\mathcal Z}$ and $\{1, 2, \cdots,
n\}$, respectively.

If we write
\begin{equation}\label{equ:8-1}
k = j n + i
\end{equation}
with $1 \leq i \leq n$ and $j \in {\mathcal Z}$, then
\begin{equation} \label{equ:7-9-10}
\begin{array}{rcl} g_n (j n + i) = j &\mbox{and}& r_n(j n
+ i) = i.
\end{array}
\end{equation}
Table \ref{fig:7-9-1} illustrates the behavior of $g_n$ and $r_n$
with $n=3$. It can be seen that $g_n(k)$ has a jump when $k$, moved from left to right, passes a multiple of $n$.

\begin{table}
\begin{center}%
\begin{tabular}[t]{c|rlllcl}
 $k$ &0
  & 1 2 3 & 4 5 6 & 7 8 9 &  10 11 12 & $\cdots$ \\ \hline
  $g_n(k)$ &-1
  & 0 0 0 & 1 1 1 & 2 2 2 & 3 3 3 & $\cdots$ \\
  $r_n(k)$&3
  & 1 2 3 & 1 2 3 & 1 2 3 & 1 2 3 & $\cdots$ \\[1.0ex]
\end{tabular}%
\caption{Simple illustration of the index functions for $n = 3$.}
\label{fig:7-9-1}
\end{center}%
\end{table}

The following properties
can be easily verified by using (\ref{equ:7-9-10}).
\\

\begin{proposition} \label{prop:1} Let $k \in {\mathcal N}$, the set of all positive integers, and $s \in
{\mathcal N}_0 := {\mathcal N} \cup \{ 0 \}$.
\begin{enumerate}
\item[(a)] $g_n(k+n) = g_n(k) + 1$ and $r_n(k + n) = r_n(k)$.

\item[(b)]  $g_n(s+1)+1 =
g_n(k + 1)$ if $\max (k - n, 0) \leq s \leq g_n(k) n - 1$.

\item[(c)] $g_n(s+1) = g_n( g_n(k)n + 1)
= g_n(k)$ if $g_n(k) n \leq s \leq k - 1$.

\item[(d)]
$g_n(k+1) = g_n(k)+1$ if $r_n(k) = n$.\\
 $g_n(k+1) = g_n(k)$ if $r_n(k) < n$.

\item[(e)] $\max(k -n, 0) > g_n(k) n -1$ if
$r_n(k) = n$ or $g_n(k) = 0$.
\end{enumerate}
\end{proposition}

\vspace{.2cm}

\section{A ML($n$)BiCG Algorithm} \label{sec:12-24-1} Parallel to the derivation of
BiCGStab from BiCG, ML($n$)BiCGStab was derived in \cite{yeungchan} from a BiCG-like
method named ML($n$)BiCG, which was built
upon a band Lanczos process with $n$ left starting vectors and a
single right starting vector. In this section, we present the
algorithm of ML($n$)BiCG from \cite{yeungchan} and summarize some
properties of it.

\subsection{The Algorithm}
Consider the solution of the linear system
\begin{equation} \label{equ:7-9-3}
{\bf A} {\bf x} = {\bf b}
\end{equation}
where ${\bf A} \in {\mathcal C}^{N \times N}$ and ${\bf b} \in {\mathcal
C}^N$.
Throughout the
paper we do not assume the matrix ${\bf A}$ is nonsingular except
where specified. The solution of singular systems has been extensively studied in area of iterative methods, see, for instance, \cite{berpl, kaas, masz, reye, wewu, zhangwei}.

Let be given $n$ vectors ${\bf q}_1, \ldots, {\bf
q}_n \in {\mathcal C}^N$, which we call {\it left starting vectors} or {\it shadow vectors}. Set
\begin{equation}
 {\bf p}_{k} = \left( {\bf A}^H \right)^{g_n(k)}
{\bf q}_{r_n(k)}
\label{equ:7-9-5}
\end{equation}
for $k = 1, 2, \cdots$. The following algorithm for solving
(\ref{equ:7-9-3}) is from
\cite{yeungchan}.\\

\begin{algorithm}{\bf ML($n$)BiCG}
\footnote{Algorithm \ref{alg:1} consists of exact mathematical
formulas for $\alpha_k, \beta_s^{(k)}, \widehat{\bf x}_k,
\widehat{\bf r}_k$ and $\widehat{\bf g}_k$ obtained in \S3 of
\cite{yeungchan}. Repeated operations should be removed in order to
make the algorithm computationally efficient. Moreover, even though
the algorithm has not been
tested,
it is believed to be numerically instable because of Line 11 in
which the left starting vectors are repeatedly multiplied by ${\bf
A}^H$, a type of operation which is highly sensitive to round-off
errors. The algorithm has been introduced only for the purpose of
developing ML($n$)BiCGStab algorithms.
}
\label{alg:1} \vspace{.2cm}
\begin{tabbing}
x\=xxx\= xxx\=xxx\=xxx\=xxx\=xxx\kill \>1. \> Choose an initial
guess $\widehat{\bf x}_0$ and $n$ vectors ${\bf q}_1, {\bf q}_2,
\cdots,
{\bf q}_n$. \\
\>2. \>  Compute $\widehat{\bf r}_0 = {\bf b} - {\bf A} \widehat{\bf
x}_0$ and
set ${\bf p}_1 = {\bf q}_1$, $\widehat{\bf g}_0 = \widehat{\bf r}_0$. \\
\>3. \>For $k = 1, 2, \cdots$, until convergence: \\
\>4. \>\>$\alpha_k = {\bf p}_k^H \widehat{\bf r}_{k-1} / {\bf p}_k^H
{\bf A} \widehat{\bf g}_{k-1}$; \\
\>5. \>\>$\widehat{\bf x}_k = \widehat{\bf x}_{k-1} + \alpha_k \widehat{\bf g}_{k-1}$; \\
\>6. \>\> $\widehat{\bf r}_k = \widehat{\bf r}_{k-1} - \alpha_k {\bf A} \widehat{\bf g}_{k-1}$; \\
\>7. \>\>For $s = \max (k - n, 0), \cdots, k - 1$ \\
\>8. \>\>\>$\beta^{(k)}_{s} = - {\bf p}^H_{s+1} {\bf A}
 \left(\widehat{\bf r}_k + \sum_{t = \max (k - n, 0) }^{s-1} \beta^{(k)}_t \widehat{\bf g}_t \right)
\big/{\bf p}^H_{s+1} {\bf A} \widehat{\bf g}_s$; \\
\>9. \>\>End \\
\>10.\>\> $\widehat{\bf g}_k = \widehat{\bf r}_k + \sum_{s = \max (k
- n, 0) }^{k -1} \beta_{s}^{(k)}
            \widehat{\bf g}_{s}$; \\
\>11. \>\>Compute ${\bf p}_{k+1}$ according to (\ref{equ:7-9-5}) \\
\>12. \> End
\end{tabbing}
\end{algorithm}
\vspace{.2cm}

The ML($n$)BiCG algorithm is a variation of the classical BiCG
algorithm. The left-hand side (shadow) Krylov subspace of BiCG is
replaced by the block Krylov subspace with $n$ starting vectors
${\bf q}_1, {\bf q}_2, \cdots, {\bf q}_n$:
$$\begin{array}{rl}
{\mathcal B}_k
&:= \mbox{the space spanned by the first } k \mbox{ columns of } [{\bf Q}, {\bf A}^H {\bf Q}, ({\bf A}^H)^2 {\bf Q}, \cdots ]\\
& = span \{{\bf p}_1, {\bf p}_2,
\cdots, {\bf p}_k\}\\
& = \sum_{i = 1}^{r_n(k)}{\mathcal K}_{g_n(k)+1}({\bf A}^H, {\bf q}_i) +
\sum_{i = r_n(k)+1}^n {\mathcal K}_{g_n(k)}({\bf A}^H, {\bf q}_i)
\end{array}
$$
where ${\bf Q} := [{\bf q}_1, {\bf q}_2, \cdots, {\bf q}_n]$ and
$$
{\mathcal K}_t({\bf M}, {\bf v}) := span \{{\bf v}, {\bf M} {\bf v},
\cdots, {\bf M}^{t-1} {\bf v}\}
$$
for ${\bf M} \in {\mathcal C}^{N \times N}, {\bf v} \in {\mathcal C}^N$ and
$t \in {\mathcal N}$. Moreover, in ML($n$)BiCG, the basis used for
${\mathcal B}_k$ is not chosen to be bi-orthogonal, but simply the set
$\{{\bf p}_1, {\bf p}_2, \cdots, {\bf p}_k\}$. Therefore, the
ML($n$)BiCG algorithm can be viewed as a generalization of a
one-sided Lanczos algorithm (see \cite{gut10, saad82}). The likely
ill-conditioning of this basis does not matter, as the algorithm is
only a technical tool for deriving ML($n$)BiCGStab and this basis
disappears in ML($n$)BiCGStab because ${\bf A}^H$ will be absorbed
by the residuals and direction vectors of ML($n$)BiCGStab. For constructing the
right-hand side basis consisting of residuals $\widehat{\bf r}_k$,
we used recurrences that generalize the coupled two-term recurrences
of BiCG, that is, direction vectors $\widehat{\bf g}_k$ are also
constructed.

\subsection{Properties}\label{sec:5-18}
Let $\nu$ be the degree of the minimal polynomial $p_{min}(\lambda;
{\bf A}, \widehat{\bf r}_0)$ of $\widehat{\bf r}_0$ with respect to
${\bf A}$ (that is, the unique monic polynomial $p(\lambda)$ of
minimum degree such that $p({\bf A})\widehat{\bf r}_0 = {\bf 0}$)
and let
$$
\widehat{\bf S}_\nu = [{\bf p}_1, {\bf p}_2, \cdots, {\bf p}_\nu]^H {\bf A}
[\widehat{\bf r}_0, {\bf A} \widehat{\bf r}_0, \cdots,
{\bf A}^{\nu-1} \widehat{\bf r}_0]
$$
and
$$
\widehat{\bf W}_\nu = [{\bf p}_1, {\bf p}_2, \cdots, {\bf p}_\nu]^H
[\widehat{\bf r}_0, {\bf A} \widehat{\bf r}_0, \cdots, {\bf
A}^{\nu-1} \widehat{\bf r}_0].
$$
Denote by $\widehat{\bf S}_l$ and $\widehat{\bf W}_l$ the $l \times
l$ leading principal submatrices of $ \widehat{\bf S}_\nu$ and
$\widehat{\bf W}_\nu$ respectively. We now summarize some useful
facts about Algorithm \ref{alg:1}. They can be derived from the
construction procedure of the algorithm.
\\

\begin{proposition}\label{prop:7-11-1}
In infinite precision arithmetic, if $
\prod_{l =
1}^\nu \det(\widehat{\bf S}_l) \det(\widehat{\bf W}_l) \ne 0
$, then Algorithm \ref{alg:1} does not break down by zero division
for $k = 1, 2, \cdots, \nu$, and $x_\nu$ is the exact solution of
(\ref{equ:7-9-3}). Moreover, the computed quantities satisfy
\begin{enumerate}
\item[(a)]
$\widehat{\bf x}_k \in \widehat{\bf x}_0 + {\mathcal K}_k ({\bf A},
\widehat{\bf r}_0)
$ and $\widehat{\bf r}_k = {\bf b} - {\bf A} \widehat{\bf x}_k \in
\widehat{\bf r}_0 + {\bf A} {\mathcal K}_k ({\bf A}, \widehat{\bf r}_0)
$
for $1 \leq k \leq \nu$.
\item[(b)] $span \{\widehat{\bf r}_0, \widehat{\bf r}_1, \cdots, \widehat{\bf r}_{k-1}\} =
{\mathcal K}_k ({\bf A}, \widehat{\bf r}_0)
$ for $1 \leq k \leq \nu$.
\item[(c)] $span \{{\bf A} \widehat{\bf r}_0, {\bf A}
\widehat{\bf r}_1, \cdots, {\bf A} \widehat{\bf r}_{\nu-1}\} = {\mathcal
K}_\nu ({\bf A}, \widehat{\bf r}_0)
$.
\item[(d)] $
\widehat{\bf r}_k \perp {\mathcal B}_k
$
and $\widehat{\bf
r}_k \not\perp {\bf p}_{k+1}$ for $0 \leq k \leq \nu
-1$.\footnote{We say that ${\bf u} \perp {\bf v}$ if ${\bf u}^H {\bf
v} = 0$.}
\item[(e)] $span \{\widehat{\bf g}_0, \widehat{\bf g}_1, \cdots, \widehat{\bf g}_{k-1}\} =
{\mathcal K}_k ({\bf A}, \widehat{\bf r}_0)
$ for $1 \leq k \leq \nu$.

\item[(f)] $span \{{\bf A} \widehat{\bf g}_0, {\bf A}
\widehat{\bf g}_1, \cdots, {\bf A} \widehat{\bf g}_{\nu-1}\} = {\mathcal
K}_\nu ({\bf A}, \widehat{\bf r}_0)
$.

\item[(g)] $
{\bf A} \widehat{\bf g}_k \perp {\mathcal B}_k
$ and ${\bf A}
\widehat{\bf g}_k \not\perp {\bf p}_{k+1}$ for $0 \leq k \leq \nu
-1$.
\end{enumerate}
\end{proposition}
\vspace{.2cm}

Because of Proposition \ref{prop:7-11-1}(a) and (d), ML($n$)BiCG is
an oblique projection Krylov subspace method\cite{saad}.

{\it Remarks}:
\begin{enumerate}
\item[(i)] The matrices $\widehat{\bf S}_l$ and $\widehat{\bf W}_l$
have already appeared in \cite{joubert1, joubert2} where they were
called moment matrices. Proposition \ref{prop:7-11-1} can be
regarded as a generalization of Theorem 2 in \cite{joubert2} from $n
= 1$ to $n
> 1$.

\item[(ii)]
Just like BiCG, ML($n$)BiCG also has two types of breakdown caused,
respectively, by the failure of the underlying Lanczos process and
the nonexistence of the $LU$ factorizations of the Hessenberg matrix
of the recurrence coefficients. Both types of breakdown are
reflected in Algorithm \ref{alg:1} by
${\bf p}_k^H
{\bf A} \widehat{\bf g}_{k-1} = 0$.
The condition $\prod_{l = 1}^\nu \det(\widehat{\bf W}_l) \ne 0$
guarantees that the underlying Lanczos process
works without
breakdown, and the condition $\prod_{l = 1}^\nu \det(\widehat{\bf
S}_l) \ne 0$ ensures that
the $LU$ factorizations 
 exist.

\item[(iii)] $\det (\widehat{\bf S}_\nu) \ne 0$ implies that $p_{min}(0;
{\bf A}, \widehat{\bf r}_0) \ne 0$ which, in turn, implies that
(\ref{equ:7-9-3}) is consistent and has a solution lying in
$\widehat{\bf x}_0 + {\mathcal K}_\nu ({\bf A}, \widehat{\bf r}_0)
$.
\end{enumerate}
\vspace{.2cm}

The derivation of ML($n$)BiCGStab will require the following result which, in the case when $n = 1$, has been used in CGS and BiCGStab. \\

\begin{corollary}\label{cor:7-23}
Let $s \in {\mathcal N}$ and
$$
\psi_{g_n(s)}(\lambda) = c_{g_n(s)}\lambda^{g_n(s)} +
c_{g_n(s)-1}\lambda^{g_n(s)-1} + \cdots + c_0
$$
be any polynomial of exact degree $g_n(s)$. Then, under the
assumptions of Proposition \ref{prop:7-11-1},
$$
\begin{array}{rcl}
\displaystyle{{\bf p}_s^H \widehat{\bf r}_k = \frac{1}{c_{g_n(s)}}
{\bf q}_{r_n(s)}^H \psi_{g_n(s)}({\bf A}) \widehat{\bf r}_k }&
\mbox{and} & \displaystyle{{\bf p}_s^H {\bf A} \widehat{\bf g}_k =
\frac{1}{c_{g_n(s)}} {\bf q}_{r_n(s)}^H {\bf A} \psi_{g_n(s)}({\bf
A}) \widehat{\bf g}_k }
\end{array}
$$
if $0 \leq k \leq \nu-1$ and $s \leq k+n$.
\end{corollary}
\vspace{.2cm}

\proof It is easy to verify that
$$
{\bf p}_s - \frac{1}{\bar{c}_{g_n(s)}} \bar{\psi}_{g_n(s)}({\bf
A}^H) {\bf q}_{r_n(s)} \in {\mathcal B}_k
$$
by Proposition \ref{prop:1}(a) and (\ref{equ:7-9-5}), where the
overbar denotes complex conjugation. The corollary then follows from
Proposition \ref{prop:7-11-1}(d) and (g).
\hfill{\rule{2mm}{2mm}} \vspace{.1cm} \vspace{.1cm}

Corollary \ref{cor:7-23} essentially says that adding to ${\bf p}_s$
a vector from ${\mathcal B}_k$ does not change the inner products ${\bf
p}_s^H \widehat{\bf r}_k$ and ${\bf p}_s^H {\bf A} \widehat{\bf
g}_k$.

Examples exist where the condition $
\prod_{l =
1}^\nu \det (\widehat{\bf W}_l) \det( \widehat{\bf S}_l) \ne 0
$ in Proposition \ref{prop:7-11-1} holds, as
shown below.\\

\begin{lemma} \label{thm:12-15} Consider the case where $n = 1,
\,\widehat{\bf r}_0 \in {\mathcal R}^N, \,\widehat{\bf r}_0 \ne {\bf
0}$\footnote{In Lemma \ref{thm:12-15} and Theorem \ref{thm:3-11-09},
$\widehat{\bf r}_0$ can be any non-zero vector in ${\mathcal R}^N$, not
necessary to be a residual vector like $\widehat{\bf r}_0 = {\bf b}
- {\bf A} \widehat{\bf x}_0$.} and ${\bf A} \in {\mathcal R}^{N \times
N}$ is nonsingular. If ${\bf q}_1 \in {\mathcal R}^N$ is a random vector
with independent and identically distributed elements from $N(0,
1)$, the normal
distribution with mean $0$ and variance $1$, then $
Prob \left( \prod_{l = 1}^\nu \det (\widehat{\bf W}_l) \det(
\widehat{\bf S}_l) = 0 \right) = 0
$.
\end{lemma}

\proof Since ${\bf p}_k = {\bf A}^{g_n (k)} {\bf q}_{r_n(k)} = {\bf
A}^{k-1} {\bf q}_1$ when $n = 1$, both $\widehat{\bf S}_\nu$ and
$\widehat{\bf W}_\nu$ are Hankel matrices
$$\begin{array}{rcl}
 \widehat{\bf S}_\nu = \left[ \begin{array}{llll}
 \widehat{s}_1 & \widehat{s}_2& \cdots & \widehat{s}_\nu\\
\widehat{s}_2 & \widehat{s}_3& \cdots & \widehat{s}_{\nu + 1}\\
\cdots & \cdots & \cdots & \cdots\\
\widehat{s}_\nu & \widehat{s}_{\nu+1}& \cdots & \widehat{s}_{2\nu-1}
\end{array} \right], & &
\widehat{\bf W}_\nu = \left[ \begin{array}{llll} \widehat{w}_1 &
\widehat{w}_2& \cdots & \widehat{w}_\nu\\
\widehat{w}_2 & \widehat{w}_3& \cdots & \widehat{w}_{\nu + 1}\\
\cdots & \cdots & \cdots & \cdots\\
\widehat{w}_\nu & \widehat{w}_{\nu+1}& \cdots & \widehat{w}_{2\nu-1}
\end{array} \right]
\end{array}
$$
where $\widehat{s}_t = {\bf q}_1^T {\bf A}^t \widehat{\bf r}_0$ and
$\widehat{w}_t = {\bf q}_1^T {\bf A}^{t-1} \widehat{\bf r}_0$ for $t
= 1, 2, \cdots, 2 \nu -1$.

We first prove
\begin{equation} \label{eq:6-28-1}
Prob \left( \det(\widehat{\bf W}_l) = 0 \right) = 0
\end{equation}
 for any fixed $l$ with
$1 \leq l \leq \nu$. It is trivial that (\ref{eq:6-28-1}) holds when
$l = 1$ and we therefore assume $l \geq 2$ in the following
discussion.

By assumption, $\nu$ is the degree of the minimal polynomial of
$\widehat{\bf r}_0$ with respect to ${\bf A}$. This implies that
${\mathcal K} := span \{{\bf A}^t \widehat{\bf r}_0\, |\, t \in {\mathcal
N}_0\}$ is a $\nu$-dimensional space with $\{\widehat{\bf r}_0, {\bf
A} \widehat{\bf r}_0, \cdots, {\bf A}^{\nu - 1} \widehat{\bf r}_0\}$
as a basis. Since ${\bf A}$ is nonsingular, $\{{\bf
A}^{l-1}\widehat{\bf r}_0, {\bf A}^l \widehat{\bf r}_0, \cdots, {\bf
A}^{l+\nu - 2} \widehat{\bf r}_0\}$ is another basis of ${\mathcal K}$.

Perform an orthogonal factorization of the $N \times \nu$ matrix
$$
\left[ {\bf A}^{l-1}\widehat{\bf r}_0, {\bf A}^l \widehat{\bf r}_0,
\cdots, {\bf A}^{l+\nu - 2} \widehat{\bf r}_0 \right] = {\bf Q} {\bf
R}
$$
where ${\bf Q} \in {\mathcal R}^{N \times N}$ is orthogonal and ${\bf R}
\in {\mathcal R}^{N \times \nu}$ is upper triangular with positive main
diagonal elements $r_{11}, r_{22}, \cdots, r_{\nu \nu}$. Clearly,
the first $\nu$ columns of ${\bf Q}$ form a basis of ${\mathcal K}$ and
the last $N - \nu$ columns belong to ${\mathcal K}^\perp$.

Write
$$\begin{array}{rl}
{\bf A}^{l-2} \widehat{\bf r}_0 &= \xi_1 {\bf A}^{l-1}\widehat{\bf
r}_0 + \xi_2 {\bf
A}^l \widehat{\bf r}_0 + \cdots + \xi_\nu {\bf A}^{l+\nu - 2} \widehat{\bf r}_0\\
& = \left[ {\bf A}^{l-1}\widehat{\bf r}_0, {\bf A}^l \widehat{\bf
r}_0, \cdots, {\bf
A}^{l+\nu - 2} \widehat{\bf r}_0 \right] {\bf \xi}\\
& = {\bf Q} {\bf R} {\bf \xi} \equiv {\bf Q} {\bf \eta}
\end{array}
$$
for some scalars $\xi_1, \xi_2, \cdots, \xi_\nu \in {\mathcal R}$, where
${\bf \xi} = [\xi_1, \xi_2, \cdots, \xi_\nu]^T \in {\mathcal R}^{\nu}$
and ${\bf \eta} = [\eta_1, \eta_2, \cdots, \eta_N]^T = {\bf R} {\bf
\xi} \in {\mathcal R}^N$. Since ${\bf A}$ is nonsingular and
$\{\widehat{\bf r}_0, {\bf A} \widehat{\bf r}_0, \cdots, {\bf
A}^{\nu - 1} \widehat{\bf r}_0\}$ linearly independent, we have
$\xi_\nu \ne 0$ and hence $\eta_\nu = r_{\nu \nu} \xi_\nu \ne 0$.
Let ${\bf \theta} = [\theta_1, \theta_2, \cdots, \theta_N]^T = {\bf
Q}^T {\bf q}_1$. Then ${\bf \theta}$ is a random vector with iid
elements from $N(0, 1)$\cite{edelman}. We now express
$\det(\widehat{\bf W}_l)$ in terms of the elements of ${\bf
\theta}$. Let us write
$$\begin{array}{rl}
&[\widehat{w}_1, \widehat{w}_2, \cdots, \widehat{w}_{l-2},
\widehat{w}_{l-1}, \widehat{w}_l, \cdots, \widehat{w}_{2l-1}]\\
 =&
{\bf q}_1^T [ \widehat{\bf r}_0, {\bf A} \widehat{\bf r}_0, \cdots,
{\bf A}^{l-3} \widehat{\bf r}_0, {\bf A}^{l-2} \widehat{\bf r}_0,
{\bf A}^{l-1} \widehat{\bf r}_0, \cdots,
{\bf A}^{2l-2} \widehat{\bf r}_0]\\
=& {\bf q}_1^T [ \widehat{\bf r}_0, {\bf A} \widehat{\bf r}_0,
\cdots, {\bf A}^{l-3}
\widehat{\bf r}_0, {\bf Q} {\bf \eta}, {\bf Q} {\bf R}^{(l)}]\\
=& {\bf q}_1^T {\bf Q} [{\bf Q}^T [ \widehat{\bf r}_0, {\bf A}
\widehat{\bf r}_0,
\cdots, {\bf A}^{l-3} \widehat{\bf r}_0], {\bf \eta}, {\bf R}^{(l)}]\\
=& {\bf \theta}^T  [{\bf Q}^T [ \widehat{\bf r}_0, {\bf A}
\widehat{\bf r}_0, \cdots, {\bf A}^{l-3} \widehat{\bf r}_0], {\bf
\eta}, {\bf R}^{(l)}]
\end{array}
$$
where ${\bf R}^{(l)}$ denotes the matrix consisting of the first $l$
columns of ${\bf R}$. Since the last $N - \nu$ columns of $\bf Q$
belong to ${\mathcal K}^\perp$, the last $N - \nu$ rows of the matrix
${\bf Q}^T [ \widehat{\bf r}_0, {\bf A} \widehat{\bf r}_0, \cdots,
{\bf A}^{l-3} \widehat{\bf r}_0]$ are zeros. Similarly, the last
$N-\nu$ elements of ${\bf \eta} = {\bf R} {\bf \xi}$ are zeros
because the last $N - \nu$ rows of $\bf R$ are zeros. We therefore
have
$$
\widehat{w}_t = \left\{ \begin{array}{l}
      \mbox{a linear combination of } \theta_1, \theta_2, \cdots,
      \theta_\nu \mbox{ if } 1 \leq t \leq l-2,\\
      \eta_1 \theta_1 + \eta_2 \theta_2 + \cdots + \eta_\nu
      \theta_\nu \mbox{ with } \eta_\nu \ne 0 \mbox{ if } t = l -1,\\
      \begin{array}{r}
r_{1, t-l+1} \theta_1 + r_{2, t-l+1} \theta_2 + \cdots + r_{t-l+1,
t-l+1} \theta_{t-l+1}\\  \mbox{ with } r_{t-l+1, t-l+1} \ne 0
 \mbox{
if } l \leq t \leq 2l-1.
\end{array}
      \end{array}\right.
$$
This shows that none of the random variables $\theta_{\nu +1},
\theta_{\nu +2}, \cdots, \theta_N$ is involved in any of the
$\widehat{w}$'s. In more detail, when $l < \nu$, $\widehat{w}_t =
\widehat{w}_t(\theta_1, \theta_2, \cdots, \theta_\nu)$ if $1 \leq t
\leq l-1$ and $\widehat{w}_t = \widehat{w}_t(\theta_1, \theta_2,
\cdots, \theta_{\nu-1})$ if $l \leq t \leq 2l-1$; when $l = \nu$,
$\widehat{w}_t = \widehat{w}_t(\theta_1, \theta_2, \cdots,
\theta_\nu)$ if $1 \leq t \leq \nu-1$ or $t = 2 \nu -1$, and
$\widehat{w}_t = \widehat{w}_t(\theta_1, \theta_2, \cdots,
\theta_{\nu-1})$ if $\nu \leq t < 2\nu-1$.

We now expand $\det (\widehat{\bf W}_l)$ by minors down its last
column and write it into a polynomial in $\theta_\nu$. This yields
$$
\begin{array}{rl}
&(-1)^{\frac{1}{2}l(l+1)+1} \det (\widehat{\bf W}_l)\\
 =& \widehat{w}_{2l-1} \widehat{w}_{l-1}^{l-1} + \cdots \\
= &\left\{ \begin{array}{ll} 
 (\sum_{s = 1}^l r_{sl} \theta_s)
\eta_\nu^{l-1} \theta_\nu^{l-1} + c_{l-2} \theta_\nu^{l-2} + \cdots
+ c_1 \theta_\nu + c_0 &\mbox{if } 2 \leq l < \nu,
\\
r_{\nu \nu} \eta_{\nu}^{\nu -1} \theta_\nu^{\nu} + d_{\nu -1}
\theta_\nu^{\nu-1} + \cdots + d_1 \theta_\nu + d_0 &
 \mbox{if } l = \nu
      \end{array}\right.
\end{array}
$$
where the coefficients $c_0, \cdots, c_{l-2}$ and $d_0, \cdots,
d_{\nu -1}$ are polynomials in $\theta_1, \theta_2, \cdots,
\theta_{\nu -1}$. Now (\ref{eq:6-28-1}) follows from the facts that
$r_{ll} \ne 0, r_{\nu \nu} \ne 0, \eta_\nu \ne 0$ and $\theta_1,
\theta_2, \cdots, \theta_\nu$ are independent random variables.

Note that $\nu$ is also the degree of the minimal polynomial of
${\bf A} \widehat{\bf r}_0$ with respect to ${\bf A}$ when ${\bf A}$
is nonsingular. With $\widehat{\bf r}_0$ replaced by ${\bf A}
\widehat{\bf r}_0$ in (\ref{eq:6-28-1}), we then have
\begin{equation} \label{eq:7-7-1}
Prob \left( \det( \widehat{\bf S}_l) = 0 \right) = 0
\end{equation}
for any $l$ with $1 \leq l \leq \nu$.

Now, (\ref{eq:6-28-1}) and (\ref{eq:7-7-1}) together imply that
$$\begin{array}{l}
Prob \left( \prod_{l = 1}^\nu \det (\widehat{\bf W}_l) \det(
\widehat{\bf S}_l) = 0 \right)\\
\leq \sum_{l = 1}^\nu Prob \left( \det (\widehat{\bf W}_l) = 0
\right) + \sum_{l = 1}^\nu Prob \left( \det ( \widehat{\bf S}_l) = 0
\right)
 = 0
\end{array}
$$
and the lemma is proved. \hfill{\rule{2mm}{2mm}}\\

The $\bf A$ in Lemma \ref{thm:12-15} is assumed to be nonsingular. For a
general
$\bf A$, we have\\

\begin{theorem}\label{thm:3-11-09} Consider the case where $n = 1,
\,\widehat{\bf r}_0 \in {\mathcal R}^N, \,\widehat{\bf r}_0 \ne {\bf 0}$
and ${\bf A} \in {\mathcal R}^{N \times N}$. If ${\bf q}_1 \in {\mathcal
R}^N$ is a random vector with independent and identically
distributed elements from $N(0, 1)$, then
 $
 Prob \left( \prod_{l = 1}^\nu
\det (\widehat{\bf W}_l) \det( \widehat{\bf S}_l) = 0 \right) = 0
$ if and only if $p_{min}(0; {\bf A}, \widehat{\bf r}_0) \ne 0$.
\end{theorem}
\vspace{.2cm}

\proof If $p_{min}(0; {\bf A}, \widehat{\bf r}_0) = 0$, then ${\bf
A}^\nu \widehat{\bf r}_0$ is a linear combination of ${\bf A}
\widehat{\bf r}_0, \cdots, {\bf A}^{\nu -1} \widehat{\bf r}_0$ or
${\bf A}^\nu \widehat{\bf r}_0 = {\bf 0}$ in the case when $\nu =
1$. Hence $\det (\widehat{\bf S}_\nu) = 0$ no matter what ${\bf
q}_1$ is and therefore $Prob \left(\prod_{l = 1}^\nu \det
(\widehat{\bf W}_l) \det( \widehat{\bf S}_l) = 0 \right) = 1$.

We now suppose $p_{min}(0; {\bf A}, \widehat{\bf r}_0) \ne 0$. By
the real version of the Schur's unitary triangularization theorem (see, for instance, 
\cite{horn}), ${\bf A}$ can be decomposed as
$$
{\bf A} = {\bf Q}^T \left[ \begin{array}{cc}
{\bf B}_{11} & {\bf B}_{12}\\
{\bf 0} & {\bf B}_{22} \end{array} \right] {\bf Q} \equiv {\bf Q}^T
{\bf B} {\bf Q}
$$
where ${\bf Q} \in {\mathcal R}^{N \times N}$ is orthogonal, ${\bf B}_{11}
\in {\mathcal R}^{N_1 \times N_1}$ nonsingular and ${\bf B}_{22} \in
{\mathcal R}^{N_2 \times N_2}$ strictly upper triangular (namely, an
upper triangular matrix with its main diagonal elements zero). Let
$\widetilde{\bf r}_0 = {\bf Q}\, \widehat{\bf r}_0 \equiv
[\widetilde{\bf r}_{01}^T, \widetilde{\bf r}_{02}^T]^T$ where
$\widetilde{\bf r}_{01} \in {\mathcal R}^{N_1}$ and $\widetilde{\bf
r}_{02} \in {\mathcal R}^{N_2}$. Then $p_{min}({\bf B}; {\bf A},
\widehat{\bf r}_0) \widetilde{\bf r}_0 = {\bf Q} \,p_{min}({\bf A};
{\bf A}, \widehat{\bf r}_0) \widehat{\bf r}_0 = {\bf 0}$. Note that
\begin{equation}\label{equ:12-15}
{\bf B}^k = \left[ \begin{array}{cc} {\bf B}_{11}^k & *\\
{\bf 0} & {\bf B}_{22}^k
\end{array} \right]
\end{equation}
for $k \in {\mathcal N}$, we have
$$
p_{min}({\bf B}; {\bf A}, \widehat{\bf r}_0) = \left[ \begin{array}{cc} p_{min}({\bf B}_{11}; {\bf A}, \widehat{\bf r}_0) & *\\
{\bf 0} & p_{min}({\bf B}_{22}; {\bf A}, \widehat{\bf r}_0)
\end{array} \right].
$$
Thus, $p_{min}({\bf B}; {\bf A}, \widehat{\bf r}_0) \widetilde{\bf
r}_0 = {\bf 0}$ implies that $p_{min}({\bf B}_{22}; {\bf A},
\widehat{\bf r}_0) \widetilde{\bf r}_{02} =  {\bf 0}$. If we write
$p_{min}(\lambda; {\bf A}, \widehat{\bf r}_0) = \sum_{t = 0}^{\nu}
c_t \lambda^t$ with $c_0 \ne 0$, then $\sum_{t = 1}^\nu c_t {\bf
B}_{22}^t$ is a strictly upper triangular matrix. Thus,
$p_{min}({\bf B}_{22}; {\bf A}, \widehat{\bf r}_0) = (\sum_{t =
1}^\nu c_t {\bf B}_{22}^t) + c_0 {\bf I}$ is an upper triangular
matrix whose main diagonal elements are $c_0$. So, $p_{min}({\bf
B}_{22}; {\bf A}, \widehat{\bf r}_0)$ is nonsingular and therefore
$p_{min}({\bf B}_{22}; {\bf A}, \widehat{\bf r}_0) \widetilde{\bf
r}_{02} =  {\bf 0}$ yields $\widetilde{\bf r}_{02} = {\bf 0}$. Since
$\widetilde{\bf r}_0 \ne {\bf 0}$ due to $\widehat{\bf r}_0 \ne {\bf
0}$ by assumption, $\widetilde{\bf r}_{02} \ne \widetilde{\bf
r}_{0}$. In other words, $N_2 < N$ or ${\bf B}_{11}$ is not a null
matrix.

Now that $\widetilde{\bf r}_{02} = {\bf 0}$, (\ref{equ:12-15})
implies that
\begin{equation}\label{equ:12-15-1}
{\bf B}^k \widetilde{\bf r}_{0} = \left[\begin{array}{c} {\bf
B}_{11}^k \widetilde{\bf r}_{01}\\
 {\bf 0}
 \end{array} \right]
\end{equation}
for $k \in {\mathcal N}$. Therefore, $p({\bf B}) \widetilde{\bf r}_{0} =
[(p({\bf B}_{11}) \widetilde{\bf r}_{01})^T, {\bf 0}^T]^T$ for any
polynomial $p(\lambda)$. Thus, the minimal polynomial of
$\widetilde{\bf r}_{0}$ with respect to ${\bf B}$ is equal to the
minimal polynomial of $\widetilde{\bf r}_{01}$ with respect to ${\bf
B}_{11}$. This implies that, $\nu$, the degree of the minimal
polynomial of $\widehat{\bf r}_0$ with respect to $\bf A$, is also
the degree of the minimal polynomial of $\widetilde{\bf r}_{01}$
with respect to ${\bf B}_{11}$.

We now set ${\bf \theta} = {\bf Q}\, {\bf q}_1 \equiv [ {\bf
\theta}_1^T, {\bf \theta}_2^T]^T$ where ${\bf \theta}_1 \in {\mathcal
R}^{N_1}$ and ${\bf \theta}_2 \in {\mathcal R}^{N_2}$. Since ${\bf q}_1$
is random with iid elements from $N(0,1)$, so is $\bf \theta$. By
(\ref{equ:12-15-1}),
$$\begin{array}{rcl}
{\bf q}_1^T {\bf A}^k \widehat{\bf r}_0 = {\bf \theta}^T {\bf B}^k
\widetilde{\bf r}_0 = {\bf \theta}^T_1 {\bf B}_{11}^k \widetilde{\bf
r}_{01} & \mbox{and} & {\bf q}_1^T  \widehat{\bf r}_0 = {\bf
\theta}^T_1 \widetilde{\bf r}_{01}
\end{array}
$$
where $k \in {\mathcal N}$. Thus
$$
\begin{array}{rcl}
\widehat{\bf S}_\nu ({\bf A}, \widehat{\bf r}_0, {\bf q}_1) =
\widehat{\bf S}_\nu ({\bf B}_{11}, \widetilde{\bf r}_{01}, {\bf
\theta}_1) &\mbox{and} & \widehat{\bf W}_\nu ({\bf A}, \widehat{\bf
r}_0, {\bf q}_1) = \widehat{\bf W}_\nu ({\bf B}_{11}, \widetilde{\bf
r}_{01}, {\bf \theta}_1).
\end{array}
$$
Now, the desired probability follows from Lemma \ref{thm:12-15}
because ${\bf B}_{11}$ is nonsingular, ${\bf \theta}_1$ is iid
$N(0,1)$ random and $\nu$ is the degree of the minimal polynomial of
$\widetilde{\bf r}_{01}$ with respect to ${\bf B}_{11}$. \hfill{\rule{2mm}{2mm}}\\

Extension of the theorem to the general case should be possible,
namely, $n \geq 1, {\bf A} \in {\mathcal C}^{N \times N}, \widehat{\bf
r}_0 \in {\mathcal C}^{N \times 1}$ and $[{\bf q}_1, \cdots, {\bf q}_n]$
is a
Gaussian matrix. We remark that, when $\bf A$ is non-defective, the
general case has been proved
 in the proof of Theorem 3 of
\cite{sonnGij}. The proof was based on the observation that, if a
polynomial $p(\lambda_1, \cdots, \lambda_l) \not\equiv 0$, then
$Prob\left(p(\lambda_1, \cdots, \lambda_l\right) = 0) = 0$ when
$\lambda_1, \cdots, \lambda_l$ are randomly chosen.

{\it Remark}: $p_{min}(0; {\bf A}, \widehat{\bf r}_0) \ne 0$ if and
only if the affine space $\widehat{\bf x}_0 + span \{ {\bf A}^t
\widehat{\bf r}_0 | t \in {\mathcal N}_0\}$ contains a solution to
(\ref{equ:7-9-3}).

 The following corollary then follows from Proposition
\ref{prop:7-11-1} and Theorem \ref{thm:3-11-09}. \\

\begin{corollary}\label{cor:4-1}
In the case where $n = 1$, (\ref{equ:7-9-3}) is a real system and
${\bf q}_1 \in {\mathcal R}^N$ is a random vector with iid elements from
$N(0, 1)$, Algorithm \ref{alg:1} almost surely works without
breakdown by zero division to find a solution from the affine space
$\widehat{\bf x}_0 + span \{ {\bf A}^t \widehat{\bf r}_0 | t \in
{\mathcal N}_0\}$ provided that $\widehat{\bf x}_0 \in {\mathcal R}^{N}$ is
chosen such that the affine space contains a solution to
(\ref{equ:7-9-3}).
\end{corollary}
\vspace{.2cm}


{\it Remarks}: \begin{enumerate} \item[(i)] The initial guess
$\widehat{\bf x}_0$ in Corollary \ref{cor:4-1} is a user-provided
vector. It may not be a random vector in some applications. For
example, in cases where a sequence of similar linear systems is
solved, the solution from the previous system may be used as the
$\widehat{\bf x}_0$ for the new system.

\item[(ii)] If we pick $\widehat{\bf x}_0 \in {\mathcal R}^N$ randomly and set ${\bf
q}_1 = {\bf b} - {\bf A} \widehat{\bf x}_0
$, then
Algorithm \ref{alg:1} with $n = 1$, or equivalently in mathematics,
the standard BiCG (see \S\ref{sec:12-24-4}), almost surely solves
(\ref{equ:7-9-3}) without breakdown by zero division for all, but a
certain small class of, nonsingular ${\bf A} \in {\mathcal R}^{N \times
N}$. For details, see \cite{joubert2}.
\end{enumerate}
\vspace{.1cm}


\section{A ML($n$)BiCGStab Algorithm} \label{sec:12-24-2} An algorithm for the
ML($n$)BiCGStab method has been derived from ML($n$)BiCG in \cite{yeungchan}
(Algorithm 2 without preconditioning and Algorithm 3 with
preconditioning in \cite{yeungchan}), but the derivation there is
complicated and less inspiring. In this section, we re-derive the
algorithm in a more systematic fashion with the help of index
functions.

\subsection{Notation and Definitions}\label{sec:9-12}
Let $\phi_k(\lambda)$ be the polynomial of degree $k$ defined by (\ref{equ:12-23}). If expressed in terms of the
power basis
\begin{equation} \label{equ:5-28-09-10}
\phi_k(\lambda) = c^{(k)}_k \lambda^k + \cdots + c^{(k)}_1 \lambda +
c_0^{(k)},
\end{equation}
it is clear that $c^{(k)}_k = \rho_1 \rho_2 \cdots \rho_k$ and
$c_0^{(k)} = 1$. Thus,
\begin{equation}\label{equ:7-27}
c^{(k)}_k = \rho_k c^{(k-1)}_{k-1}.
\end{equation}

In ML($n$)BiCGStab, we construct the following vectors: for $k \in {\mathcal N}$,
\begin{equation} \label{equ:7-29}
\begin{array}{lcl}
{\bf r}_k = \phi_{g_n(k)+1}({\bf A}) \,\widehat{\bf r}_k,& &{\bf
u}_k = \phi_{g_n(k)} ({\bf A}) \,\widehat{\bf r}_k,\\
{\bf g}_k = \phi_{g_n(k)+1}({\bf A}) \,\widehat{\bf g}_k, & & {\bf
d}_k = \rho_{g_n(k)+1}{\bf A} \phi_{g_n(k)}({\bf A})\, \widehat{\bf
g}_k,\\
{\bf w}_k = {\bf A} {\bf g}_k
\end{array}
\end{equation}
and for $k = 0$, set
\begin{equation}\label{equ:9-13}
\begin{array}{lcl}
{\bf r}_0 = \widehat{\bf r}_0 & \mbox{and} & {\bf g}_0 =
\widehat{\bf g}_0.
\end{array}
\end{equation}
The vectors ${\bf r}_k$ will be the residual vectors of the
approximate solutions ${\bf x}_k$ computed in the following
ML($n$)BiCGStab algorithm.

\subsection{Algorithm Derivation}
The derivation parallels the one of BiCGStab from BiCG. We first
replace all the inner products ${\bf p}^H \widehat{\bf r}$ and ${\bf
p}^H {\bf A} \widehat{\bf g}$ in ML($n$)BiCG respectively by inner
products of the forms ${\bf q}^H \phi({\bf A}) \widehat{\bf r}$ and
${\bf q}^H {\bf A} \phi({\bf A}) \widehat{\bf g}$, where $\phi$ is
the polynomial (\ref{equ:12-23}). Corollary \ref{cor:7-23}
guarantees that the inner products remain unchanged with such
replacements. Then we compile recurrences for the new residuals
${\bf r}_k$ and the corresponding iterates. The overall derivation
is best described and verified in stages, and depends on Proposition
\ref{prop:1} and Corollary \ref{cor:7-23}.

The derivation is complicated by the fact that the recurrences in
the $k$th iteration in ML($n$)BiCG involve $n$ terms which stretch
from $k -n$ to $k - 1$.
Note that $k-n \leq g_n(k) n \leq k-1$. The degrees of the
$\phi_{g_n(s)}$ and $\phi_{g_n(s)+1}$ in (\ref{equ:7-29}) are
increased at $g_n(k) n + 1$ as $s$ runs from $k-n$ to $k-1$ (see,
for example, Table \ref{fig:7-9-1}). Therefore, our first task in
the derivation is to split up in ML($n$)BiCG the loops and the sums
of length $n$ into two parts, one from $k -n$ to $g_n(k)n -1$ and
the other from $g_n(k)n + 1$ to $k-1$. The following Derivation
Stage (DS) \#1 is computationally equivalent to
Algorithm \ref{alg:1} (forgetting Lines 1, 2, 5 and 11).\\

{\bf Derivation Stage \#1.}
\begin{tabbing}
x\=xxx\= xxx\=xxx\=xxx\=xxx\=xxx\kill
\>1. \>For $k = 1, 2, \cdots$, until convergence: \\
\>2.\>\>If $r_n(k) = 1$\\
\>3. \>\>\>$\alpha_k = {\bf p}_k^H \widehat{\bf r}_{k-1} / {\bf
p}_k^H
{\bf A} \widehat{\bf g}_{k-1}$; \\
\>4. \>\>\> $\widehat{\bf r}_k = \widehat{\bf r}_{k-1} - \alpha_k {\bf A} \widehat{\bf g}_{k-1}$; \\
\>5.\>\>Else\\
\>6. \>\>\>$\alpha_k = {\bf p}_k^H \widehat{\bf r}_{k-1} / {\bf
p}_k^H
{\bf A} \widehat{\bf g}_{k-1}$; \\
\>7. \>\>\> $\widehat{\bf r}_k = \widehat{\bf r}_{k-1} - \alpha_k {\bf A} \widehat{\bf g}_{k-1}$; \\
\>8.\>\>End\\
\>9.\>\> If $r_n(k) < n$\\
\>10. \>\>\>For $s = \max (k - n, 0), \cdots, g_n(k)n-1$ \\
\>11. \>\>\>\>$\beta^{(k)}_{s} = - {\bf p}^H_{s+1} {\bf A}
 \left(\widehat{\bf r}_k + \sum_{t = \max (k - n, 0) }^{s-1} \beta^{(k)}_t \widehat{\bf g}_t \right)
\big/{\bf p}^H_{s+1} {\bf A} \widehat{\bf g}_s$; \\
\>12. \>\>\>End \\
\>13. \>\>\>$\beta^{(k)}_{g_n(k)n} = - {\bf p}^H_{g_n(k)n+1} {\bf A}
 \left(\widehat{\bf r}_k + \sum_{t = \max (k - n, 0) }^{g_n(k)n-1} \beta^{(k)}_t \widehat{\bf g}_t
 \right)
\big/{\bf p}^H_{g_n(k)n+1} {\bf A} \widehat{\bf g}_{g_n(k)n}$; \\
\>14. \>\>\>For $s = g_n(k)n+1, \cdots, k-1$ \\
\>15. \>\>\>\>$\beta^{(k)}_{s} = - {\bf p}^H_{s+1} {\bf A}
 \left(\widehat{\bf r}_k + \sum_{t = \max (k - n, 0) }^{g_n(k)n} \beta^{(k)}_t \widehat{\bf g}_t
+ \sum_{t = g_n(k)n+1}^{s-1} \beta^{(k)}_t \widehat{\bf g}_t
 \right)
\big/{\bf p}^H_{s+1} {\bf A} \widehat{\bf g}_s$; \\
\>16. \>\>\>End \\
\>17.\>\>\> $\widehat{\bf g}_k = \widehat{\bf r}_k + \sum_{s = \max
(k - n, 0) }^{g_n(k)n} \beta_{s}^{(k)}
            \widehat{\bf g}_{s} +
\sum_{s = g_n(k)n+1 }^{k -1} \beta_{s}^{(k)}
            \widehat{\bf g}_{s}$; \\
\>18.\>\>Else\\
\>19. \>\>\>$\beta^{(k)}_{g_n(k)n} = - {\bf p}^H_{g_n(k)n+1} {\bf A}
 \widehat{\bf r}_k
\big/{\bf p}^H_{g_n(k)n+1} {\bf A} \widehat{\bf g}_{g_n(k)n}$; \\
\>20. \>\>\>For $s = g_n(k)n+1, \cdots, k-1$ \\
\>21. \>\>\>\>$\beta^{(k)}_{s} = - {\bf p}^H_{s+1} {\bf A}
 \left(\widehat{\bf r}_k +
\beta^{(k)}_{g_n(k)n} \widehat{\bf g}_{g_n(k)n} +
 \sum_{t = g_n(k)n+1}^{s-1} \beta^{(k)}_t \widehat{\bf g}_t
 \right)
\big/{\bf p}^H_{s+1} {\bf A} \widehat{\bf g}_s$; \\
\>22. \>\>\>End \\
\>23.\>\>\> $\widehat{\bf g}_k = \widehat{\bf r}_k +
\beta_{g_n(k)n}^{(k)}
            \widehat{\bf g}_{g_n(k)n} +
 \sum_{s = g_n(k)n+1}^{k -1} \beta_{s}^{(k)}
            \widehat{\bf g}_{s}$; \\
\>24.\>\>End\\
\>25. \> End
\end{tabbing}

\vspace{.2cm} We have adopted the conventions: {\it empty
loops are skipped and empty sums are zero}.
These conventions will also be applied
in the sequel.

In the next stage of the derivation, we replace inner products ${\bf
p}^H \widehat{\bf r}$ and ${\bf p}^H {\bf A} \widehat{\bf g}$ by
inner products of the forms ${\bf q}^H \phi({\bf A}) \widehat{\bf
r}$ and ${\bf q}^H {\bf A} \phi({\bf A}) \widehat{\bf g}$
respectively. That is, the factor $({\bf A}^H)^{g_n(k)}$ that is
hidden in the left basis vector ${\bf p}_k$ is moved to the
right-hand side space and replaced by the factor $\phi_{g_n(k)}({\bf
A})$. Formally, by Corollary \ref{cor:7-23} together with
(\ref{equ:7-9-5}), (\ref{equ:7-27})
 and Proposition
\ref{prop:1}(a), DS\#1 can be further
transformed into the version below. Explanations are given after listing.\\

{\bf Derivation Stage \#2.}
\begin{tabbing}
x\=xxx\=xxx\=xxx\=xxx\=xxx\=xxx\=xxx\=xxx\=xxx\=xxx\=xxx\=xxx\kill
\>1. \>For $k = 1, 2, \cdots$, until convergence: \\
\>2.\>\>If $r_n(k) = 1$\\
\>3. \>\>\>$\displaystyle{\alpha_k = 
{\bf q}_{r_n(k)}^H \phi_{g_n(k)}({\bf A}) \widehat{\bf r}_{k-1}
 /
{\bf q}_{r_n(k)}^H {\bf A}
\phi_{g_n(k)}({\bf A}) \widehat{\bf g}_{k-1}
}$; \\
\>4. \>\>\> $\phi_{g_n(k)}({\bf A}) \widehat{\bf r}_k =
\phi_{g_n(k)}({\bf A}) \widehat{\bf r}_{k-1} - \alpha_k {\bf A}
\phi_{g_n(k)}({\bf A})
\widehat{\bf g}_{k-1}$; \\
\>5. \>\>\> $\phi_{g_n(k)+1}({\bf A}) \widehat{\bf r}_k =
(\rho_{g_n(k)+1} {\bf A} + {\bf I}) \phi_{g_n(k)}({\bf A})
\widehat{\bf r}_k
$; \\
\>6.\>\>Else\\
\>7. \>\>\>$\displaystyle{\alpha_k = 
{\bf q}_{r_n(k)}^H \phi_{g_n(k)}({\bf A}) \widehat{\bf r}_{k-1}
 /
{\bf q}_{r_n(k)}^H {\bf A}
\phi_{g_n(k)}({\bf A}) \widehat{\bf g}_{k-1}
}$; \\
\>8. \>\>\> $\phi_{g_n(k)}({\bf A}) \widehat{\bf r}_k =
\phi_{g_n(k)}({\bf A}) \widehat{\bf r}_{k-1} - \alpha_k {\bf A}
\phi_{g_n(k)}({\bf A})
\widehat{\bf g}_{k-1}$; \\
\>9. \>\>\> $\phi_{g_n(k)+1}({\bf A}) \widehat{\bf r}_k =
\phi_{g_n(k)+1}({\bf A}) \widehat{\bf r}_{k-1} - \alpha_k {\bf A}
\phi_{g_n(k)+1}({\bf A})
\widehat{\bf g}_{k-1}$; \\
\>10.\>\>End\\
\>11.\>\> If $r_n(k) < n$\\
\>12. \>\>\>For $s = \max (k - n, 0), \cdots, g_n(k)n-1$ \\
\>13. \>\>\>\>$
\beta^{(k)}_{s} = - {\bf
q}_{r_n(s+1)}^H
 \left(
\phi_{g_n(s+1)+1}({\bf A}) \widehat{\bf r}_k
  + \right.
  $\\
\>\>\>\>\> $
\left. \sum_{t = \max (k - n, 0) }^{s-1}
\beta^{(k)}_t \rho_{g_n(s+1)+1}
  {\bf A} \phi_{g_n(s+1)}({\bf A})
\widehat{\bf g}_t
 \right)
\big/ \rho_{g_n(s+1)+1}
 {\bf q}_{r_n(s+1)}^H {\bf A} \phi_{g_n(s+1)}({\bf A})
\widehat{\bf g}_s
$; \\
\>14. \>\>\>End\\
\>15. \>\>\>$\beta^{(k)}_{g_n(k)n} = - {\bf q}_{1}^H
 \left(
\phi_{g_n(k)+1}({\bf A}) \widehat{\bf r}_k
  + \right.$\\
\>\>\>\>\>\>\>$\left. \sum_{t = \max (k - n, 0) }^{g_n(k)n-1}
\beta^{(k)}_t \rho_{g_n(k)+1}
  {\bf A}
\phi_{g_n(k)}({\bf A}) \widehat{\bf g}_t
 \right)
\big/ \rho_{g_n(k)+1}
 {\bf q}_{1}^H {\bf A} \phi_{g_n(k)}({\bf A}) \widehat{\bf
g}_{g_n(k)n}
$; \\
\>16. \>\>\>For $s = g_n(k)n+1, \cdots, k-1$ \\
\>17. \>\>\>\>$\beta^{(k)}_{s} = - {\bf q}_{r_n(s+1)}^H
 \left(
\phi_{g_n(s+1)+1}({\bf A}) \widehat{\bf r}_k
 +
 \sum_{t = \max (k - n, 0) }^{g_n(k)n}
\beta^{(k)}_t \rho_{g_n(s+1)+1} {\bf A} \phi_{g_n(s+1)}({\bf A})
\widehat{\bf g}_t + \right.
$\\
 \>\>\>\>\>\>$\left. \sum_{t =
g_n(k)n+1}^{s-1} \beta^{(k)}_t \rho_{g_n(s+1)+1} {\bf A}
\phi_{g_n(s+1)}({\bf A}) \widehat{\bf g}_t
 \right)
\big/ \rho_{g_n(s+1)+1}
 {\bf q}_{r_n(s+1)}^H {\bf A} \phi_{g_n(s+1)}({\bf A})
\widehat{\bf g}_s
$; \\
\>18. \>\>\>End \\
\>19.\>\>\> $\rho_{g_n(k)+1} {\bf A} \phi_{g_n(k)}({\bf A})
\widehat{\bf g}_k = \rho_{g_n(k)+1} {\bf A} \phi_{g_n(k)}({\bf A})
\widehat{\bf r}_k +$\\
\>\>\>\>\>\>\>\>\>\>\> $
 \sum_{s = \max (k - n, 0) }^{g_n(k)n}
\beta_{s}^{(k)} \rho_{g_n(k)+1} {\bf A} \phi_{g_n(k)}({\bf A})
 \widehat{\bf g}_{s} +$\\
\>\>\>\>\>\>\>\>\>\>\> $\sum_{s = g_n(k)n+1 }^{k -1} \beta_{s}^{(k)}
\rho_{g_n(k)+1} {\bf A} \phi_{g_n(k)}({\bf A})
 \widehat{\bf g}_{s}$; \\
\>20.\>\>\> $ \phi_{g_n(k)+1}({\bf A})
 \widehat{\bf g}_k =
\phi_{g_n(k)+1}({\bf A})
 \widehat{\bf r}_k +
\sum_{s = \max (k - n, 0) }^{g_n(k)n} \beta_{s}^{(k)}
\phi_{g_n(k)+1}({\bf A}) \widehat{\bf g}_{s} +$\\
\>\>\>\>\>\>\>\>\>    $\sum_{s = g_n(k)n+1 }^{k -1} \beta_{s}^{(k)}
\phi_{g_n(k)+1}({\bf A})
 \widehat{\bf g}_{s}$; \\
\>21.\>\>Else\\
\>22. \>\>\>$\displaystyle{\beta^{(k)}_{g_n(k)n} = -
 {\bf q}_{1}^H
\phi_{g_n(k)+1}({\bf A}) \widehat{\bf r}_k \big/ \rho_{g_n(k)+1}
 {\bf q}_{1}^H {\bf A}
\phi_{g_n(k)}({\bf A}) \widehat{\bf g}_{g_n(k)n}
}$; \\
\>23. \>\>\>For $s = g_n(k)n+1, \cdots, k-1$ \\
\>24. \>\>\>\>$\beta^{(k)}_{s} = - {\bf q}_{r_n(s+1)}^H
 \left(
\phi_{g_n(s+1)+1}({\bf A}) \widehat{\bf r}_k
 +
 \beta^{(k)}_{g_n(k)n}
\rho_{g_n(s+1)+1}
  {\bf A}
\phi_{g_n(s+1)}({\bf A}) \widehat{\bf g}_{g_n(k)n}
 + \right.$\\
\>\>\>\>\>\>$\left.
 \sum_{t = g_n(k)n+1}^{s-1} \beta^{(k)}_t
\rho_{g_n(s+1)+1}
  {\bf A}
\phi_{g_n(s+1)}({\bf A}) \widehat{\bf g}_t
 \right)
\big/ \rho_{g_n(s+1)+1}
 {\bf q}_{r_n(s+1)}^H {\bf A} \phi_{g_n(s+1)}({\bf A})
\widehat{\bf g}_s
$; \\
\>25. \>\>\>End \\
\>26.\>\>\> $\phi_{g_n(k)+1}({\bf A}) \widehat{\bf g}_k =
\phi_{g_n(k)+1}({\bf A}) \widehat{\bf r}_k + \beta_{g_n(k)n}^{(k)}
\phi_{g_n(k)+1}({\bf A}) \widehat{\bf g}_{g_n(k)n} +$\\
\>\>\>\>\>\>\>\>\>
 $\sum_{s = g_n(k)n+1}^{k -1} \beta_{s}^{(k)}
\phi_{g_n(k)+1}({\bf A}) \widehat{\bf g}_{s}$; \\
\>27.\>\>End\\
\>28. \> End
\end{tabbing}

\vspace{.2cm}

Lines 4, 8, 9, 19, 20 and 26, DS\#2, were obtained from Lines 4, 7,
17 and 23, DS\#1, through a multiplication by $\phi_{g_n(k)}({\bf
A}), \phi_{g_n(k)+1}({\bf A})$ and $\rho_{g_n(k)+1} {\bf A}
\phi_{g_n(k)}({\bf A})$ respectively. Line 5, DS\#2, is a direct
result of the definition (\ref{equ:12-23}) of $\phi$. These lines
are prepared for the updates of the vectors defined in
(\ref{equ:7-29}).

To help understand how DS\#1 is turned into DS\#2, let us
demonstrate (i) the transformation of Line 3, DS\#1, into Line 3,
DS\#2 and (ii) the transformation of the term ${\bf p}^H_{g_n(k)n+1}
{\bf A} \widehat{\bf r}_k$ on Line 13, DS\#1, into the term ${\bf
q}^H_{1} \phi_{g_n(k)+1}({\bf A}) \widehat{\bf r}_k$ on Line 15,
DS\#2, as follows.

\begin{enumerate}
\item[(i)] By Corollary \ref{cor:7-23},
$$
\alpha_k = \frac{{\bf p}_k^H \widehat{\bf r}_{k-1}}{{\bf p}_k^H {\bf
A} \widehat{\bf g}_{k-1}} = \frac{\frac{1}{c_{g_n(k)}^{(g_n(k))}}
{\bf q}_{r_n(k)}^H \phi_{g_n(k)}({\bf A}) \widehat{\bf r}_{k-1}}{
\frac{1}{c_{g_n(k)}^{(g_n(k))}} {\bf q}_{r_n(k)}^H {\bf A}
\phi_{g_n(k)}({\bf A})  \widehat{\bf g}_{k-1}} = \frac{ {\bf
q}_{r_n(k)}^H \phi_{g_n(k)}({\bf A}) \widehat{\bf r}_{k-1}}{  {\bf
q}_{r_n(k)}^H {\bf A} \phi_{g_n(k)}({\bf A}) \widehat{\bf g}_{k-1}}
$$
where $c_{g_n(k)}^{(g_n(k))}$ is the leading coefficient of
$\phi_{g_n(k)}(\lambda)$ (see (\ref{equ:5-28-09-10})).

\item[(ii)]
By (\ref{equ:7-9-5}) and Proposition \ref{prop:1}(a), we have
$$\begin{array}{rl}
{\bf A}^H {\bf p}_{g_n(k)n+1}
& = ({\bf A}^H)^{g_n(g_n(k)n+1)+1} {\bf q}_{r_n(g_n(k)n+1)}\\
&= ({\bf A}^H)^{g_n((g_n(k)+1)n+1)} {\bf
q}_{r_n((g_n(k)+1)n+1)}\\
& = {\bf p}_{(g_n(k)+1)n+1}.
\end{array}
$$
Hence ${\bf p}^H_{g_n(k)n+1} {\bf A} \widehat{\bf r}_k = {\bf
p}_{(g_n(k)+1)n+1}^H  \widehat{\bf r}_k$. Since $(g_n(k)+1)n+1 \leq
k + n$, an application of Corollary \ref{cor:7-23} to ${\bf
p}_{(g_n(k)+1)n+1}^H  \widehat{\bf r}_k$ thus yields
$$\begin{array}{rl}
{\bf p}^H_{g_n(k)n+1} {\bf A} \widehat{\bf r}_k &
=
\frac{1}{c_{g_n((g_n(k)+1)n+1)}^{(g_n((g_n(k)+1)n+1))}} {\bf
q}_{r_n((g_n(k)+1)n+1)}^H \phi_{g_n((g_n(k)+1)n+1)}({\bf A})
\widehat{\bf r}_k \\
&
= \frac{1}{c_{g_n(k)+1}^{(g_n(k)+1)}} {\bf q}_{1}^H
\phi_{g_n(k)+1}({\bf A}) \widehat{\bf r}_k.
\end{array}
$$
The second equation above follows from (\ref{equ:7-9-10}). The
coefficient $1/c_{g_n(k)+1}^{(g_n(k)+1)}$ is missed from Line 15,
DS\#2, because it was canceled out by the coefficient from the
denominator.\\
\end{enumerate}

Our goal is to establish updating relations for the quantities
introduced in (\ref{equ:7-29}).
To this end, we further transform DS\#2 into the following version.
This time, we work on the index function $g_n$ with the aid of
Proposition \ref{prop:1} so that the definitions in (\ref{equ:7-29})
can be applied. Again, further explanations are given after the
listing.\\

{\bf Derivation Stage \#3.}
\begin{tabbing}
x\=xxx\=xxx\=xxx\=xxx\=xxx\=xxx\=xxx\=xxx\=xxx\=xxx\=xxx\=xxx\kill
\>1. \>For $k = 1, 2, \cdots$, until convergence: \\
\>2.\>\>If $r_n(k) = 1$\\
\>3. \>\>\>$\displaystyle{\alpha_k = 
{\bf q}_{r_n(k)}^H \phi_{g_n(k-1)+1}({\bf A}) \widehat{\bf r}_{k-1}
 /
{\bf q}_{r_n(k)}^H {\bf A} \phi_{g_n(k-1)+1}({\bf A}) \widehat{\bf
g}_{k-1}
}$; \\
\>4. \>\>\> $\phi_{g_n(k)}({\bf A}) \widehat{\bf r}_k =
\phi_{g_n(k-1)+1}({\bf A}) \widehat{\bf r}_{k-1} - \alpha_k {\bf A}
\phi_{g_n(k-1)+1}({\bf A})
\widehat{\bf g}_{k-1}$; \\
\>5. \>\>\> $\phi_{g_n(k)+1}({\bf A}) \widehat{\bf r}_k =
(\rho_{g_n(k)+1} {\bf A} + {\bf I}) \phi_{g_n(k)}({\bf A})
\widehat{\bf r}_k
$; \\
\>6.\>\>Else\\
\>7. \>\>\>$\displaystyle{\alpha_k = 
{\bf q}_{r_n(k)}^H \phi_{g_n(k-1)}({\bf A}) \widehat{\bf r}_{k-1}
 /
{\bf q}_{r_n(k)}^H {\bf A} \phi_{g_n(k-1)}({\bf A}) \widehat{\bf
g}_{k-1}
}$; \\
\>8. \>\>\> $\phi_{g_n(k)}({\bf A}) \widehat{\bf r}_k =
\phi_{g_n(k-1)}({\bf A}) \widehat{\bf r}_{k-1} - \alpha_k {\bf A}
\phi_{g_n(k-1)}({\bf A})
\widehat{\bf g}_{k-1}$; \\
\>9. \>\>\> $\phi_{g_n(k)+1}({\bf A}) \widehat{\bf r}_k =
\phi_{g_n(k-1)+1}({\bf A}) \widehat{\bf r}_{k-1} - \alpha_k {\bf A}
\phi_{g_n(k-1)+1}({\bf A})
\widehat{\bf g}_{k-1}$; \\
\>10.\>\>End\\
\>11.\>\> If $r_n(k) < n$\\
\>12. \>\>\>For $s = \max (k - n, 0), \cdots, g_n(k)n-1$ \\
\>13. \>\>\>\>$
\beta^{(k)}_{s} = - {\bf q}_{r_n(s+1)}^H
 \left(
\phi_{g_n(k)}({\bf A}) \widehat{\bf r}_k
  + \right.
  $\\
\>\>\>\>\> $
\left. \sum_{t = \max (k - n, 0) }^{s-1} \beta^{(k)}_t
\rho_{g_n(t)+1}
  {\bf A} \phi_{g_n(t)}({\bf A})
\widehat{\bf g}_t
 \right)
\big/ \rho_{g_n(s)+1}
 {\bf q}_{r_n(s+1)}^H {\bf A} \phi_{g_n(s)}({\bf A})
\widehat{\bf g}_s
$; \\
\>14. \>\>\>End\\
\>15. \>\>\>$\beta^{(k)}_{g_n(k)n} = - {\bf q}_{1}^H
 \left(
\phi_{g_n(k)+1}({\bf A}) \widehat{\bf r}_k
  + \right.$\\
\>\>\>\>\>$\left. \sum_{t = \max (k - n, 0) }^{g_n(k)n-1}
\beta^{(k)}_t \rho_{g_n(k)+1}
  {\bf A}
\phi_{g_n(t)+1}({\bf A}) \widehat{\bf g}_t
 \right)
\big/ \rho_{g_n(k)+1}
 {\bf q}_{1}^H {\bf A} \phi_{g_n(g_n(k)n)+1}({\bf A}) \widehat{\bf
g}_{g_n(k)n}
$; \\
\>16. \>\>\>For $s = g_n(k)n+1, \cdots, k-1$ \\
\>17. \>\>\>\>$\beta^{(k)}_{s} = - {\bf q}_{r_n(s+1)}^H
 \left(
\phi_{g_n(k)+1}({\bf A}) \widehat{\bf r}_k
 +
 \sum_{t = \max (k - n, 0) }^{g_n(k)n}
\beta^{(k)}_t \rho_{g_n(k)+1} {\bf A} \phi_{g_n(t)+1}({\bf A})
\widehat{\bf g}_t + \right.
$\\
\>\>\>\>\>\>$\left. \sum_{t = g_n(k)n+1}^{s-1} \beta^{(k)}_t
\rho_{g_n(t)+1} {\bf A} \phi_{g_n(t)}({\bf A}) \widehat{\bf g}_t
 \right)
\big/ \rho_{g_n(s)+1}
 {\bf q}_{r_n(s+1)}^H {\bf A} \phi_{g_n(s)}({\bf A})
\widehat{\bf g}_s
$; \\
\>18. \>\>\>End \\
\>19.\>\>\> $\rho_{g_n(k)+1} {\bf A} \phi_{g_n(k)}({\bf A})
\widehat{\bf g}_k = \rho_{g_n(k)+1} {\bf A} \phi_{g_n(k)}({\bf A})
\widehat{\bf r}_k +$\\
\>\>\>\>\>\>\>\>\>\>\> $
 \sum_{s = \max (k - n, 0) }^{g_n(k)n}
\beta_{s}^{(k)} \rho_{g_n(k)+1} {\bf A} \phi_{g_n(s)+1}({\bf A})
 \widehat{\bf g}_{s} +$\\
\>\>\>\>\>\>\>\>\>\>\> $\sum_{s = g_n(k)n+1 }^{k -1} \beta_{s}^{(k)}
\rho_{g_n(s)+1} {\bf A} \phi_{g_n(s)}({\bf A})
 \widehat{\bf g}_{s}$; \\
\>20.\>\>\> $\phi_{g_n(k)+1}({\bf A})
 \widehat{\bf g}_k =
\phi_{g_n(k)+1}({\bf A})
 \widehat{\bf r}_k +
\sum_{s = \max (k - n, 0) }^{g_n(k)n} \beta_{s}^{(k)}
(\rho_{g_n(k)+1} {\bf A} + {\bf I}) \phi_{g_n(s)+1}({\bf A}) \widehat{\bf g}_{s} +$\\
\>\>\>\>\>\>\>\>\>    $\sum_{s = g_n(k)n+1 }^{k -1} \beta_{s}^{(k)}
\phi_{g_n(s)+1}({\bf A})
 \widehat{\bf g}_{s}$; \\
\>21.\>\>Else\\
\>22. \>\>\>$\displaystyle{\beta^{(k)}_{g_n(k)n} = -
 {\bf q}_{1}^H
\phi_{g_n(k)+1}({\bf A}) \widehat{\bf r}_k \big/ \rho_{g_n(k)+1}
 {\bf q}_{1}^H {\bf A}
\phi_{g_n(g_n(k)n)+1}({\bf A}) \widehat{\bf g}_{g_n(k)n}
}$; \\
\>23. \>\>\>For $s = g_n(k)n+1, \cdots, k-1$ \\
\>24. \>\>\>\>$\beta^{(k)}_{s} = - {\bf q}_{r_n(s+1)}^H
 \left(
\phi_{g_n(k)+1}({\bf A}) \widehat{\bf r}_k
 +
 \beta^{(k)}_{g_n(k)n}
\rho_{g_n(k)+1}
  {\bf A}
\phi_{g_n(g_n(k)n)+1}({\bf A}) \widehat{\bf g}_{g_n(k)n}
 + \right.$\\
\>\>\>\>\>\>$\left.
 \sum_{t = g_n(k)n+1}^{s-1} \beta^{(k)}_t
\rho_{g_n(t)+1}
  {\bf A}
\phi_{g_n(t)}({\bf A}) \widehat{\bf g}_t
 \right)
\big/ \rho_{g_n(s)+1}
 {\bf q}_{r_n(s+1)}^H {\bf A} \phi_{g_n(s)}({\bf A})
\widehat{\bf g}_s
$; \\
\>25. \>\>\>End \\
\>26.\>\>\> $\phi_{g_n(k)+1}({\bf A}) \widehat{\bf g}_k =
\phi_{g_n(k)+1}({\bf A}) \widehat{\bf r}_k + \beta_{g_n(k)n}^{(k)}
(\rho_{g_n(k)+1} {\bf A} + {\bf I})    \phi_{g_n(g_n(k)n)+1}({\bf A}) \widehat{\bf g}_{g_n(k)n} +$\\
\>\>\>\>\>\>\>\>\>
 $\sum_{s = g_n(k)n+1}^{k -1} \beta_{s}^{(k)}
\phi_{g_n(s)+1}({\bf A}) \widehat{\bf g}_{s}$; \\
\>27.\>\>End\\
\>28. \> End
\end{tabbing}

\vspace{.2cm}

As an example, let us show how the $g_n(s+1)$ inside the sum
$\sum_{t = \max(k-n, 0)}^{s-1} \cdots$ on Line 13, DS\#2, was
written as the $g_n(t)$ on Line 13, DS\#3.

If $g_n(k) = 0$, Line 13 of DS\#2 is not implemented because of the
conventions immediately following DS\#1. So, we assume that $g_n(k)
> 0$. Since
$$\max(k-n, 0) \leq s,\, t \leq g_n(k)n-1,$$
we have
$$
g_n(s+1) = g_n(k + 1) - 1 = g_n(t+1)
$$
by Proposition \ref{prop:1}(b). Now that $g_n(k) > 0$, $\max(k-n, 0)
= k-n$ and hence
\begin{equation} \label{equ:8-1-1}
k-n \leq  t \leq g_n(k)n-1.
\end{equation}
Let $k = jn + i$ as in (\ref{equ:8-1}). Then (\ref{equ:8-1-1}) is
$$
(j-1)n + i \leq  t \leq (j-1)n + n-1
$$
which implies that $r_n(t) < n$. Now, Proposition \ref{prop:1}(d)
yields $g_n(t+1) = g_n(t)$ and therefore we have $g_n(s+1) =
g_n(t)$.

Now we are ready to use the vectors defined in (\ref{equ:7-29}) and
(\ref{equ:9-13}). Substituting these vectors into DS\#3 leads to the
following stage.
\\

{\bf Derivation Stage \#4.}
\begin{tabbing}
x\=xxx\=xxx\=xxx\=xxx\=xxx\=xxx\=xxx\=xxx\=xxx\=xxx\=xxx\=xxx\kill
\>1. \>For $k = 1, 2, \cdots$, until convergence: \\
\>2.\>\>If $r_n(k) = 1$\\
\>3. \>\>\>$\displaystyle{\alpha_k = {\bf q}_{r_n(k)}^H {\bf
r}_{k-1}
 /
{\bf q}_{r_n(k)}^H {\bf A} {\bf g}_{k-1}
}$; \\
\>4. \>\>\> ${\bf u}_k = {\bf r}_{k-1} - \alpha_k {\bf A}
{\bf g}_{k-1}$; \\
\>5. \>\>\> ${\bf r}_k = \rho_{g_n(k)+1} {\bf A} {\bf u}_k
 + {\bf u}_k
$; \\
\>6.\>\>Else\\
\>7. \>\>\>$\displaystyle{\alpha_k = 
\rho_{g_n(k-1)+1} {\bf q}_{r_n(k)}^H {\bf u}_{k-1}
 /
{\bf q}_{r_n(k)}^H {\bf d}_{k-1}
}$; \\
\>8. \>\>\> ${\bf u}_k = {\bf u}_{k-1} - (\alpha_k/
\rho_{g_n(k-1)+1})  {\bf d}_{k-1}$; \\
\>9. \>\>\> ${\bf r}_k = {\bf r}_{k-1} - \alpha_k {\bf A}
{\bf g}_{k-1}$; \\
\>10.\>\>End\\
\>11.\>\> If $r_n(k) < n$\\
\>12. \>\>\>For $s = \max (k - n, 0), \cdots, g_n(k)n-1$ \\
\>13. \>\>\>\>$
\beta^{(k)}_{s} = - {\bf q}_{r_n(s+1)}^H
 \left(
{\bf u}_k
  +  \sum_{t = \max (k - n, 0) }^{s-1} \beta^{(k)}_t {\bf d}_t
 \right)
\big/
 {\bf q}_{r_n(s+1)}^H {\bf d}_s
$; \\
\>14. \>\>\>End\\
\>15. \>\>\>$\beta^{(k)}_{g_n(k)n} = - {\bf q}_{1}^H
 \left(
{\bf r}_k
  +
\rho_{g_n(k)+1} \sum_{t = \max (k - n, 0) }^{g_n(k)n-1}
\beta^{(k)}_t
  {\bf A}
{\bf g}_t
 \right)
\big/ \rho_{g_n(k)+1}
 {\bf q}_{1}^H {\bf A} {\bf
g}_{g_n(k)n}
$; \\
\>16. \>\>\>For $s = g_n(k)n+1, \cdots, k-1$ \\
\>17. \>\>\>\>$\beta^{(k)}_{s} = - {\bf q}_{r_n(s+1)}^H
 \left(
{\bf r}_k
 + \rho_{g_n(k)+1}
 \sum_{t = \max (k - n, 0) }^{g_n(k)n}
\beta^{(k)}_t  {\bf A} {\bf g}_t + \right.
$\\
\>\>\>\>\>\>\>$\left. \sum_{t = g_n(k)n+1}^{s-1} \beta^{(k)}_t {\bf
d}_t
 \right)
\big/
 {\bf q}_{r_n(s+1)}^H {\bf d}_s
$; \\
\>18. \>\>\>End \\
\>19.\>\>\> ${\bf d}_k = {\bf r}_k - {\bf u}_k + \rho_{g_n(k)+1}
 \sum_{s = \max (k - n, 0) }^{g_n(k)n}
\beta_{s}^{(k)} {\bf A} {\bf g}_{s} + \sum_{s = g_n(k)n+1 }^{k -1}
\beta_{s}^{(k)}
{\bf d}_{s}$; \\
\>20.\>\>\> ${\bf g}_k = {\bf r}_k + \sum_{s = \max (k - n, 0)
}^{g_n(k)n} \beta_{s}^{(k)} (\rho_{g_n(k)+1} {\bf A} + {\bf I}) {\bf
g}_{s} + \sum_{s = g_n(k)n+1 }^{k -1} \beta_{s}^{(k)}
{\bf g}_{s}$; \\
\>21.\>\>Else\\
\>22. \>\>\>$\displaystyle{\beta^{(k)}_{g_n(k)n} = -
 {\bf q}_{1}^H
{\bf r}_k \big/ \rho_{g_n(k)+1}
 {\bf q}_{1}^H {\bf A}
{\bf g}_{g_n(k)n}
}$; \\
\>23. \>\>\>For $s = g_n(k)n+1, \cdots, k-1$ \\
\>24. \>\>\>\>$\beta^{(k)}_{s} = - {\bf q}_{r_n(s+1)}^H
 \left(
{\bf r}_k
 +
\rho_{g_n(k)+1}   \beta^{(k)}_{g_n(k)n}
  {\bf A}
{\bf g}_{g_n(k)n}
 + 
 \sum_{t = g_n(k)n+1}^{s-1} \beta^{(k)}_t
{\bf d}_t
 \right)
\big/
 {\bf q}_{r_n(s+1)}^H {\bf d}_s
$; \\
\>25. \>\>\>End \\
\>26.\>\>\> ${\bf g}_k = {\bf r}_k + \beta_{g_n(k)n}^{(k)}
(\rho_{g_n(k)+1} {\bf A} + {\bf I}) {\bf g}_{g_n(k)n} + 
\sum_{s = g_n(k)n+1}^{k -1} \beta_{s}^{(k)}
{\bf g}_{s}$; \\
\>27.\>\>End\\
\>28. \> End
\end{tabbing}

\vspace{.2cm}

We consider ${\bf r}_k$ to be the residual of the $k$th approximate
solution ${\bf x}_k$. Updating relations for ${\bf x}_k$ can be
obtained from Lines 4, 5 and 9 respectively:
\begin{equation}\label{equ:8-8}
\begin{array}{rl}
{\bf x}_k = & \left\{ \begin{array}{lcl}
{\bf x}_{k-1} -
\rho_{g_n(k)+1} {\bf u}_k +
\alpha_k {\bf g}_{k-1}, & &\mbox{if } r_n(k) = 1 \\
 {\bf x}_{k-1} + {\alpha}_k {\bf g}_{k-1}, & & \mbox{if } r_n(k) > 1.
\end{array} \right.
\end{array}
\end{equation}
After adding (\ref{equ:8-8}) to DS\#4 and simplifying the operations
appropriately, we arrive at the following ML($n$)BiCGStab algorithm.
Just like BiCGStab, the free parameter $\rho_{g_n(k)+1}$ on Line 5,
DS\#4, is chosen to minimize the $2$-norm of ${\bf r}_k$.
\\

\begin{algorithm}\label{alg:2}
{\bf ML($n$)BiCGStab without preconditioning associated with definition (\ref{equ:7-29})} 
\begin{tabbing}
x\=xxx\= xxx\=xxx\=xxx\=xxx\=xxx\kill \>1. \> Choose an initial
guess ${\bf x}_0$ and $n$ vectors ${\bf q}_1, {\bf q}_2, \cdots,
{\bf q}_n$. \\
\>2. \>  Compute ${\bf r}_0 = {\bf b} - {\bf A} {\bf x}_0$ and set
${\bf g}_0 = {\bf r}_0$. Compute ${\bf w}_0 = {\bf A}
{\bf g}_0,\,\, c_0 = {\bf q}_1^H {\bf w}_0$. \\
\>3. \>  For $k = 1, 2, \cdots$, until convergence: \\
\>4. \>\>If $r_n(k) = 1$\\
\>5. \>\>\> $\alpha_k = {\bf q}^H_{r_n(k)} {\bf r}_{k-1} / c_{k-1}$;
\\
\>6. \>\>\>$
 {\bf u}_k =
{\bf r}_{k-1} - \alpha_k {\bf w}_{k-1}$; \\
\>7. \>\>\>${\bf x}_k = {\bf x}_{k-1} +
\alpha_k {\bf g}_{k-1}$; \\
 \>8. \>\>\>$
\rho_{g_n(k)+1} = -({\bf A} {\bf u}_k)^H {\bf u}_k / \| {\bf A} {\bf
u}_k \|^2_2
$;\\
\>9. \>\>\>${\bf x}_k = {\bf x}_{k} - \rho_{g_n(k)+1} {\bf u}_k$; \\
\>10. \>\>\>$ {\bf r}_k = \rho_{g_n(k)+1} {\bf A} {\bf u}_k
+ {\bf u}_k$;\\
\>11. \>\>Else\\
\>12. \>\>\> $\widetilde{\alpha}_k = {\bf q}^H_{r_n(k)} {\bf
u}_{k-1} / c_{k-1}$;  \,\,\,\,\,\,\,\,\,\,\,\,\,\,\,\,\,\,  \%
$\widetilde{\alpha}_k = \alpha_k / \rho_{g_n(k-1)+1}$
\\
\>13. \>\>\>If $r_n(k) < n$\\
\>14. \>\>\>\>${\bf u}_k = {\bf u}_{k-1} - \widetilde{\alpha}_k
{\bf d}_{k-1}$; \\
\>15. \>\>\>End\\
\>16. \>\>\>$  {\bf x}_k =
 {\bf x}_{k-1} + \rho_{g_n(k-1) + 1}
\widetilde{\alpha}_k {\bf g}_{k-1}$; \\
\>17. \>\>\>${\bf r}_k = {\bf r}_{k-1} - \rho_{g_n(k-1) + 1}
\widetilde{\alpha}_k
{\bf w}_{k-1}$; \\
\>18. \>\>End\\
\>19. \>\> If $r_n(k) < n$\\
\>20. \>\>\> ${\bf z}_d = {\bf u}_k,\,\, {\bf g}_k = {\bf 0}, \,\,
{\bf z}_w
= {\bf 0}$;\\
\>21. \>\>\>   For $s = k - n, \cdots, g_n(k) n - 1$
and $g_n(k) \geq 1$\\
\>22. \>\>\>\> $\beta^{(k)}_{s} = - {\bf q}^H_{r_n(s+1)} {\bf z}_d
\big/
c_s$;\\
\>23. \>\>\>\>${\bf z}_d = {\bf z}_d + \beta^{(k)}_{s} {\bf d}_s$;\\
\>24. \>\>\>\>${\bf g}_k = {\bf g}_k + \beta^{(k)}_{s} {\bf g}_s$;\\
\>25. \>\>\>\>${\bf z}_w = {\bf z}_w + \beta^{(k)}_{s} {\bf w}_s$;\\
\>26. \>\>\>            End \\
\>27. \>\>\>${\bf z}_w = {\bf r}_k + \rho_{g_n(k)+1}
{\bf z}_w$;\\
\>28. \>\>\> $\tilde{\beta}^{(k)}_{g_n(k) n} = - {\bf q}^H_{1} {\bf z}_w
 \big/
 c_{g_n(k) n}
$; \,\,\,\,\,\,\,\, \% $\tilde{\beta}^{(k)}_{g_n(k) n} = \rho_{g_n(k)+1} \beta^{(k)}_{g_n(k) n}$\\
\>29. \>\>\> ${\bf z}_w = {\bf z}_w +
\tilde{\beta}^{(k)}_{g_n(k) n} {\bf w}_{g_n(k) n}$;\\
\>30. \>\>\> ${\bf g}_k = {\bf g}_k + {\bf z}_w  +
(\tilde{\beta}^{(k)}_{g_n(k) n}/\rho_{g_n(k)+1}) {\bf g}_{g_n(k) n}$;\\
\>31. \>\>\>   For $s = g_n(k) n + 1, \cdots, k - 1$ \\
\>32. \>\>\>\> $\beta^{(k)}_{s} = - {\bf q}^H_{r_n(s+1)} {\bf z}_w
 \big/
c_s
$; \\
\>33. \>\>\>\>${\bf g}_k = {\bf g}_k + \beta_s^{(k)} {\bf g}_s$;\\
\>34. \>\>\>\>${\bf z}_w = {\bf z}_w + \beta_s^{(k)} {\bf d}_s$;\\
\>35. \>\>\>            End \\
\>36.\>\>\> $ {\bf d}_k =
{\bf z}_w - {\bf u}_k$;\\
\>37.\>\>\> $c_k = {\bf q}_{r_n(k+1)}^H {\bf d}_k$;\\
\>38.\>\>\> ${\bf w}_k = {\bf A}
{\bf g}_k$;\\
\>39. \>\> Else \\
\>40. \>\>\> $\tilde{\beta}^{(k)}_{g_n(k) n} = - {\bf q}^H_{1} {\bf r}_k
\big/
c_{g_n(k) n}$; \,\,\,\,\,\,\,\, \% $\tilde{\beta}^{(k)}_{g_n(k) n} = \rho_{g_n(k)+1} \beta^{(k)}_{g_n(k) n}$\\
\>41. \>\>\> ${\bf z}_w = {\bf r}_k +
\tilde{\beta}^{(k)}_{g_n(k)n}
{\bf w}_{g_n(k)n}$;\\
\>42. \>\>\> ${\bf g}_k = {\bf z}_w + (\tilde{\beta}^{(k)}_{g_n(k)n}/ \rho_{g_n(k)+1})
{\bf g}_{g_n(k)n}$;\\
\>43. \>\>\>   For $s = g_n(k) n+1, \cdots, k - 1$ \\
\>44. \>\>\>\> $\beta^{(k)}_{s} = - {\bf q}^H_{r_n(s+1)} {\bf z}_w
\big/
c_s$;\\
\>45. \>\>\>\> ${\bf g}_k = {\bf g}_k + \beta_s^{(k)} {\bf g}_s$;\\
\>46. \>\>\>\> ${\bf z}_w = {\bf z}_w + \beta_s^{(k)} {\bf d}_s$;\\
\>47. \>\>\>            End \\
\>48.\>\>\> ${\bf w}_k = {\bf A}
{\bf g}_k$;\\
\>49.\>\>\> $c_k = {\bf q}^H_{r_n(k+1)} {\bf w}_k$;\\
\>50. \>\> End\\
\>51. \> End
\end{tabbing}
\end{algorithm}
\vspace{.2cm}

{\it Remarks}: \begin{enumerate} \item[(i)] Algorithm \ref{alg:2}
does not compute the quantities ${\bf u}_k$ and ${\bf d}_k$ when
$r_n(k) = n$ (see Lines 13-15 and Lines 39-50).
 \item[(ii)] if the ${\bf
u}_k$ on Line 6 happens to be zero, then the $\rho_{g_n(k)+1}$ on
Line 8 and therefore the ${\bf x}_k$ and ${\bf r}_k$ on Lines 9 and
10 will not be computable. In this case, however, the ${\bf x}_k$ on
Line 7 will be the exact solution to system (\ref{equ:7-9-3}) and
Algorithm \ref{alg:2} stops there.\\
\end{enumerate}

We now compare Algorithm \ref{alg:2} with the ML($n$)BiCGStab
algorithm 
in \cite{yeungchan}. First, the definitions of ${\bf r}_k$, ${\bf
u}_k$ and ${\bf g}_k$ are the same in both algorithms, but ${\bf
d}_k$ is defined differently. In \cite{yeungchan}, ${\bf d}_k =
\phi_{g_n(k)}({\bf A}) \widehat{\bf g}_k$. In exact arithmetic,
however, both algorithms compute the same $\rho_{g_n(k)+1}$, ${\bf
r}_k$ and ${\bf x}_k$. Second, the derivation of Algorithm
\ref{alg:2} has been made simpler by using index functions. As a
result, some redundant operations in Algorithm 2 of \cite{yeungchan}
can been seen and removed and some arithmetics are simplified. For
example, the vectors ${\bf d}_k, {\bf u}_k$ are computed in every
iteration in Algorithm 2 of \cite{yeungchan}. They are now computed
only when $r_n(k) < n$. Also, the expression of
$\beta^{(k)}_{g_n(k)n}$ on Line 39 of Algorithm \ref{alg:2} is
simpler. Some other minor changes were also made so that the
algorithm becomes more efficient.

\begin{table}[tbp]
\footnotesize \caption{Average cost per ($k$-)iteration of Algorithm
\ref{alg:10-22} and its storage requirement.}
\begin{center}
\begin{tabular}{|c|c|c|c|}  \hline
Preconditioning (${\bf M}^{-1} {\bf v}$) & $\displaystyle{1 +
\frac{1}{n}}$ & Vector addition (${\bf u} \pm {\bf v}$) &
$\displaystyle{2 - \frac{2}{n}}$ \\ \hline Matvec (${\bf A} {\bf
v}$) & $\displaystyle{1 + \frac{1}{n}}$ & Saxpy (${\bf u} + \alpha
{\bf v}$) & $\displaystyle{\max (2.5 n + 2.5 - \frac{2}{n}, 6)}$
 \\ \hline
dot product ($\displaystyle{ {\bf u}^H {\bf v}}$) & $\displaystyle{n
+ 1+ \frac{2}{n}}$  &Storage  &${\bf A} + {\bf M} + (4 n+5) N +
O(n)$
\\ \hline
\end{tabular}
\end{center} \label{tab:10-28-1}
\end{table}

Computational cost and storage requirement of Algorithm \ref{alg:2},
obtained based on its preconditioned version, Algorithm
\ref{alg:10-22} in \S\ref{sec:apen}, are summarized in Table
\ref{tab:10-28-1}. Since the vectors $\{{\bf q}_1, \ldots, {\bf
q}_n\}$, $\{{\bf d}_{k - n}, \ldots, {\bf d}_{g_n(k)n-1}, {\bf
d}_{g_n(k)n+1}, \ldots, {\bf d}_{k-1}\}$, $\{{\bf g}_{k-n}, \ldots,
{\bf g}_{k-1}\}$ and $\{{\bf w}_{k-n}, \ldots, {\bf w}_{g_n(k)n},
{\bf w}_{k-1}\}$ are required
in iteration $k$, they must be stored. When $n$ is large, this
storage is dominant. So, the storage requirement of the algorithm is
about $4 n N$.

\subsection{Properties} We summarize the properties of Algorithm
\ref{alg:2} in the following proposition. Since ${\bf r}_0 =
\widehat{\bf r}_0$ by (\ref{equ:9-13}), $\nu$ (see
\S\ref{sec:5-18}) is also the degree of the minimal polynomial of
${\bf r}_0$ with respect to ${\bf A}$.  \\

\begin{proposition}
\label{prop:10-1} Under the assumptions of Proposition
\ref{prop:7-11-1}, if $\rho_{g_n(k)+1} \ne 0$ and $-
1/\rho_{g_n(k)+1} \not\in \sigma({\bf A})$ for $1 \leq k \leq
\nu-1$, where $\sigma({\bf A})$ is the spectrum of $\bf A$, then
Algorithm \ref{alg:2} does not break down by zero division for $k =
1, 2, \cdots, \nu$, and ${\bf x}_\nu$ is the exact solution of
(\ref{equ:7-9-3}). Moreover, the computed quantities satisfy
\begin{enumerate}
\item[(a)]
${\bf x}_k \in {\bf x}_0 + {\mathcal K}_{g_n(k)+k+1}({\bf A}, {\bf r}_0)
$ and ${\bf r}_k = {\bf b} - {\bf A} {\bf x}_k \in {\bf r}_0 + {\bf
A} {\mathcal K}_{g_n(k)+k+1}({\bf A}, {\bf r}_0)
$ for $1 \leq k \leq \nu-1$.

\item[(b)] ${\bf r}_k \ne {\bf 0}$ for $1 \leq k \leq \nu -1$ and ${\bf
r}_\nu = {\bf 0}$.

\item[(c)] ${\bf r}_k \not\perp {\bf
q}_{1}$ for $1 \leq k \leq \nu -1$ with $r_n(k) = n$.

\item[(d)] ${\bf u}_k \perp span \{ {\bf q}_1, {\bf q}_2, \cdots,
{\bf q}_{r_n(k)}\}$ and ${\bf u}_k \not\perp {\bf q}_{r_n(k)+1}$ for
$1 \leq k \leq \nu -1$ with $r_n(k) < n$.

\item[(e)] ${\bf d}_k \perp span \{ {\bf q}_1, {\bf q}_2, \cdots,
{\bf q}_{r_n(k)}\}$ and ${\bf d}_k
 \not\perp {\bf q}_{r_n(k)+1}$ for
$1 \leq k \leq \nu -1$ with $r_n(k) < n$.
\end{enumerate}
\end{proposition}

\vspace{.2cm}

\proof The divisors in Algorithm \ref{alg:2} are $c_k, \|{\bf A}
{\bf u}_k\|^2_2$ and $\rho_{g_n(k)+1}$ respectively, where the
$\rho$'s have been assumed to be nonzero. By Proposition
\ref{prop:7-11-1}(c), we have ${\bf A} \widehat{\bf r}_k \ne {\bf
0}$ for $1 \leq k \leq \nu -1$. Since $- 1/\rho \not\in \sigma({\bf
A})$ by assumption, $\phi_{g_n(k)} ({\bf A})$ is nonsingular. Hence
${\bf A} {\bf u}_k = \phi_{g_n(k)} ({\bf A}) \,{\bf A} \widehat{\bf
r}_k \ne {\bf 0}$ (see (\ref{equ:7-29}) for the first equation).
Therefore, $\|{\bf A} {\bf u}_k\|_2 \ne 0$ for $1 \leq k \leq \nu
-1$.

$c_k$ is defined respectively on Lines 37 and 49 in the algorithm.
When $r_n(k) < n$, we have $c_k = {\bf q}_{r_n(k+1)}^H {\bf d}_k$.
In this case, $c_k = \rho_{g_n(k)+1} {\bf q}_{r_n(k+1)}^H {\bf A}
\phi_{g_n(k)}({\bf A}) \widehat{\bf g}_k = \rho_{g_n(k)+1} {\bf
q}_{r_n(k+1)}^H {\bf A} \phi_{g_n(k+1)}({\bf A}) \widehat{\bf g}_k =
\rho_{g_n(k)+1} c_{g_n(k+1)}^{(g_n(k+1))} {\bf p}^H_{k+1} {\bf A}
\widehat{\bf g}_k = \rho_{g_n(k)+1} c_{g_n(k)}^{(g_n(k))} {\bf
p}^H_{k+1} {\bf A} \widehat{\bf g}_k$ $= c_{g_n(k)+1}^{(g_n(k)+1)}
{\bf p}^H_{k+1} {\bf A} \widehat{\bf g}_k$ by (\ref{equ:7-29}),
Proposition \ref{prop:1}(d), Corollary \ref{cor:7-23},
(\ref{equ:5-28-09-10}) and (\ref{equ:7-27}). Since the $\rho$'s are
nonzero and ${\bf p}^H_{k+1} {\bf A} \widehat{\bf g}_k \ne 0$ by
Proposition
\ref{prop:7-11-1}(g), we have $c_{g_n(k)+1}^{(g_n(k)+1)} 
\ne 0$ 
and hence $c_k \ne 0$. When $r_n(k) = n$, on the other hand, $c_k =
{\bf q}^H_{r_n(k+1)} {\bf w}_k = {\bf q}^H_{r_n(k+1)} {\bf A} {\bf
g}_k$. In this case, $c_k = {\bf q}^H_{r_n(k+1)} {\bf A}
\phi_{g_n(k)+1}({\bf A}) \,\widehat{\bf g}_k = {\bf q}^H_{r_n(k+1)}
{\bf A} \phi_{g_n(k+1)}({\bf A}) \,\widehat{\bf g}_k$ $=
c_{g_n(k+1)}^{(g_n(k+1))} {\bf p}^H_{k+1} {\bf A} \widehat{\bf g}_k
= c_{g_n(k)+1}^{(g_n(k)+1)} {\bf p}^H_{k+1} {\bf A} \widehat{\bf
g}_k \ne 0$.
Therefore, in either case, we
always have $c_k \ne 0$ for $1 \leq k \leq \nu -1$. Moreover, $c_0 =
{\bf q}_1^H {\bf w}_0 = {\bf q}_1^H {\bf A} {\bf g}_0$ according to
Line 2 of the algorithm. Since ${\bf p}_1 = {\bf q}_1$ by
(\ref{equ:7-9-5}) and ${\bf g}_0 = \widehat{\bf g}_0$ by
(\ref{equ:9-13}), $c_0 \ne 0$ by Proposition \ref{prop:7-11-1}(g).

Now that $\|{\bf A} {\bf u}_k\|_2 \ne 0$ and $\rho_{g_n(k)+1} \ne 0$
for $1 \leq k \leq \nu-1$ and $c_k \ne 0$ for $0 \leq k \leq \nu-1$,
Algorithm \ref{alg:2} does not break down by zero division in the
first $\nu -1$ iterations. When $k = \nu$, ${\bf u}_k = {\bf u}_\nu
= \phi_{g_n(\nu)}({\bf A}) \widehat{\bf r}_\nu = {\bf 0}$ and ${\bf
r}_k = {\bf r}_\nu = \phi_{g_n(\nu)+1}({\bf A}) \widehat{\bf r}_\nu
= {\bf 0}$ due to $\widehat{\bf r}_\nu = {\bf 0}$ by Proposition
\ref{prop:7-11-1}. If it happens that $r_n(\nu) = 1$, then the ${\bf
x}_k ( = {\bf x}_\nu)$ on Line 7 is the exact solution to system
(\ref{equ:7-9-3}) because its residual ${\bf u}_\nu$ is zero. So,
the algorithm stops there. Otherwise, the ${\bf x}_k (= {\bf
x}_\nu)$ on Line 16 will be exact with residual ${\bf r}_\nu = {\bf
0}$ and where the algorithm stops.

Part (a) follows from the definition of ${\bf r}_k$ in
(\ref{equ:7-29}) and Proposition \ref{prop:7-11-1}(a).

Since $\widehat{\bf r}_k \ne {\bf 0}$ for $1 \leq k \leq \nu -1$ by
Proposition \ref{prop:7-11-1}(b) and $\phi_{g_n(k)+1}({\bf A})$ is
nonsingular due to $-1/\rho \not\in \sigma({\bf A})$, we have ${\bf
r}_k = \phi_{g_n(k)+1}({\bf A}) \widehat{\bf r}_k \ne {\bf 0}$.
Therefore, Part (b) holds.

For Part (c), write $k = j n + n$ with $0 \leq j$. By
(\ref{equ:7-29}), (\ref{equ:5-28-09-10}) and Corollary
\ref{cor:7-23}, we have
 ${\bf q}_1^H {\bf r}_k = {\bf
q}_1^H \phi_{g_n(k)+1}({\bf A}) \widehat{\bf r}_k = {\bf
q}_{r_n((j+1)n+1)}^H \phi_{g_n((j+1)n+1)}({\bf A}) \widehat{\bf r}_k
= {\bf q}_{r_n(k+1)}^H \phi_{g_n(k+1)}({\bf A}) \widehat{\bf r}_k$
$= c_{g_n(k+1)}^{(g_n(k+1))} {\bf p}_{k+1}^H \widehat{\bf r}_k =
c_{g_n(k)+1}^{(g_n(k)+1)}
 \,{\bf p}_{k+1}^H
\widehat{\bf r}_k$. Now Part (c) follows from Proposition
\ref{prop:7-11-1}(d) and $c_{g_n(k)+1}^{(g_n(k)+1)} \ne 0$.

For the proof of Part (d), we first note that Algorithm \ref{alg:2}
does not compute ${\bf u}_k$ when $r_n(k) = n$ (see Lines 13 - 15).
Write $k = j n + i$ as in (\ref{equ:8-1}) and let $1 \leq t \leq i <
n$. Then $r_n(k) = i, \,g_n(k) = j = g_n (jn+t)$ and $r_n(jn+t) =
t$. Now, by (\ref{equ:7-29}) and Corollary \ref{cor:7-23}, we have
 ${\bf q}_t^H {\bf u}_k = {\bf
q}_t^H \phi_{g_n(k)}({\bf A}) \widehat{\bf r}_k = {\bf
q}_{r_n(jn+t)}^H \phi_{g_n(jn+t)}({\bf A}) \widehat{\bf r}_k =
c_{g_n(jn+t)}^{(g_n(jn+t))} \,{\bf p}_{jn+t}^H \widehat{\bf r}_k$.
Since ${\bf p}_{jn+t}^H \widehat{\bf r}_k = 0$ by Proposition
\ref{prop:7-11-1}(d), ${\bf q}_t^H {\bf u}_k = 0$ for $1 \leq t \leq
i$. Similarly, ${\bf q}_{i+1}^H {\bf u}_k = c_{g_n(k)}^{(g_n(k))}
 \,{\bf p}_{jn+i+1}^H \widehat{\bf
r}_k = c_{g_n(k)}^{(g_n(k))} \,{\bf p}_{k+1}^H \widehat{\bf r}_k$
(the validity of the first equation requires $i < n$). Because of
Proposition \ref{prop:7-11-1}(d) and $c_{g_n(k)}^{(g_n(k))} \ne 0$,
${\bf q}_{i+1}^H {\bf u}_k \ne 0$.

Similar to the quantity ${\bf u}_k$, Algorithm \ref{alg:2} does not
compute ${\bf d}_k$ when $r_n(k) = n$ (see Lines 40 - 49). By
(\ref{equ:7-29}), ${\bf d}_k = \rho_{g_n(k)+1} {\bf A} \phi_{g_n(k)}({\bf A})
\,\widehat{\bf g}_k$ and the proof of Part (e) is parallel to that
of Part (d). \hfill{\rule{2mm}{2mm}} \vspace{.2cm}

The conditions of $\rho_{g_n(k)+1} \ne 0$ and $- 1/\rho_{g_n(k)+1}
\not\in \sigma({\bf A})$ can be easily made satisfied. For example, one
can add some small random noise (e.g., $N(0, \delta)$ with $\delta
\ll 1$) to $\rho_{g_n(k)+1}$ after it is
computed.\\

\begin{corollary} \label{cor:3-28}
Consider the case where $n = 1$, (\ref{equ:7-9-3}) is a real system
and ${\bf q}_1 \in {\mathcal R}^N$ is a random vector with iid elements
from $N(0, 1)$. If some small random number
is added to $\rho_{g_n(k)+1}$ after
it is computed so that
 $\rho_{g_n(k)+1} \ne 0$ and
$-1/\rho_{g_n(k)+1} \not\in \sigma({\bf A})$, then Algorithm
\ref{alg:2} will work almost surely without breakdown by zero
division to find a solution of (\ref{equ:7-9-3}) from the affine
space ${\bf x}_0 + span \{ {\bf A}^t {\bf r}_0 | t \in {\mathcal N}_0\}$
provided that ${\bf x}_0 \in {\mathcal R}^{N}$ is chosen such that the
affine space contains a solution to (\ref{equ:7-9-3}).
\end{corollary}
\vspace{.2cm}

Proposition \ref{prop:10-1} indicates that exact solution can only
be found at iteration $k = \nu$. It is possible, however, that
$\|{\bf r}_k\|_2$ can become very small for some $k < \nu$. In
practice, we terminate the algorithm when $\|{\bf r}_k\|_2$ falls
within a given tolerance.


As in the case of BiCGStab, ML($n$)BiCGStab can encounter a breakdown in its implementation.
ML($n$)BiCGStab, besides the two
types of breakdown of ML($n$)BiCG, has one more type of breakdown
caused by $\rho_{g_n(k)+1}$. In more detail, the divisors in
Algorithm \ref{alg:2} are $c_k, \|{\bf A} {\bf u}_k\|^2_2$ and
$\rho_{g_n(k)+1}$. If $\|{\bf A} {\bf u}_k\|_2 = 0$, then
$\rho_{g_n(k)+1} = \infty$ and a breakdown due to the overflow
of $\rho_{g_n(k)+1}$ occurs. Under the assumptions of Proposition
\ref{prop:10-1}, on the other hand, it can be shown (see the proof
of the proposition) that
$c_k = c_{g_n(k)+1}^{(g_n(k)+1)} {\bf p}^H_{k+1} {\bf A}
\widehat{\bf g}_k$, where $c_{g_n(k)+1}^{(g_n(k)+1)}$ is the leading
coefficient of $\phi_{g_n(k)+1}(\lambda)$ (see
(\ref{equ:5-28-09-10})). So, $c_k$ is a quantity that relates to
$\rho_{g_n(k)+1}$ and the ML($n$)BiCG divisor ${\bf p}^H_{k+1} {\bf
A} \widehat{\bf g}_k$. Thus, either $\rho_{g_n(k)+1} = 0$ or
${\bf p}^H_{k+1} {\bf A} \widehat{\bf g}_k = 0$ can cause $c_k
= 0$.


\section{A Second ML($n$)BiCGStab Algorithm}\label{sec:12-24-3}
If we write $k = j n + i$ as in (\ref{equ:8-1}), the ${\bf r}_k$
defined by (\ref{equ:9-13}) then becomes
\begin{equation}\label{equ:9-24}
{\bf r}_{jn+i} = \phi_{j+1}({\bf A}) \,\widehat{\bf r}_{jn+i}
\end{equation}
where $i = 1, 2, \cdots, n$ and $j = 0, 1, 2, \cdots$.

Starting with $k = 1$, let us call every $n$ consecutive
$k$-iterations a ``cycle'', namely, iterations $k = 1, 2, \cdots, n$
form the first cycle, iterations $k = n+1, n+2, \cdots, n+n$ the
second cycle and so on. Then (\ref{equ:9-24}) increases the degree
of the polynomial $\phi$ by $1$ at the beginning of every cycle. For
example, consider $n = 3$. Then (\ref{equ:9-24}) implies that
$$
\begin{array}{ccccc}
{\bf r}_{1} = \phi_{1}({\bf A}) \,\widehat{\bf r}_{1}, & & {\bf r}_{4} = \phi_{2}({\bf A}) \,\widehat{\bf r}_{4},
& &{\bf r}_{7} = \phi_{3}({\bf A}) \,\widehat{\bf r}_{7},\\
{\bf r}_{2} = \phi_{1}({\bf A}) \,\widehat{\bf r}_{2},&  & {\bf
r}_{5} = \phi_{2}({\bf A}) \,\widehat{\bf r}_{5},& & {\bf r}_{8} =
\phi_{3}({\bf A}) \,\widehat{\bf r}_{8},
\\
{\bf r}_{3} = \phi_{1}({\bf A}) \,\widehat{\bf r}_{3},& &{\bf r}_{6}
= \phi_{2}({\bf A}) \,\widehat{\bf r}_{6},& &{\bf r}_{9} =
\phi_{3}({\bf A}) \,\widehat{\bf r}_{9}.
\end{array}
$$
Iteration $k = 4$ is the first iteration of the second cycle and the
degree of $\phi$ is increased from $1$ to $2$ there.

One can define ${\bf r}_k$ by increasing the degree of $\phi$ by one
anywhere within a cycle. Correspondingly, (we believe) the
definition will lead to a different algorithm of ML($n$)BiCGStab. As
an illustration, let us increase the degree of $\phi$ at the end of
every cycle and derive the algorithm associated with it.

\subsection{Notation and Definitions}
Let $\phi_k(\lambda)$ be defined as in (\ref{equ:12-23}). For $k
\in {\mathcal N}$, define
\begin{equation}\label{equ:9-24-1}
\begin{array}{lcl}
{\bf r}_k = \phi_{g_n(k+1)}({\bf A}) \,\widehat{\bf r}_k, & & {\bf
g}_k = \phi_{g_n(k+1)}({\bf A}) \widehat{\bf g}_k,\\
 {\bf u}_k = \phi_{g_n(k)}({\bf A}) \widehat{\bf r}_k, & & {\bf w}_k = {\bf A} {\bf g}_k.
\end{array}
\end{equation}
and set
\begin{equation}\label{equ:10-1-5}
\begin{array}{lcl}
{\bf r}_0 = \widehat{\bf r}_0 & \mbox{and} & {\bf g}_0 =
\widehat{\bf g}_0.
\end{array}
\end{equation}
The vector ${\bf r}_k$ is considered to be the residual of the
approximate solution ${\bf x}_k$ computed. We remark that ${\bf r}_k
= {\bf u}_k$ when $r_n(k) < n$ since $g_n(k+1) = g_n(k)$ in this
case.

Definition (\ref{equ:9-24-1}) increases the degree of $\phi$ at the
end of a cycle. To see this, let $n = 3$. Then (\ref{equ:9-24-1})
yields
$$
\begin{array}{ccccc}
{\bf r}_{1} = \phi_{0}({\bf A}) \,\widehat{\bf r}_{1}, & & {\bf
r}_{4} = \phi_{1}({\bf A}) \,\widehat{\bf r}_{4},
& &{\bf r}_{7} = \phi_{2}({\bf A}) \,\widehat{\bf r}_{7},\\
{\bf r}_{2} = \phi_{0}({\bf A}) \,\widehat{\bf r}_{2},&  & {\bf
r}_{5} = \phi_{1}({\bf A}) \,\widehat{\bf r}_{5},& & {\bf r}_{8} =
\phi_{2}({\bf A}) \,\widehat{\bf r}_{8},
\\
{\bf r}_{3} = \phi_{1}({\bf A}) \,\widehat{\bf r}_{3},& &{\bf r}_{6}
= \phi_{2}({\bf A}) \,\widehat{\bf r}_{6},& &{\bf r}_{9} =
\phi_{3}({\bf A}) \,\widehat{\bf r}_{9}.
\end{array}
$$

\subsection{Algorithm Derivation} To derive the algorithm associated
with (\ref{equ:9-24-1}), we first transform Algorithm \ref{alg:1}
(forgetting Lines 1, 2, 5 and 11) into the following version which
is computationally equivalent to Algorithm \ref{alg:1}, but is more
convenient for us to apply Proposition \ref{prop:1}.\\

{\bf Derivation Stage \#5.} \vspace{.2cm}
\begin{tabbing}
x\=xxx\= xxx\=xxx\=xxx\=xxx\=xxx\kill
\>1. \>For $k = 1, 2, \cdots$, until convergence: \\
\>2. \>\>$\alpha_k = {\bf p}_k^H \widehat{\bf r}_{k-1} / {\bf p}_k^H
{\bf A} \widehat{\bf g}_{k-1}$; \\
\>3.\>\> If $r_n(k) < n$\\
\>4. \>\>\> $\widehat{\bf r}_k = \widehat{\bf r}_{k-1} - \alpha_k {\bf A} \widehat{\bf g}_{k-1}$; \\
\>5. \>\>\>For $s = \max (k - n, 0), \cdots, g_n(k)n - 1$ \\
\>6. \>\>\>\>$\beta^{(k)}_{s} = - {\bf p}^H_{s+1} {\bf A}
 \left(\widehat{\bf r}_k + \sum_{t = \max (k - n, 0) }^{s-1} \beta^{(k)}_t \widehat{\bf g}_t \right)
\big/{\bf p}^H_{s+1} {\bf A} \widehat{\bf g}_s$; \\
\>7. \>\>\>End \\
\>8. \>\>\>For $s = g_n(k)n, \cdots, k - 1$ \\
\>9. \>\>\>\>$\beta^{(k)}_{s} = - {\bf p}^H_{s+1} {\bf A}
 \left(\widehat{\bf r}_k + \sum_{t = \max (k - n, 0) }^{g_n(k)n-1} \beta^{(k)}_t \widehat{\bf g}_t
+ \sum_{t = g_n(k)n }^{s-1} \beta^{(k)}_t \widehat{\bf g}_t
 \right)
\big/{\bf p}^H_{s+1} {\bf A} \widehat{\bf g}_s$; \\
\>10. \>\>\>End \\
\>11.\>\>\> $\widehat{\bf g}_k = \widehat{\bf r}_k + \sum_{s = \max
(k - n, 0) }^{g_n(k)n -1} \beta_{s}^{(k)}
            \widehat{\bf g}_{s}
+ \sum_{s = g_n(k)n }^{k -1} \beta_{s}^{(k)}
            \widehat{\bf g}_{s}
$; \\
\>12.\>\>Else\\
\>13. \>\>\> $\widehat{\bf r}_k = \widehat{\bf r}_{k-1} - \alpha_k {\bf A} \widehat{\bf g}_{k-1}$; \\
\>14. \>\>\>For $s = g_n(k)n, \cdots, k - 1$ \\
\>15. \>\>\>\>$\beta^{(k)}_{s} = - {\bf p}^H_{s+1} {\bf A}
 \left(\widehat{\bf r}_k +
\sum_{t = g_n(k)n }^{s-1} \beta^{(k)}_t \widehat{\bf g}_t
 \right)
\big/{\bf p}^H_{s+1} {\bf A} \widehat{\bf g}_s$; \\
\>16. \>\>\>End \\
\>17.\>\>\> $\widehat{\bf g}_k = \widehat{\bf r}_k + \sum_{s =
g_n(k)n }^{k -1} \beta_{s}^{(k)}
            \widehat{\bf g}_{s}$; \\
\>18.\>\>End\\
\>19. \> End
\end{tabbing}
\vspace{.2cm}

Then we transform DS\#5 as follows by Corollary
\ref{cor:7-23}.\\

{\bf Derivation Stage \#6.} \vspace{.2cm}
\begin{tabbing}
x\=xxx\= xxx\=xxx\=xxx\=xxx\=xxx\=xxx\=xxx\=xxx\kill
\>1. \>For $k = 1, 2, \cdots$, until convergence: \\
\>2. \>\>$\alpha_k = {\bf q}_{r_n(k)}^H \phi_{g_n(k)} ({\bf A})
\widehat{\bf r}_{k-1} / {\bf q}_{r_n(k)}^H
{\bf A} \phi_{g_n(k)}({\bf A}) \widehat{\bf g}_{k-1}$; \\
\>3.\>\> If $r_n(k) < n$\\
\>4. \>\>\> $\phi_{g_n(k)}({\bf A}) \widehat{\bf r}_k =
\phi_{g_n(k)}({\bf A}) \widehat{\bf r}_{k-1} - \alpha_k {\bf A}
\phi_{g_n(k)}({\bf A})
\widehat{\bf g}_{k-1}$; \\
\>5. \>\>\>For $s = \max (k - n, 0), \cdots, g_n(k)n - 1$ \\
\>6. \>\>\>\>$\beta^{(k)}_{s} = - {\bf q}^H_{r_n(s+1)}
 \left(
  \phi_{g_n(s+1)+1}({\bf
A}) \widehat{\bf r}_k \right.$\\
\>\>\>\>\> $\left. + \,\rho_{g_n(s+1)+1} \sum_{t = \max (k - n, 0)
}^{s-1} \beta^{(k)}_t {\bf A} \phi_{g_n(s+1)}({\bf A}) \widehat{\bf
g}_t \right) \big/\rho_{g_n(s+1)+1}
 {\bf q}^H_{r_n(s+1)} {\bf A} \phi_{g_n(s+1)} ({\bf A}) \widehat{\bf g}_s$; \\
\>7. \>\>\>End
\\
\>8. \>\>\>For $s = g_n(k)n, \cdots, k - 1$ \\
\>9. \>\>\>\>$\beta^{(k)}_{s} = - {\bf q}^H_{r_n(s+1)} {\bf A}
 \left(
 \phi_{g_n(s+1)}({\bf A}) \widehat{\bf
r}_k + \sum_{t = \max (k - n, 0) }^{g_n(k)n-1} \beta^{(k)}_t
 \phi_{g_n(s+1)}
({\bf A})
 \widehat{\bf g}_t \right.$\\
\>\>\>\>\>\>\> $\left. + \sum_{t = g_n(k)n }^{s-1} \beta^{(k)}_t
 \phi_{g_n(s+1)}
({\bf A}) \widehat{\bf g}_t
 \right)
\big/ {\bf q}^H_{r_n(s+1)} {\bf A} \phi_{g_n(s+1)} ({\bf A}) \widehat{\bf g}_s$; \\
\>10. \>\>\>End
\\
\>11.\>\>\> $\phi_{g_n(k+1)} ({\bf A}) \widehat{\bf g}_k =
\phi_{g_n(k+1)} ({\bf A}) \widehat{\bf r}_k + (\rho_{g_n(k+1)} {\bf
A} + {\bf I}) \sum_{s = \max (k - n, 0) }^{g_n(k)n -1}
\beta_{s}^{(k)} \phi_{g_n(k+1)-1} ({\bf A})
            \widehat{\bf g}_{s}
+$\\
\>\>\>\>\>\>\>\>\> $\sum_{s = g_n(k)n }^{k -1} \beta_{s}^{(k)}
\phi_{g_n(k+1)} ({\bf A})
 \widehat{\bf g}_{s}
$; \\
\>12.\>\>Else\\
\>13. \>\>\> $ \phi_{g_n(k)}({\bf A}) \widehat{\bf r}_k =
\phi_{g_n(k)}({\bf A}) \widehat{\bf r}_{k-1} - \alpha_k {\bf A}
\phi_{g_n(k)}({\bf A}) \widehat{\bf g}_{k-1}$;
\\
\>14. \>\>\>$ \phi_{g_n(k+1)}({\bf A}) \widehat{\bf r}_k =
(\rho_{g_n(k+1)} {\bf A} + {\bf I}) \phi_{g_n(k+1)-1} ({\bf
A})\widehat{\bf r}_k$;\\
\>15. \>\>\>For $s = g_n(k)n, \cdots, k - 1$ \\
\>16. \>\>\>\>$\beta^{(k)}_{s} = - {\bf q}^H_{r_n(s+1)} \left(
\phi_{g_n(s+1)+1} ({\bf A})
 \widehat{\bf r}_k + \right.$\\
\>\>\>\>\>\>\>$\left. \rho_{g_n(s+1)+1} \sum_{t = g_n(k)n }^{s-1}
\beta^{(k)}_t
  {\bf A} \phi_{g_n(s+1)} ({\bf
A}) \widehat{\bf g}_t
 \right)
\big/\rho_{g_n(s+1)+1}
 {\bf q}^H_{r_n(s+1)} {\bf A} \phi_{g_n(s+1)}({\bf A}) \widehat{\bf g}_s$; \\
\>17. \>\>\>End \\
\>18.\>\>\> $\phi_{g_n(k+1)} ({\bf A}) \widehat{\bf g}_k =
\phi_{g_n(k+1)} ({\bf A}) \widehat{\bf r}_k + (\rho_{g_n(k+1)}{\bf
A} + {\bf I})\sum_{s = g_n(k)n }^{k -1} \beta_{s}^{(k)}
\phi_{g_n(k+1)-1} ({\bf A}) \widehat{\bf g}_{s}$; \\
\>19.\>\>End\\
\>20. \> End
\end{tabbing}

\vspace{.2cm}

Lines 4, 11, 13 and 18, DS\#6, were obtained from Lines 4, 11, 13
and 17, DS\#5, by multiplying them with $\phi_{g_n(k)}({\bf A})$ and
$\phi_{g_n(k+1)}({\bf A})$ respectively. Line 14, DS\#6, is a direct
result of the definition (\ref{equ:12-23}) of $\phi$.

Now we use Proposition \ref{prop:1} to write DS\#6 as\\

{\bf Derivation Stage \#7.} \vspace{.2cm}
\begin{tabbing}
x\=xxx\= xxx\=xxx\=xxx\=xxx\=xxx\=xxx\=xxx\=xxx\kill
\>1. \>For $k = 1, 2, \cdots$, until convergence: \\
\>2. \>\>$\alpha_k = {\bf q}_{r_n(k)}^H \phi_{g_n(k)} ({\bf A})
\widehat{\bf r}_{k-1} / {\bf q}_{r_n(k)}^H
{\bf A} \phi_{g_n(k)}({\bf A}) \widehat{\bf g}_{k-1}$;\\
\>3.\>\> If $r_n(k) < n$\\
\>4. \>\>\> $\phi_{g_n(k+1)}({\bf A}) \widehat{\bf r}_k =
\phi_{g_n(k)}({\bf A}) \widehat{\bf r}_{k-1} - \alpha_k {\bf A}
\phi_{g_n(k)}({\bf A}) \widehat{\bf g}_{k-1}$;
\\
\>5. \>\>\>For $s = \max (k - n, 0), \cdots, g_n(k)n - 1$ \\
\>6. \>\>\>\>$\beta^{(k)}_{s} = - {\bf q}^H_{r_n(s+1)}
 \left(
  \phi_{g_n(k+1)}({\bf
A}) \widehat{\bf r}_k \right.$\\
\>\>\>\>\> $\left. + \,\rho_{g_n(k+1)} \sum_{t = \max (k - n, 0)
}^{s-1} \beta^{(k)}_t {\bf A} \phi_{g_n(t+1)}({\bf A}) \widehat{\bf
g}_t \right) \big/\rho_{g_n(k+1)}
 {\bf q}^H_{r_n(s+1)} {\bf A} \phi_{g_n(s+1)} ({\bf A}) \widehat{\bf g}_s$;\\
\>7. \>\>\>End
\\
\>8. \>\>\>For $s = g_n(k)n, \cdots, k - 1$ \\
\>9. \>\>\>\>$\beta^{(k)}_{s} = - {\bf q}^H_{r_n(s+1)} {\bf A}
 \left(
 \phi_{g_n(k+1)}({\bf A}) \widehat{\bf
r}_k + \sum_{t = \max (k - n, 0) }^{g_n(k)n-1} \beta^{(k)}_t
 \phi_{g_n(t+1)+1}
({\bf A})
 \widehat{\bf g}_t \right.$\\
\>\>\>\>\>\>\> $\left. + \sum_{t = g_n(k)n }^{s-1} \beta^{(k)}_t
 \phi_{g_n(t+1)}
({\bf A}) \widehat{\bf g}_t
 \right)
\big/ {\bf q}^H_{r_n(s+1)} {\bf A} \phi_{g_n(s+1)} ({\bf A}) \widehat{\bf g}_s$;\\
\>10. \>\>\>End
\\
\>11.\>\>\> $\phi_{g_n(k+1)} ({\bf A}) \widehat{\bf g}_k =
\phi_{g_n(k+1)} ({\bf A}) \widehat{\bf r}_k + (\rho_{g_n(k+1)} {\bf
A} + {\bf I}) \sum_{s = \max (k - n, 0) }^{g_n(k)n -1}
\beta_{s}^{(k)} \phi_{g_n(s+1)} ({\bf A})
            \widehat{\bf g}_{s}
+$\\
\>\>\>\>\>\>\>\>\> $\sum_{s = g_n(k)n }^{k -1} \beta_{s}^{(k)}
\phi_{g_n(s+1)} ({\bf A})
 \widehat{\bf g}_{s}
$; \\
\>12.\>\>Else\\
\>13. \>\>\> $ \phi_{g_n(k)}({\bf A}) \widehat{\bf r}_k =
\phi_{g_n(k)}({\bf A}) \widehat{\bf r}_{k-1} - \alpha_k {\bf A}
\phi_{g_n(k)}({\bf A}) \widehat{\bf g}_{k-1}$;
\\
\>14. \>\>\>$ \phi_{g_n(k+1)}({\bf A}) \widehat{\bf r}_k =
(\rho_{g_n(k+1)} {\bf A} + {\bf I}) \phi_{g_n(k)} ({\bf
A})\widehat{\bf r}_k$;\\
\>15. \>\>\>For $s = g_n(k)n, \cdots, k - 1$ \\
\>16. \>\>\>\>$\beta^{(k)}_{s} = - {\bf q}^H_{r_n(s+1)} \left(
\phi_{g_n(k+1)} ({\bf A})
 \widehat{\bf r}_k + \right.$\\
\>\>\>\>\>\>\>$\left. \rho_{g_n(k+1)} \sum_{t = g_n(k)n }^{s-1}
\beta^{(k)}_t
  {\bf A} \phi_{g_n(t+1)} ({\bf
A}) \widehat{\bf g}_t
 \right)
\big/\rho_{g_n(k+1)}
 {\bf q}^H_{r_n(s+1)} {\bf A} \phi_{g_n(s+1)}({\bf A}) \widehat{\bf g}_s$; \\
\>17. \>\>\>End \\
\>18.\>\>\> $\phi_{g_n(k+1)} ({\bf A}) \widehat{\bf g}_k =
\phi_{g_n(k+1)} ({\bf A}) \widehat{\bf r}_k + (\rho_{g_n(k+1)}{\bf
A} + {\bf I})\sum_{s = g_n(k)n }^{k -1} \beta_{s}^{(k)}
\phi_{g_n(s+1)} ({\bf A}) \widehat{\bf g}_{s}$; \\
\>19.\>\>End\\
\>20. \> End
\end{tabbing}

\vspace{.2cm} We remark that the term $\phi_{g_n(t+1)+1} ({\bf A})
 \widehat{\bf g}_t$ in the first sum on Line 9 can be further
 written as
\begin{equation} \label{equ:9-29}
\begin{array}{rl} \phi_{g_n(t+1)+1} ({\bf A})
 \widehat{\bf g}_t &= (\rho_{g_n(t+1)+1} {\bf A} + {\bf I})
 \phi_{g_n(t+1)}({\bf A}) \widehat{\bf g}_t\\
 & = (\rho_{g_n(k+1)} {\bf A} + {\bf I})
 \phi_{g_n(t+1)}({\bf A})
 \widehat{\bf g}_t.
\end{array}
\end{equation}

Substituting (\ref{equ:9-29}) and (\ref{equ:9-24-1}) into DS\#7 then
yields a set of updating relations of the vectors defined by
(\ref{equ:9-24-1}).\\

{\bf Derivation Stage \#8.} \vspace{.2cm}
\begin{tabbing}
x\=xxx\= xxx\=xxx\=xxx\=xxx\=xxx\=xxx\=xxx\=xxx\kill
\>1. \>For $k = 1, 2, \cdots$, until convergence: \\
\>2. \>\>$\alpha_k = {\bf q}_{r_n(k)}^H {\bf r}_{k-1} / {\bf
q}_{r_n(k)}^H
 {\bf w}_{k-1}$;\\
\>3.\>\> If $r_n(k) < n$\\
\>4. \>\>\> ${\bf r}_k = {\bf r}_{k-1} - \alpha_k {\bf w}_{k-1}$;
\\
\>5. \>\>\>For $s = \max (k - n, 0), \cdots, g_n(k)n - 1$ \\
\>6. \>\>\>\>$\beta^{(k)}_{s} = - {\bf q}^H_{r_n(s+1)}
 \left(
  {\bf r}_k + \,\rho_{g_n(k+1)} \sum_{t = \max (k - n, 0)
}^{s-1} \beta^{(k)}_t {\bf w}_t \right) \big/\rho_{g_n(k+1)}
 {\bf q}^H_{r_n(s+1)} {\bf w}_s$; \\
\>7. \>\>\>End
\\
\>8. \>\>\>For $s = g_n(k)n, \cdots, k - 1$ \\
\>9. \>\>\>\>$\beta^{(k)}_{s} = - {\bf q}^H_{r_n(s+1)} \left({\bf A}
 {\bf
r}_k + \sum_{t = \max (k - n, 0) }^{g_n(k)n-1} \beta^{(k)}_t
(\rho_{g_n(k+1)} {\bf A} + {\bf I})
 {\bf w}_t
  \right.$\\
  \>\>\>\>\>\>\> $\left. + \sum_{t = g_n(k)n }^{s-1} \beta^{(k)}_t
 {\bf w}_t
 \right)
\big/ {\bf q}^H_{r_n(s+1)} {\bf w}_s$; \\
\>10. \>\>\>End
\\
\>11.\>\>\> ${\bf g}_k = {\bf r}_k + \rho_{g_n(k+1)} \sum_{s = \max
(k - n, 0) }^{g_n(k)n -1} \beta_{s}^{(k)} {\bf w}_{s} + \sum_{s =
\max (k - n, 0) }^{g_n(k)n -1} \beta_{s}^{(k)} {\bf g}_{s} + \sum_{s
= g_n(k)n }^{k -1} \beta_{s}^{(k)} {\bf g}_{s}
$;\\
\>12.\>\>Else\\
\>13. \>\>\> $ {\bf u}_k = {\bf r}_{k-1} - \alpha_k {\bf w}_{k-1}$;
\\
\>14. \>\>\>${\bf r}_k = (\rho_{g_n(k+1)} {\bf A} + {\bf I})
{\bf u}_k$;\\
\>15. \>\>\>For $s = g_n(k)n, \cdots, k - 1$ \\
\>16. \>\>\>\>$\beta^{(k)}_{s} = - {\bf q}^H_{r_n(s+1)} \left( {\bf
r}_k +  \rho_{g_n(k+1)} \sum_{t = g_n(k)n }^{s-1} \beta^{(k)}_t
  {\bf w}_t
 \right)
\big/\rho_{g_n(k+1)}
 {\bf q}^H_{r_n(s+1)} {\bf w}_s$;\\
\>17. \>\>\>End \\
\>18.\>\>\> ${\bf g}_k = {\bf r}_k + \rho_{g_n(k+1)} \sum_{s =
g_n(k)n }^{k -1} \beta_{s}^{(k)} {\bf w}_{s} + \sum_{s = g_n(k)n
}^{k -1} \beta_{s}^{(k)} {\bf g}_{s}
$; \\
\>19.\>\>End\\
\>20. \> End
\end{tabbing}

\vspace{.2cm}

DS\#8 does not contain any update about ${\bf w}_k$. For the
updates, we multiply the equations on Lines 11 and 18 by ${\bf A}$
to get
\begin{equation} \label{equ:10-1-1}
\begin{array}{rl} {\bf w}_k =& {\bf A} ({\bf r}_k +
\rho_{g_n(k+1)} \sum_{s = \max (k - n, 0) }^{g_n(k)n -1}
\beta_{s}^{(k)} {\bf w}_{s}) + \sum_{s = \max (k - n, 0) }^{g_n(k)n
-1} \beta_{s}^{(k)} {\bf w}_{s} \\
&+ \sum_{s = g_n(k)n }^{k -1} \beta_{s}^{(k)} {\bf w}_{s}
\end{array}
\end{equation}
if $r_n(k) < n$, and
\begin{equation}
\label{equ:10-1-2}
\begin{array}{l} {\bf w}_k = {\bf A} ({\bf r}_k +
\rho_{g_n(k+1)} \sum_{s = g_n(k)n }^{k -1} \beta_{s}^{(k)} {\bf
w}_{s}) + \sum_{s = g_n(k)n }^{k -1} \beta_{s}^{(k)} {\bf w}_{s}
\end{array}
\end{equation}
if $r_n(k) = n$.

Again, we consider ${\bf r}_k$ to be a residual. To be consistent
with Lines 4, 13 and 14, we update the solution vector ${\bf x}_k$
as
\begin{equation}\label{equ:10-1-3}
\begin{array}{rl}
{\bf x}_k = & \left\{ \begin{array}{lcl} {\bf x}_{k-1} + \alpha_{k}
 {\bf g}_{k-1}, & &\mbox{if } r_n(k) < n \\
 -\rho_{g_n(k+1)} {\bf u}_k + {\bf x}_{k-1} + \alpha_k {\bf g}_{k-1}, & & \mbox{if } r_n(k) = n.
\end{array} \right.
\end{array}
\end{equation}
Now adding (\ref{equ:10-1-1}), (\ref{equ:10-1-2}) and
(\ref{equ:10-1-3}) to DS\#8 and simplifying the operations
appropriately, we then arrive at the following algorithm. The free
parameter $\rho_{g_n(k+1)}$ is chosen to minimize the $2$-norm of
${\bf r}_k$.
\\

\begin{algorithm}{\bf ML($n$)BiCGStab without preconditioning
associated with definition (\ref{equ:9-24-1}) } \label{alg:3}
\vspace{.2cm}
\begin{tabbing}
x\=xxx\= xxx\= xxx\= xxx\= xxx\= xxx\= xxx\=
xxx\=xxx\=xxx\=xxx\=xxx\kill \>1. \> Choose an initial guess ${\bf
x}_0$ and $n$ vectors ${\bf q}_1, {\bf q}_2, \cdots,
{\bf q}_n$. \\
\>2. \>  Compute ${\bf r}_0 = {\bf b} - {\bf A} {\bf x}_0$ and ${\bf
g}_0 = {\bf r}_0,\,\, {\bf w}_0 = {\bf A}
{\bf g}_0,\,\, c_0 = {\bf q}^H_{1} {\bf w}_0$. \\
\>3. \>  For $k = 1, 2, \cdots$, until convergence: \\
\>4. \>\>
           $  \alpha_k = {\bf q}^H_{r_n(k)}
{\bf r}_{k-1} /
c_{k-1}$; \\
\>5. \>\> If $r_n(k) < n$\\
\>6. \>\>\>  $ {\bf x}_k = {\bf x}_{k-1} + \alpha_k
{\bf g}_{k-1}$; \\
\>7. \>\>\>  $ {\bf r}_k = {\bf r}_{k-1} - \alpha_k
{\bf w}_{k-1}$; \\
\>8. \>\>\>  ${\bf z}_w = {\bf r}_k,\,\, {\bf g}_k = {\bf 0}$;\\
\>9. \>\>\>   For $s = \max (k - n, 0), \cdots, g_n(k) n - 1$ \\
\>10. \>\>\>\>        $\tilde{\beta}^{(k)}_{s} = - {\bf q}^H_{r_n(s+1)}
 {\bf z}_w \big/
c_s$; \,\,\,\,\,\,\,\,\, \% $\tilde{\beta}^{(k)}_{s} = \beta^{(k)}_{s} \rho_{g_n(k+1)}$\\
\>11. \>\>\>\>${\bf z}_w = {\bf z}_w + \tilde{\beta}^{(k)}_s
{\bf w}_s$;\\
\>12. \>\>\>\>${\bf g}_k = {\bf g}_k + \tilde{\beta}^{(k)}_s
{\bf g}_s$;\\
\>13. \>\>\>            End \\
\>14. \>\>\>$\displaystyle{{\bf g}_k = {\bf z}_w + \frac{1}{\rho_{g_n(k+1)}} {\bf g}_k}$;\\
\>15. \>\>\>${\bf w}_k = {\bf A} {\bf g}_k$;\\
\>16. \>\>\>   For $s = g_n(k) n, \cdots, k - 1$ \\
\>17. \>\>\>\>        $\beta^{(k)}_{s} = - {\bf q}^H_{r_n(s+1)} {\bf
w}_k
 \big/
c_s$; \\
\>18. \>\>\>\>${\bf w}_k = {\bf w}_k + \beta^{(k)}_s {\bf w}_s$;\\
\>19. \>\>\>\>${\bf g}_k = {\bf g}_k + \beta^{(k)}_s
{\bf g}_s$;\\
\>20. \>\>\>            End \\
\>21. \>\> Else \\
\>22. \>\>\>  ${\bf x}_k = {\bf x}_{k-1} + \alpha_k
{\bf g}_{k-1}$;\\
\>23. \>\>\>  $ {\bf u}_k = {\bf r}_{k-1} - \alpha_k
{\bf w}_{k-1}$; \\
\>24. \>\>\> $\rho_{g_n(k+1)} = - ({\bf A} {\bf u}_k)^H {\bf u}_k /
\| {\bf A} {\bf u}_k \|^2_2$;\\
\>25. \>\>\> ${\bf x}_k = {\bf x}_k - \rho_{g_n(k+1)}
{\bf u}_k$; \\
\>26. \>\>\> ${\bf r}_k = \rho_{g_n(k+1)} {\bf A} {\bf u}_k +
{\bf u}_k$; \\
\>27. \>\>\> ${\bf z}_w = {\bf r}_k,\,\, {\bf g}_k = {\bf 0}$;\\
\>28. \>\>\>   For $s = g_n(k) n, \cdots, k - 1$ \\
\>29. \>\>\>\>        $\tilde{\beta}^{(k)}_{s} = - {\bf q}^H_{r_n(s+1)} {\bf
z}_w \big/
c_s$; \,\,\,\,\,\,\,\,\, \% $\tilde{\beta}^{(k)}_{s} = \beta^{(k)}_{s} \rho_{g_n(k+1)}$\\
\>30. \>\>\>\>${\bf z}_w = {\bf z}_w +
\tilde{\beta}^{(k)}_s {\bf w}_s$; \\
\>31. \>\>\>\>${\bf g}_k = {\bf g}_k + \tilde{\beta}^{(k)}_s
{\bf g}_s$; \\
\>32. \>\>\>            End \\
\>33. \>\>\>$\displaystyle{{\bf g}_k =
 {\bf z}_w + \frac{1}{\rho_{g_n(k+1)}}{\bf g}_k}$;\\
\>34. \>\>\>${\bf w}_k = {\bf A} {\bf g}_k
$;\\
\>35. \>\> End\\
\>36. \>\>$c_k = {\bf q}^H_{r_n(k+1)} {\bf w}_k$;\\
\>37. \> End
\end{tabbing}
\end{algorithm}

\vspace{.2cm}

We remark that (i) the algorithm does not compute ${\bf u}_k$ when
$r_n(k) < n$. In fact, ${\bf u}_k = {\bf r}_k$ when $r_n(k) < n$
(see the remark right after (\ref{equ:9-24-1})); (ii) if the ${\bf
u}_k$ on Line 23 happens to be zero, then the ${\bf x}_k$ on Line 22
will be the exact solution to system (\ref{equ:7-9-3}) and the
algorithm stops there.

\begin{table}[tbp]
\footnotesize \caption{Average cost per ($k$-)iteration step of
Algorithm \ref{alg:10-27} and its storage requirement.}
\begin{center}
\begin{tabular}{|c|c|c|c|}  \hline
Preconditioning (${\bf M}^{-1} {\bf v}$) & $\displaystyle{1 +
\frac{1}{n}}$  & Vector addition (${\bf u} \pm {\bf v}$) &
$\displaystyle{1}$ \\ \hline Matvec (${\bf A} {\bf v}$) &
$\displaystyle{1 + \frac{1}{n}}$  & Saxpy (${\bf u} + \alpha {\bf
v}$) & $\displaystyle{2 n + 2 + \frac{2}{n}}$  \\ \hline dot product
($\displaystyle{ {\bf u}^H {\bf v}}$) & $\displaystyle{n + 1+
\frac{2}{n}}$ & Storage &${\bf A} + {\bf M} + (3 n+5) N +O(n)$
\\ \hline
\end{tabular}
\end{center} \label{tab:10-29-1}
\end{table}

The cost and storage requirement, obtained from its preconditioned
version, Algorithm \ref{alg:10-27} in \S\ref{sec:apen}, are listed
in Table \ref{tab:10-29-1}. Compared to Algorithm \ref{alg:2},
Algorithm \ref{alg:3} saves about $20\%$ in saxpy.
Since only three sets of vectors $\{ {\bf q}_1, \ldots, {\bf
q}_n\}$, $\{{\bf g}_{k-n}, \ldots, {\bf g}_{k-1}\}$ and $\{{\bf
w}_{k-n}, \ldots, {\bf w}_{k-1}\}$ are needed in iteration $k$, the
storage is about $3 n N$ besides storing ${\bf A}$ and ${\bf M}$.


\subsection{Properties}
We summarize the properties about Algorithm \ref{alg:3} below. Their
proofs are similar to those in Proposition \ref{prop:10-1}. Since
${\bf r}_0 = \widehat{\bf r}_0$ by (\ref{equ:10-1-5}), $\nu$ is also
the degree of the minimal polynomial of ${\bf
r}_0$ with respect to ${\bf A}$.\\

\begin{proposition} Under the assumptions of Proposition
\ref{prop:7-11-1}, if $\rho_{g_n(k+1)} \ne 0$ and $-
1/\rho_{g_n(k+1)} \not\in \sigma({\bf A})$ for $1 \leq k \leq
\nu-1$, 
then
Algorithm \ref{alg:3} does not break down by zero division for $k =
1, 2, \cdots, \nu$, and the approximate solution ${\bf x}_\nu$ at
step $k = \nu$ is exact to the system (\ref{equ:7-9-3}). Moreover,
the computed quantities satisfy
\begin{enumerate}
\item[(a)]
${\bf x}_k \in {\bf x}_0 + span \{ {\bf r}_0, {\bf A} {\bf r}_0,
\ldots, {\bf A}^{g_n(k+1)+k-1} {\bf r}_0 \}$ and ${\bf r}_k = {\bf
b} - {\bf A} {\bf x}_k \in {\bf r}_0 + span \{ {\bf A} {\bf r}_0,
{\bf A}^2 {\bf r}_0, \ldots, {\bf A}^{g_n(k+1)+k} {\bf r}_0 \}$ for
$1 \leq k \leq \nu-1$.

\item[(b)] ${\bf r}_k \ne {\bf 0}$ for $1 \leq k \leq \nu -1$; ${\bf
r}_\nu = {\bf 0}$.

\item[(c)] ${\bf r}_k \perp span \{ {\bf q}_1, {\bf q}_2, \cdots,
{\bf q}_{r_n(k)}\}$ and ${\bf r}_k \not\perp {\bf q}_{r_n(k)+1}$ for
$1 \leq k \leq \nu -1$ with $r_n(k) < n$; ${\bf r}_k \not\perp {\bf
q}_{1}$ for $1 \leq k \leq \nu -1$ with $r_n(k) = n$.

\item[(d)] ${\bf u}_k \perp span \{ {\bf q}_1, {\bf q}_2, \cdots,
{\bf q}_{n}\}$ for $1 \leq k \leq \nu$ with $r_n(k) = n$.

\item[(e)] ${\bf A}{\bf g}_k \perp span \{ {\bf q}_1, {\bf q}_2, \cdots,
{\bf q}_{r_n(k)}\}$ and ${\bf A}{\bf g}_k
 \not\perp {\bf q}_{r_n(k)+1}$ for
$1 \leq k \leq \nu -1$ with $r_n(k) < n$; ${\bf A} {\bf g}_k
\not\perp {\bf q}_{1}$ for $1 \leq k \leq \nu -1$ with $r_n(k) = n$.
\end{enumerate}
\end{proposition}

\vspace{.2cm}

\section{Relations to Some Other Methods} \label{sec:12-24-4} In this section, we discuss the
relations of ML($n$)BiCGStab with the FOM, BiCGStab and IDR($n$)
methods under the exact arithmetic environment.

\subsection{Algorithm \ref{alg:2}}\label{sec:10-9}
\begin{enumerate}
\item {\it Relation with FOM\cite{saad2}}. Consider the case where $n \geq \nu$. In this case, $g_n(k) = 0$ and
$r_n(k) = k$ for $k = 1, 2, \cdots, \nu$. Hence ${\bf p}_k = {\bf
q}_k$ by (\ref{equ:7-9-5}). If we choose ${\bf q}_k = \widehat{\bf
r}_{k-1}$ in Algorithm \ref{alg:1} (it is possible since
$\widehat{\bf r}_{k-1}$ is computed before ${\bf q}_k$ is used),
then the $\widehat{\bf x}_k$ and $\widehat{\bf r}_k$ computed by the
algorithm satisfy
\begin{equation} \label{equ:10-3-10}
\left\{ \begin{array}{l} \widehat{\bf x}_k \in \widehat{\bf x}_0 +
span \{ \widehat{\bf r}_0, {\bf A} \widehat{\bf r}_0, \ldots, {\bf
A}^{k-1} \widehat{\bf r}_0 \}, \\
\widehat{\bf r}_k \perp span \{
\widehat{\bf r}_0, \widehat{\bf r}_1, \ldots, \widehat{\bf
r}_{k-1}\}
\end{array} \right.
\end{equation}
for $1 \leq k \leq \nu$ by Proposition \ref{prop:7-11-1}(a), (d).
(\ref{equ:10-3-10}) is what the FOM approximate solution ${\bf
x}_k^{FOM}$ needs to satisfy. Therefore, when $n \geq \nu$ and with
the choice ${\bf q}_k = \widehat{\bf r}_{k-1}$, Algorithm
\ref{alg:1} is mathematically equivalent to FOM. \\

Now, from (\ref{equ:7-29}), the ${\bf r}_k$ computed by Algorithm
\ref{alg:2} satisfies
$$
{\bf r}_k = \phi_{g_n(k)+1}({\bf A}) \,\widehat{\bf r}_k =
\phi_{1}({\bf A}) \,\widehat{\bf r}_k = (\rho_1 {\bf A} + {\bf I})
\,\widehat{\bf r}_k.
$$
Note that ${\bf u}_k = \phi_{g_n(k)}({\bf A}) \widehat{\bf r}_k =
\phi_{0}({\bf A}) \widehat{\bf r}_k = \widehat{\bf r}_k$. Thus, for
$1 \leq k \leq \nu$, ${\bf r}_k$ is the factor $\rho_1 {\bf A} +
{\bf I}$ times the FOM residual ${\bf u}_k$ if we set ${\bf q}_1 =
{\bf r}_0$ and ${\bf q}_{k+1} = {\bf u}_k$ in Algorithm
\ref{alg:2}.\footnote{In \cite{yeungchan}, a remark immediately
following Theorem 4.1 states that, when $n \geq \nu$ and with the
choice that ${\bf q}_1 = \phi_1({\bf A}^H) \phi_1({\bf A}) {\bf
r}_0$ and ${\bf q}_k = \phi_1({\bf A}^H) {\bf r}_{k-1}$ for $k \geq
2$, the ${\bf x}_k$ and ${\bf r}_k$ computed by Algorithm 2 (which
is mathematically equivalent to Algorithm \ref{alg:2} of this paper)
will satisfies (\ref{equ:10-3-10}) and therefore Algorithm 2 is a
FOM. The argument 
there about this remark is not correct.
The author remembers that the referees of \cite{yeungchan} were skeptical about the argument.}\\

\item \underline{Relation with BiCGStab\cite{van}}. When $n = 1$, we have $g_n(k) =
k-1$ and $r_n(k) = 1$ for $k \in {\mathcal N}$. Hence ${\bf p}_{k} =
\left( {\bf A}^H \right)^{k-1} {\bf q}_{1}$ by (\ref{equ:7-9-5}). By
Proposition \ref{prop:7-11-1}(a) and (d), the $\widehat{\bf x}_k$
and $\widehat{\bf r}_k$ computed by Algorithm \ref{alg:1} satisfy
\begin{equation} \label{equ:10-5}
\left\{
\begin{array}{l}
\widehat{\bf x}_k \in \widehat{\bf x}_0 + span \{ \widehat{\bf r}_0,
{\bf A} \widehat{\bf r}_0, \ldots, {\bf A}^{k-1} \widehat{\bf r}_0
\}\\
\widehat{\bf r}_k \perp span \{ {\bf q}_1, {\bf A}^H {\bf q}_1,
\ldots, \left( {\bf A}^H \right)^{k-1} {\bf q}_{1}\}
\end{array}\right.
\end{equation}
for $1 \leq k \leq \nu$. (\ref{equ:10-5}) is what the BiCG
approximate solution ${\bf x}_k^{BiCG}$ needs to satisfy. Therefore,
when $n = 1$, Algorithm
\ref{alg:1} is mathematically equivalent to BiCG.\\

Now, from (\ref{equ:7-29}), the ${\bf r}_k$ computed by Algorithm
\ref{alg:2} satisfies
$$
{\bf r}_k = \phi_{g_n(k)+1}({\bf A}) \,\widehat{\bf r}_k =
\phi_{k}({\bf A}) \,\widehat{\bf r}_k
$$
which is the definition of the BiCGStab residuals.  Thus
Algorithm \ref{alg:2} is mathematically equivalent to BiCGStab when
$n = 1$.\\

\item \underline{Relation with IDR($n$)\cite{sonn1}}. Write $k = j n + i$ as in
(\ref{equ:8-1}) with $1 \leq i \leq n, 0 \leq j$. Let ${\mathcal G}_0 =
K({\bf A}, {\bf r}_0)$ be the complete Krylov space and let ${\mathcal
S} = span \{{\bf q}_1, {\bf q}_2, \cdots, {\bf q}_n\}^\perp$. Define the Sonneveld spaces
$$
{\mathcal G}_{j+1} = (\rho_{j+1} {\bf A} + {\bf I}) ({\mathcal G}_j \cap
{\mathcal S}) = (\rho_{g_n(k)+1} {\bf A} + {\bf I}) ({\mathcal G}_j \cap
{\mathcal S})
$$
for $j = 0, 1, 2, \cdots$. By (\ref{equ:7-29}), we have
$$
{\bf r}_{jn+i} = \phi_{j+1}({\bf A}) \widehat{\bf r}_{jn+i} =
(\rho_{j+1} {\bf A} + {\bf I}) \phi_{j}({\bf A})\widehat{\bf
r}_{jn+i} = (\rho_{j+1} {\bf A} + {\bf I}) {\bf u}_{jn+i}.
$$
From Proposition \ref{prop:10-1}(d), ${\bf u}_{jn+i} \not\perp {\bf
q}_{i+1}$ if $i < n$. Hence ${\bf u}_{jn+i} \not\in {\mathcal G}_j \cap
{\mathcal S}$ and therefore ${\bf r}_{jn+i} \not\in {\mathcal G}_{j+1}$ when
$i < n$. From this point of view, Algorithm \ref{alg:2} is not a
IDR($n$) algorithm.
\end{enumerate}

\subsection{Algorithm \ref{alg:3}}

\begin{enumerate}
\item \underline{Relation with FOM}.
When $n \geq \nu$, $g_n(k) = 0$ and $r_n(k) = k$ for $1 \le k \leq
\nu$ and Algorithm \ref{alg:1}, with the choice ${\bf q}_k =
\widehat{\bf r}_{k-1}$, is a FOM algorithm as seen in
\S\ref{sec:10-9}. Now, from (\ref{equ:9-24-1}), the ${\bf r}_k$
computed by Algorithm \ref{alg:3} satisfies
$$
{\bf r}_k = \phi_{g_n(k+1)}({\bf A}) \,\widehat{\bf r}_k =
\phi_0({\bf A}) \,\widehat{\bf r}_k = \widehat{\bf r}_k.
$$
Thus Algorithm \ref{alg:3} is a FOM algorithm when we set ${\bf q}_k
= {\bf r}_{k-1}$.\\

\item \underline{Relation with BiCGStab}. When $n = 1$, we have $g_n(k) =
k-1$ and $r_n(k) = 1$ for $k \in {\mathcal N}$ and Algorithm \ref{alg:1}
is a BiCG algorithm. Now, from (\ref{equ:9-24-1}), the ${\bf r}_k$
computed by Algorithm \ref{alg:3} satisfies
$$
{\bf r}_k = \phi_{g_n(k+1)}({\bf A}) \,\widehat{\bf r}_k =
\phi_{k}({\bf A}) \,\widehat{\bf r}_k
$$
which is the definition of the BiCGStab residuals.  Thus
Algorithm \ref{alg:3} is mathematically equivalent to BiCGStab.\\

\item \underline{Relation with IDR($n$)}. Write $k = j n + i$ as in
(\ref{equ:8-1}) with $1 \leq i \leq n, 0 \leq j$. By
(\ref{equ:9-24-1}), we have
\begin{equation}\label{equ:11-23-10}
{\bf r}_{jn+i} = \phi_{g_n(jn+i+1)}({\bf A}) \widehat{\bf r}_{jn+i}
= \left\{ \begin{array}{lcl} \phi_{j}({\bf A})
\widehat{\bf r}_{jn+i} & & \mbox{ if } 1 \leq i < n,\\
\phi_{j+1}({\bf A}) \widehat{\bf r}_{jn+i} & & \mbox{ if } i = n.
\end{array}\right.
\end{equation}
By (\ref{equ:7-9-5}) and
Proposition \ref{prop:7-11-1}(d), we have
$$\left\{ \begin{array}{lcl}
\phi_{t}({\bf A}) \widehat{\bf r}_{jn+i} \in {\mathcal S} & &\mbox{ if }
1 \leq i < n \mbox{ and } 0 \leq t < j, \\
\phi_{t}({\bf A}) \widehat{\bf r}_{jn+n} \in {\mathcal S} & &\mbox{ if }
0 \leq t \leq j.
\end{array}\right.
$$
Thus, by induction on $t$, we have
$$
\left\{ \begin{array}{lcl} \phi_{t}({\bf A}) \widehat{\bf r}_{jn+i}
\in {\mathcal G}_t \cap {\mathcal S} & &\mbox{ if }
1 \leq i < n \mbox{ and } 0 \leq t < j, \\
\phi_{t}({\bf A}) \widehat{\bf r}_{jn+n} \in {\mathcal G}_t \cap {\mathcal
S} & &\mbox{ if } 0 \leq t \leq j.
\end{array}\right.
$$ 
Therefore, by (\ref{equ:11-23-10}),
$$
\left\{ \begin{array}{ll} {\bf r}_{jn+i} = \widehat{\bf r}_i \in
{\mathcal G}_0& \mbox{ if }
1 \leq i < n, j = 0,\\
{\bf r}_{jn+i} = (\rho_j {\bf A} + {\bf I}) \phi_{j-1}({\bf A})
\widehat{\bf r}_{jn+i} \in {\mathcal G}_j & \mbox{ if }
1 \leq i < n, 1 \leq j,\\
{\bf r}_{jn+n} = (\rho_{j+1}{\bf A} + {\bf I}) \phi_j({\bf A})
\widehat{\bf r}_{jn+n} \in {\mathcal G}_{j+1}.
\end{array}\right.
$$
So, the residuals in (\ref{equ:11-23-10}) lie in the Sonneveld spaces ${\mathcal G}$ and therefore Algorithm \ref{alg:3} is a IDR($n$) algortithm.
\end{enumerate}

\section{Implementation Issues}\label{sec:12-19}
A preconditioned ML($n$)BiCGStab algorithm can be obtained by
applying either Algorithm \ref{alg:2} or Algorithm \ref{alg:3} to
the system
$$
{\bf A} {\bf M}^{-1} {\bf y} = {\bf b}
$$
where ${\bf M}$ is nonsingular, then recovering ${\bf x}$ through
${\bf x} = {\bf M}^{-1} {\bf y}$. The resulting algorithms,
Algorithm \ref{alg:10-22} and Algorithm \ref{alg:10-27}, together
with their Matlab codes are presented in \S\ref{sec:11-19-1} and
\S\ref{sec:11-19-2} respectively. To avoid calling the index
functions $r_n(k)$ and $g_n(k)$ every $k$-iteration, we have split
the $k$-loop into a $i$-loop and a $j$-loop where $i, j, k$ are
related by (\ref{equ:8-1}) with $1 \leq i \leq n, 0 \leq j$.
Moreover, we have optimized the operations as possible as we can in
the resulting preconditioned algorithms.

Since we have compared ML($n$)BiCGStab with some existing methods in
\cite{yeungchan}, we will only concentrate on the performance of
ML($n$)BiCGStab itself. The following test data were downloaded from
Matrix Market.\footnote{ http://math.nist.gov/MatrixMarket/data/}\\

\begin{enumerate}
\item {\it e20r0100}, DRIVCAV Fluid Dynamics.  {\it e20r0100}
contains a $4241 \times 4241$ real unsymmetric matrix $\bf A$ with
$131,556$ nonzero entries and a real right-hand side $\bf b$.

\item {\it qc2534}, H2PLUS Quantum Chemistry, NEP Collection.
{\it qc2534} contains a $2534 \times 2534$ complex symmetric
indefinite matrix with $463,360$ nonzero entries, but does not
provide the right-hand side ${\bf b}$. Following \cite{ssg},  we set
${\bf b} = {\bf A} {\bf 1}$ with ${\bf 1} = [1, ,1, \cdots, 1]^T$.

\item {\it utm5940}, TOKAMAK Nuclear Physics (Plasmas). {\it utm5940}
contains a $5940 \times 5940$ real unsymmetric matrix ${\bf A}$ with
$83,842$ nonzero entries and a real right-hand side $\bf b$.\\
\end{enumerate}

Experiments were performed in Matlab Version 7.1 on a Windows XP
machine with a Pentium 4 processor. $ILU(0)$ preconditioner (p.294,
\cite{saad}) has been used in all the experiments. For {\it
e20r0100}, the $U$-factor of the $ILU(0)$ decomposition of $\bf A$
has some zeros along its main diagonal. In that experiment, we
replaced those zeros with $1$ so that the $U$-factor was invertible.

In all the experiments, initial guess ${\bf x}_0 = {\bf 0}$ and
stopping criterion is
$$
\| {\bf r}_k \|_2 / \| {\bf b}\|_2 < 10^{-7}
$$
where ${\bf r}_k$ is the computed residual. Except where specified,
auxiliary vectors ${\bf Q} \equiv [{\bf q}_1, {\bf q}_2, \cdots,
{\bf q}_n]$ are chosen to be ${\bf Q} = [{\bf r}_0, randn(N,n-1)]$
for {\it e20r0100} and {\it utm5940} and ${\bf Q} = [{\bf r}_0,
randn(N,n-1) + sqrt(-1)* randn(N,n-1)]$ for {\it qc2534}.

Moreover, for the convenience of our presentation, we introduce the
following functions:

\begin{enumerate}
\item[(a)]
$T_{conv}(n)$ is the time that a ML($n$)BiCGStab algorithm takes to
converge.

\item[(b)] $I_{conv}(n)$ is the number of iterations that
a ML($n$)BiCGStab algorithm takes to converge.


\item[(c)] $E(n) := \| {\bf b} - {\bf
A} {\bf x}\|_2 / \|{\bf b}\|_2$ is the true relative error of ${\bf
x}$ where ${\bf x}$ is the computed solution output by a
ML($n$)BiCGStab algorithm when it converges.
\end{enumerate}

\subsection{Stability} 
\begin{figure}[htbp]
\centerline{\hbox{ \psfig{figure=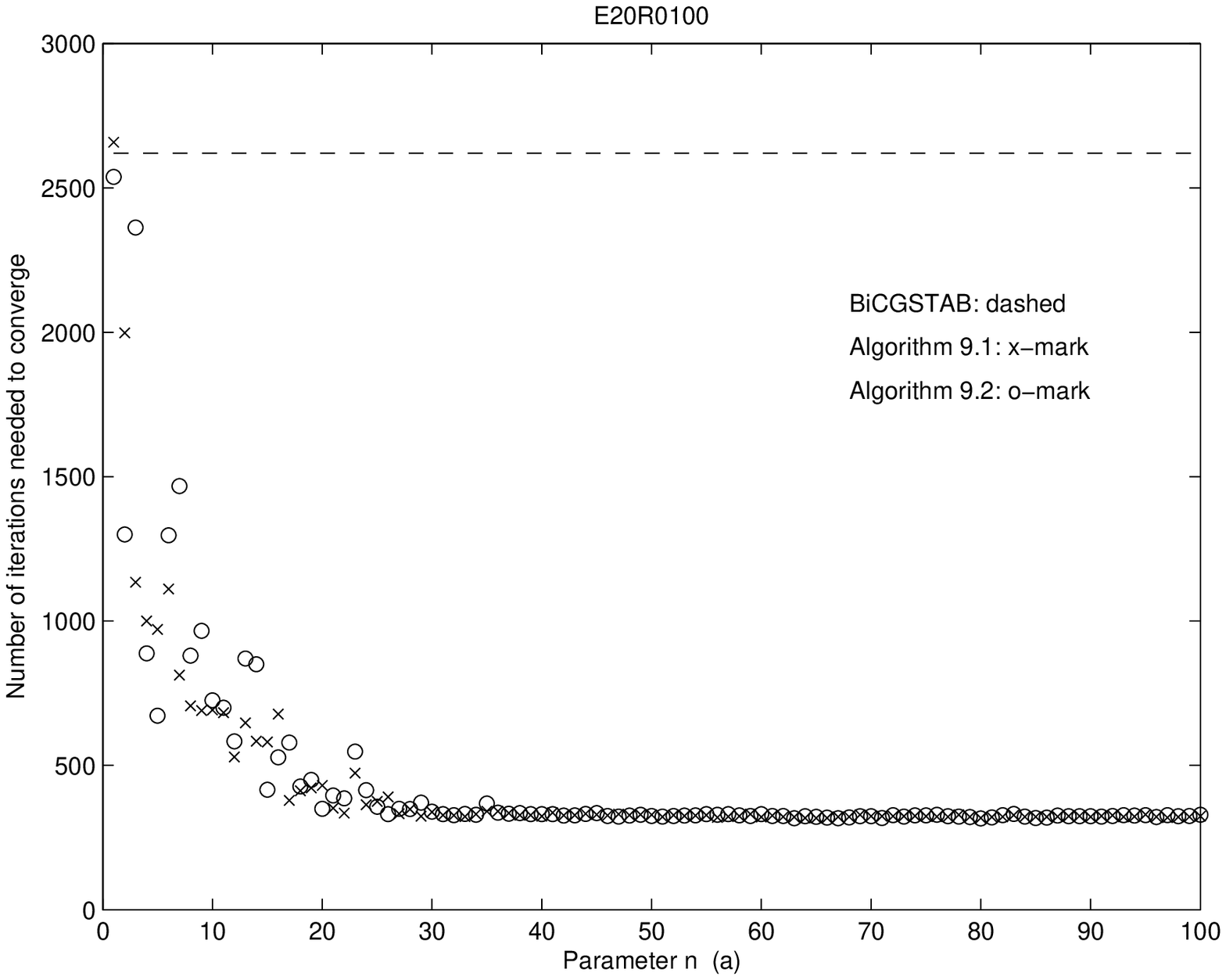,height=6cm }
\psfig{figure=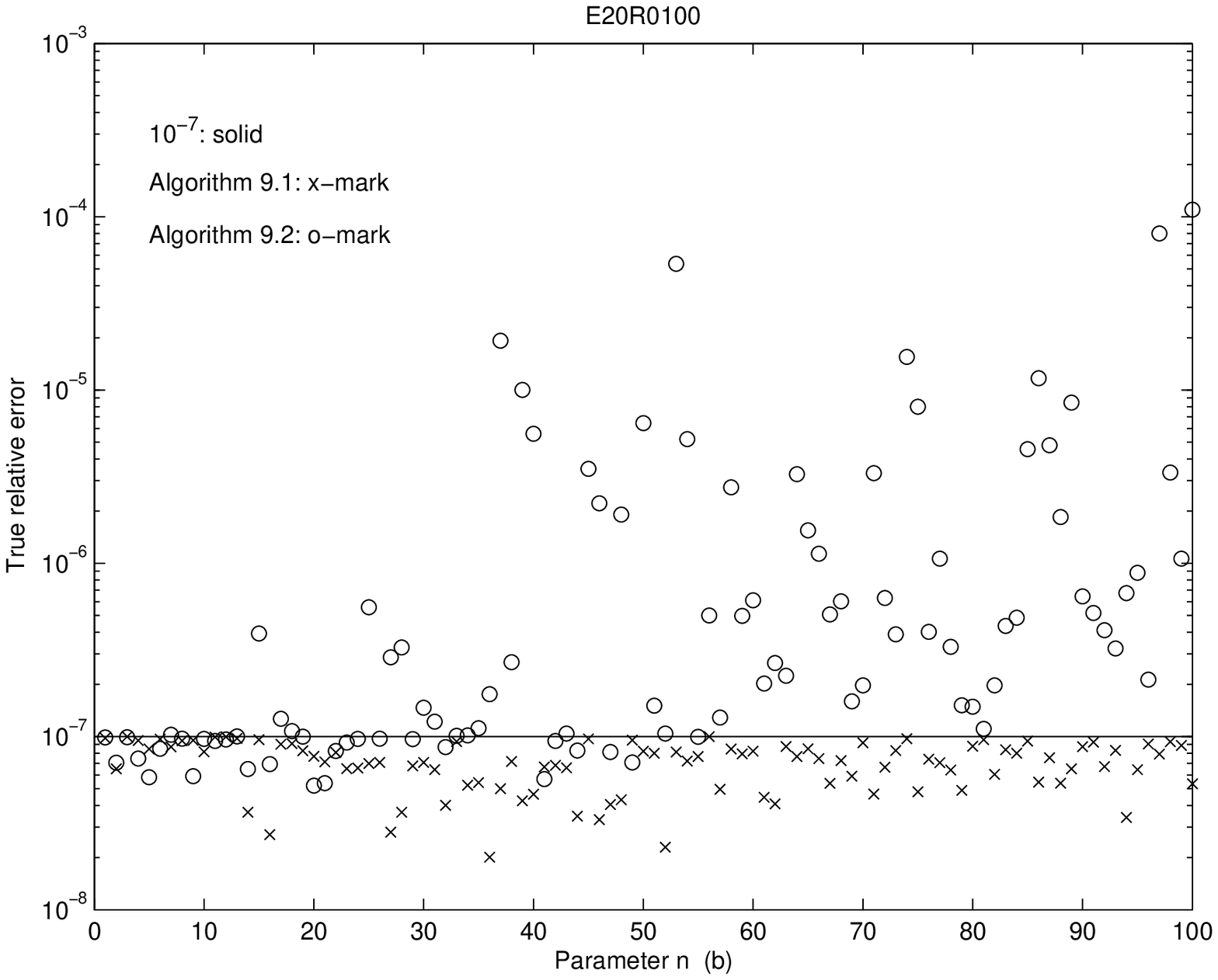,height=6cm}
 }} \caption{e20r0100: (a) Graphs of $I_{conv}(n)$ against $n$. BiCGStab took $2620$ iterations/$5240$
 MVs to converge. Full GMRES
 converged with $308$ MVs.
 (b)
Graphs of $E(n)$ against $n$. 
}
\label{fig:11-24-1}
\end{figure}

\begin{figure}[htbp]
\centerline{\hbox{ \psfig{figure=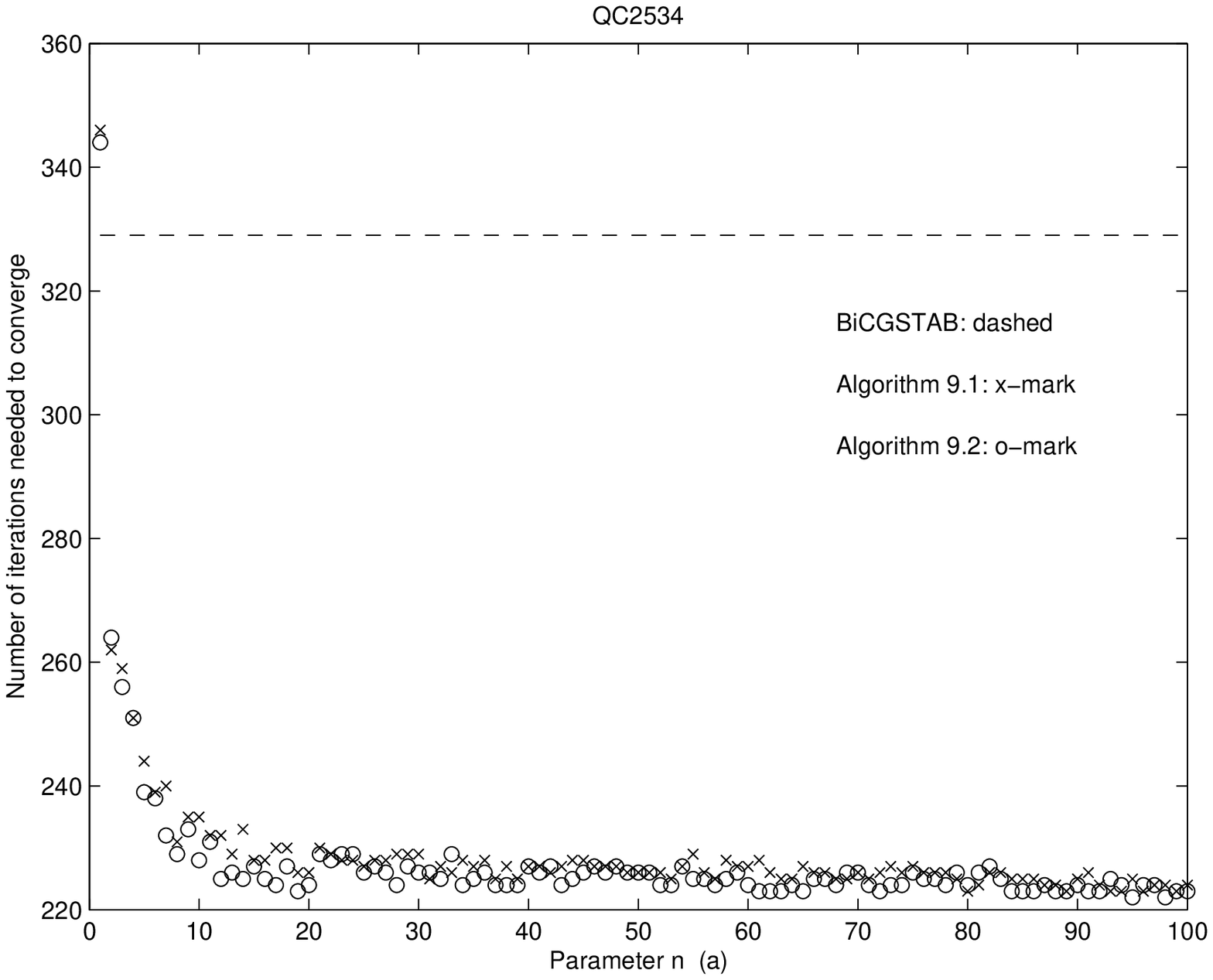,height=6cm }
\psfig{figure=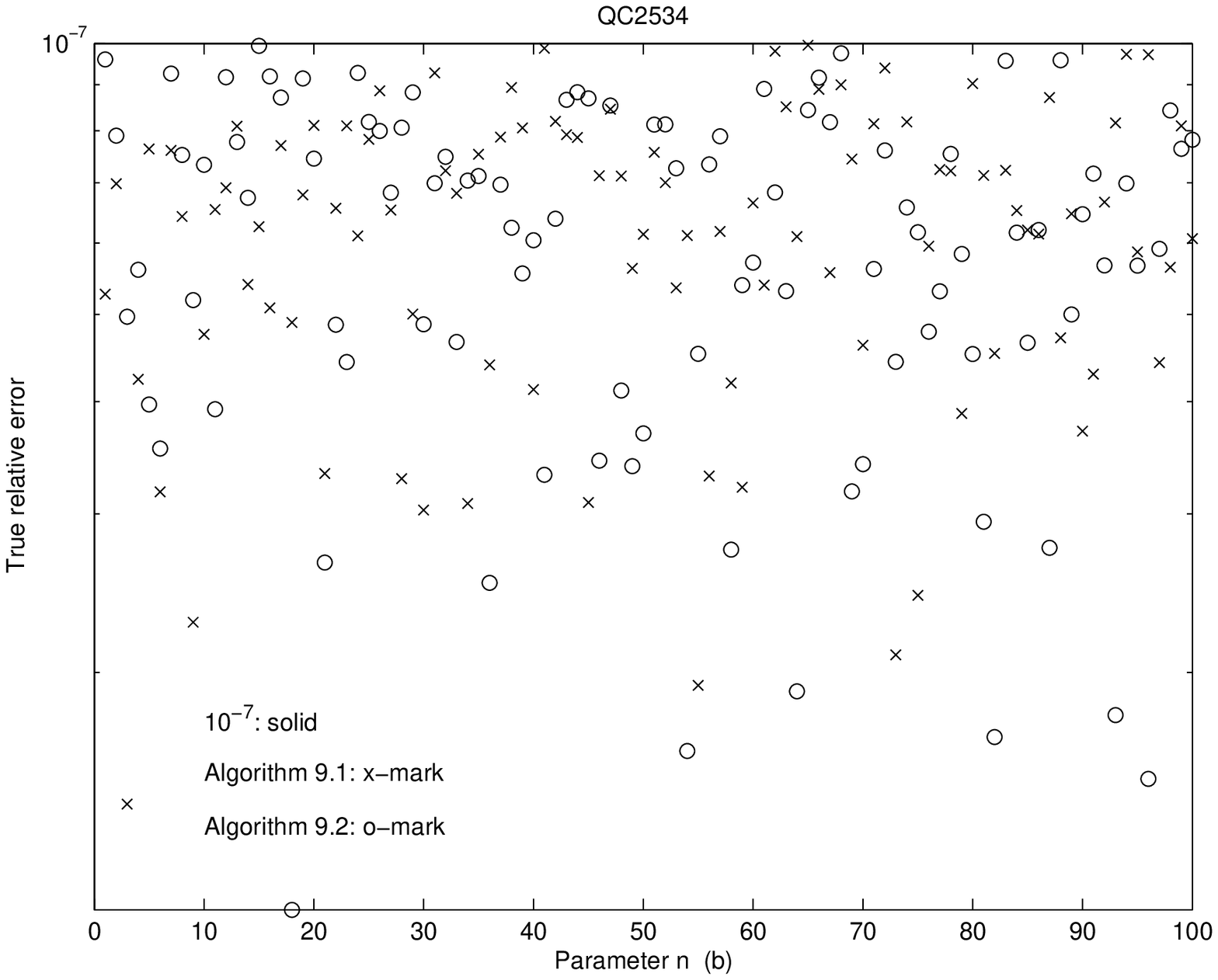,height=6cm}
 }} \caption{qc2534: (a) Graphs of $I_{conv}(n)$ against $n$. BiCGStab took $329$ iterations/$658$
 MVs to converge. Full GMRES
 converged with $439$ MVs.
 (b)
Graphs of $E(n)$ against $n$.
} \label{fig:11-24-2}
\end{figure}

\begin{figure}[htbp]
\centerline{\hbox{ \psfig{figure=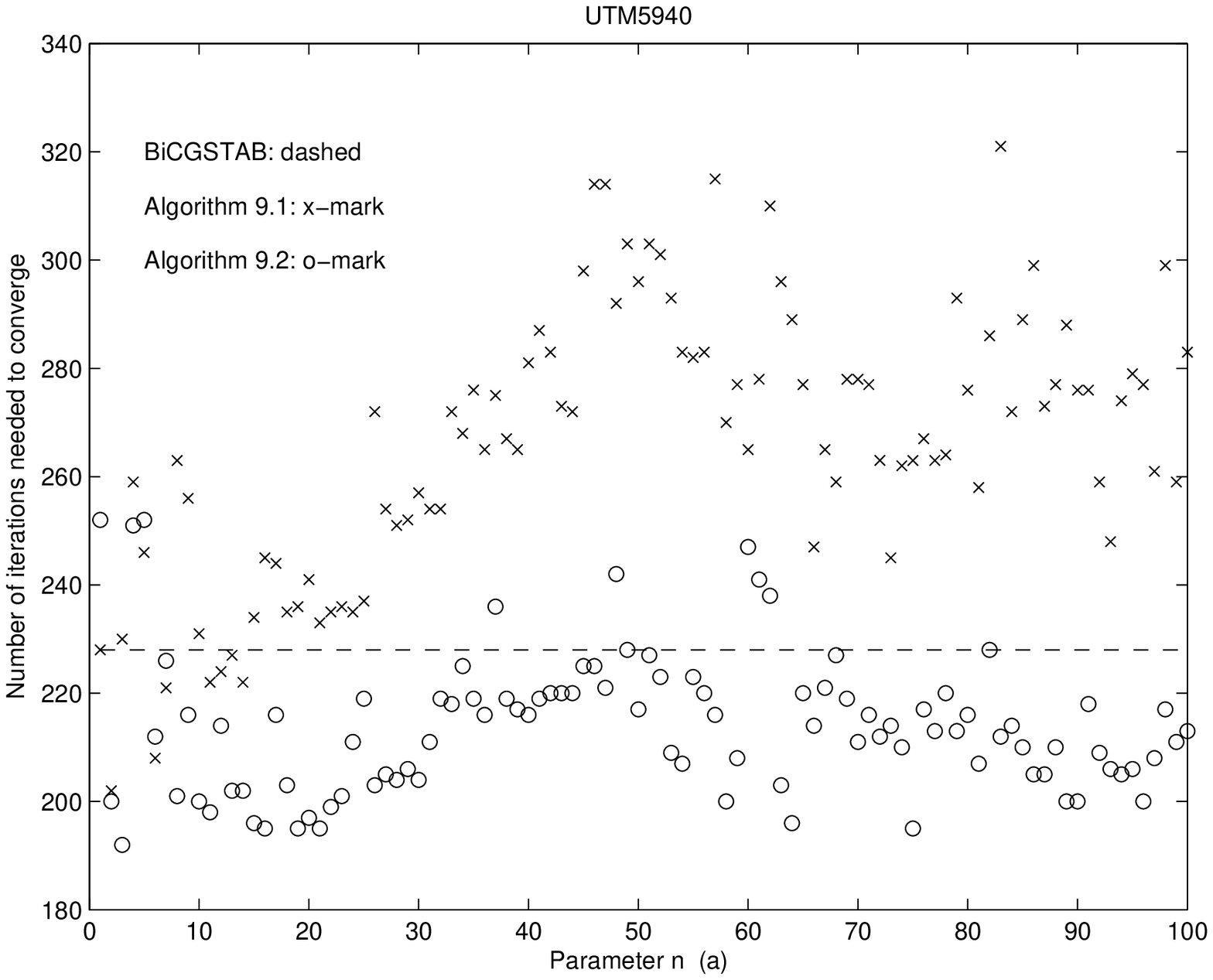,height=6cm }
\psfig{figure=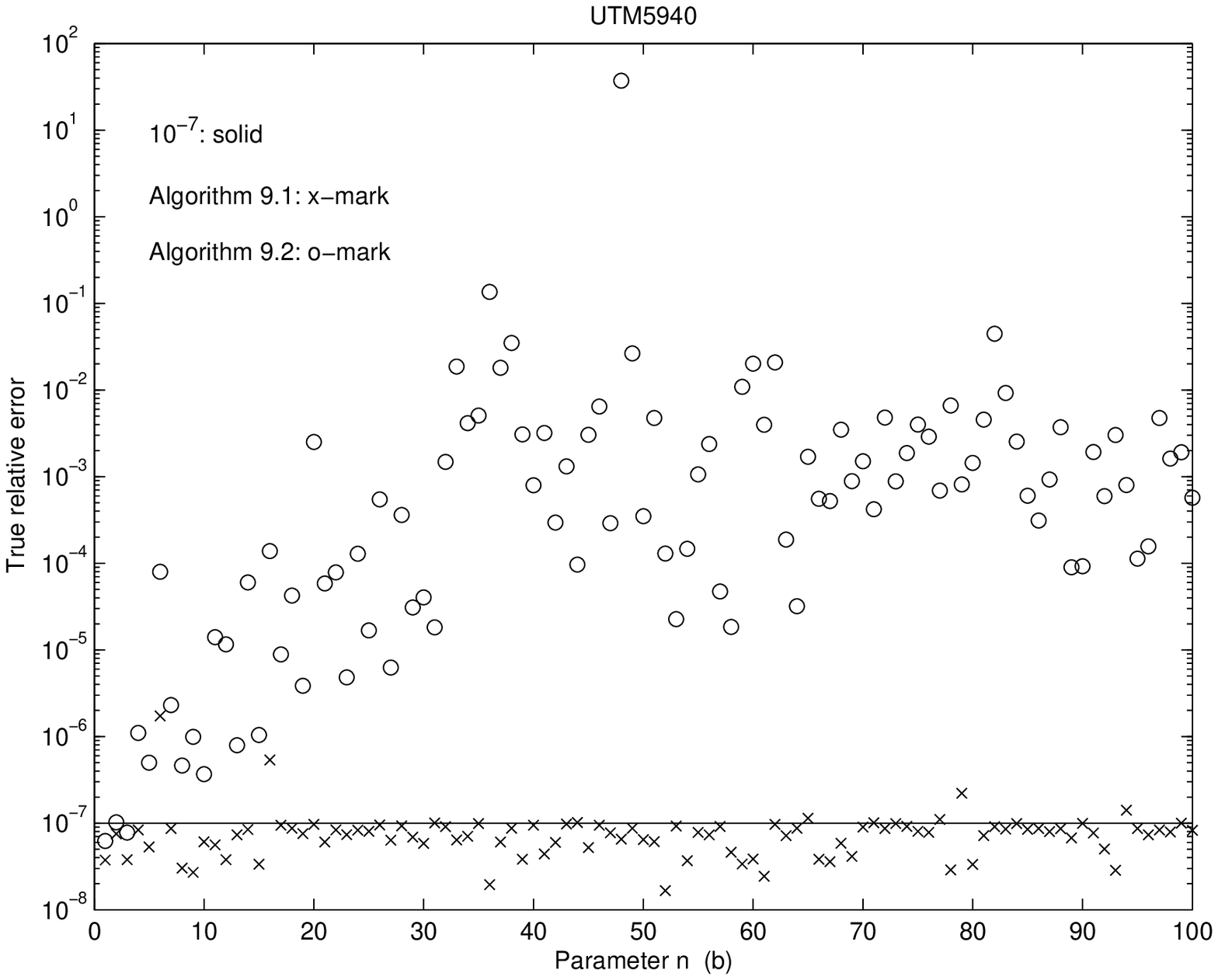,height=6cm}
 }} \caption{utm5940: (a) Graphs of $I_{conv}(n)$ against $n$. BiCGStab took $228$ iterations/$455$
 MVs to converge. Full GMRES
 converged with $176$ MVs.
 (b)
Graphs of $E(n)$ against $n$.
} \label{fig:11-24-3}
\end{figure}

We plot the graphs of $I_{conv}(n)$
in Figures \ref{fig:11-24-1}(a), \ref{fig:11-24-2}(a) and
\ref{fig:11-24-3}(a). For {\it e20r0100} and {\it qc2534},
$I_{conv}(n)$ decreases as $n$ increases. However, the $I_{conv}(n)$
for {\it utm5940} behaves very irregularly due to some of the
$\rho$'s are too small. Recall that ML($n$)BiCGStab performs $1 +
1/n$ matrix-vector multiplications (MVs) per iteration on average.
In terms of the number of MVs, both Algorithms \ref{alg:10-22} and
\ref{alg:10-27} are considerably faster than BiCGStab in all the three experiments.

The graphs of $E(n)$ are plotted in Figures \ref{fig:11-24-1}(b),
\ref{fig:11-24-2}(b) and \ref{fig:11-24-3}(b). It can be seen that
the computed ${\bf r}_k$ in Algorithm \ref{alg:10-27} easily
diverges from its exact counterpart ${\bf b} - {\bf A} {\bf x}_k$.
This divergence
becomes significant when $n \geq  15$ for {\it e20r0100} and $n \geq
4$ for {\it utm5940}.
By contrast, the computed relative errors $\|{\bf r}_k\|_2 / \| {\bf
b}\|_2$ by Algorithm \ref{alg:10-22} well approximate their
corresponding true ones.
Thus, from this point of view, we consider that Algorithm
\ref{alg:10-22} is numerically more stable than Algorithm
\ref{alg:10-27}.
However, Algorithm \ref{alg:10-27}
taken twice to form a predictor-corrector pair can be an efficient and stable algorithm.
We remark that the issues of divergence of computed
residuals and corresponding remedy techniques have been discussed in
detail in \cite{moriya, Svan, vy}.

\subsection{Choice of $n$}\label{subsec:5-20}
In this and the following
subsections, we will focus on Algorithm \ref{alg:10-22}.

From the experiments in \cite{yeungchan} and this paper, we have
observed that ML($n$)BiCGStab behaves more and more robust as $n$ is
increased. So, for an ill-conditioned problem, we would tend to
suggest a large $n$ for ML($n$)BiCGStab.

ML($n$)BiCGStab minimizes $\|{\bf r}_k\|_2$ once every $n$
iterations. The convergence of a well-conditioned problem is usually
accelerated by the minimization steps. So, when a problem is
well-conditioned, we would suggest a small $n$, say, $n \leq 3$. $n
= 2$ may be a good choice since it reduces the MV cost by $25\%$ per
iteration while keeping the minimization performed with a high
frequency.

\begin{figure}[htbp]
\centerline{\hbox{
\psfig{figure=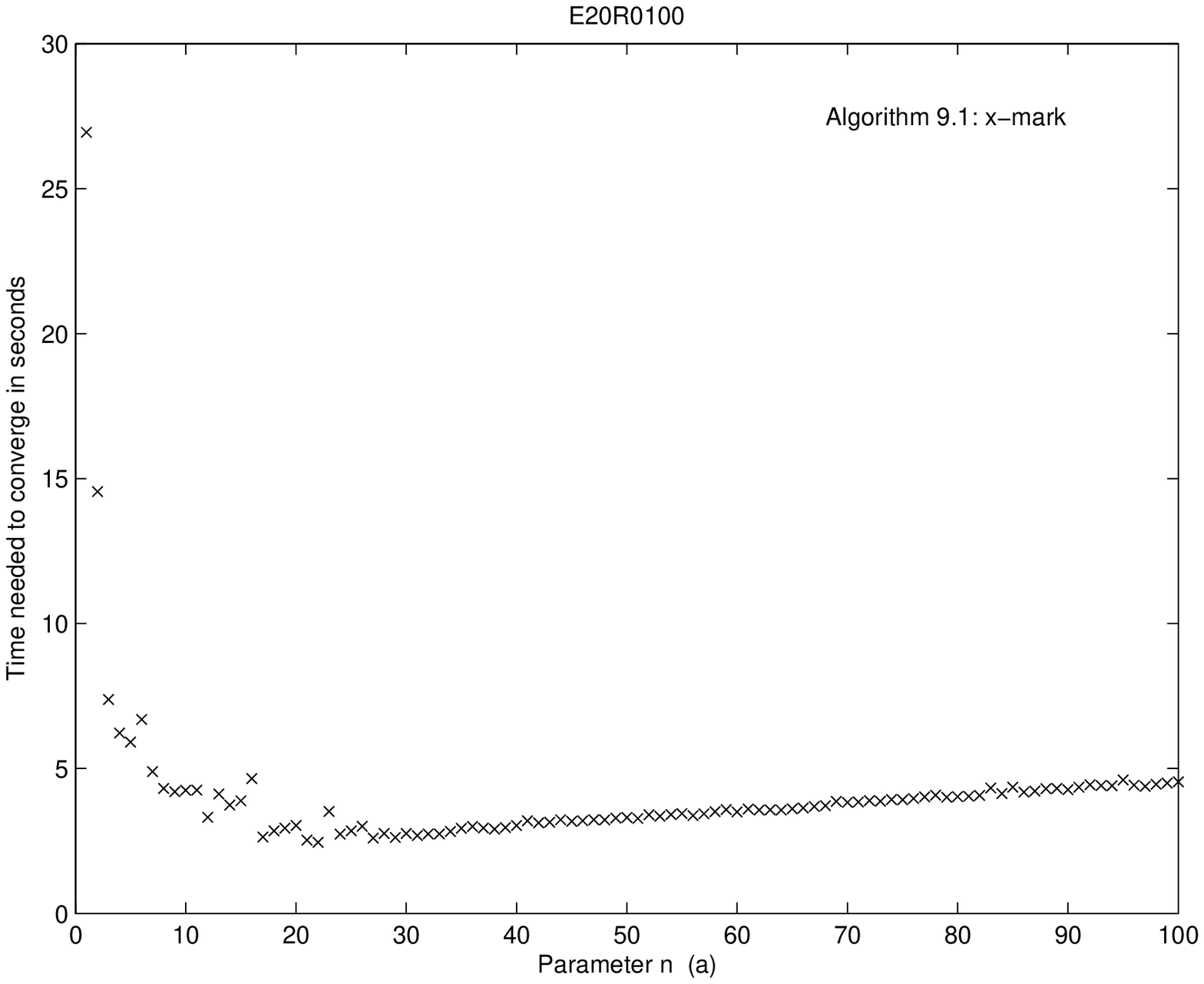,height=6cm}
\psfig{figure=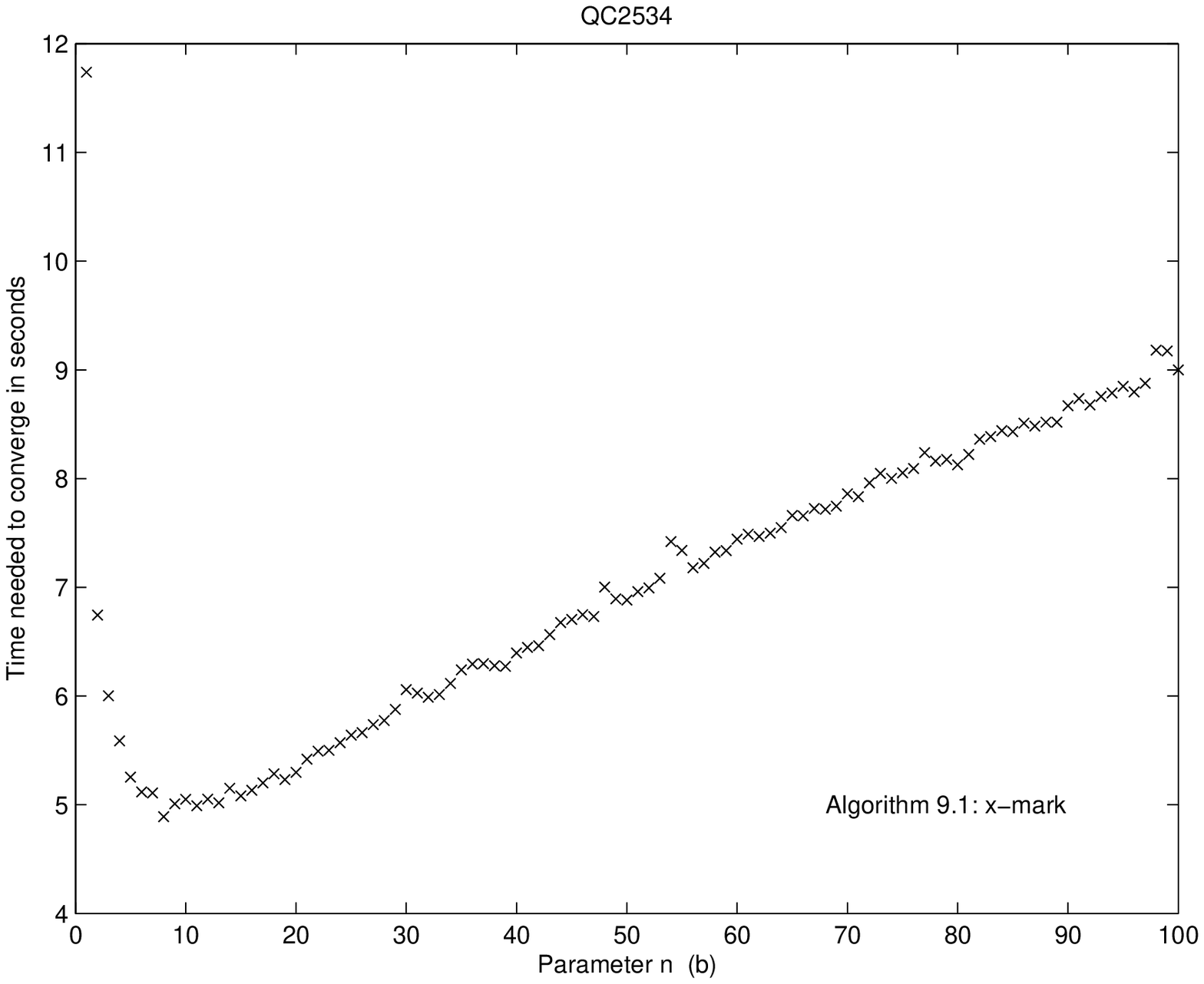,height=6cm}
 }} \caption{(a) e20r0100:
Graph of $T_{conv}(n)$ of Algorithm \ref{alg:10-22} against $n$.
$T_{conv}(n)$ reaches its
 minimum at $n = 22$. (b) qc2534: Graph of $T_{conv}(n)$ of Algorithm \ref{alg:10-22} against $n$.
$T_{conv}(n)$ reaches its
 minimum at $n = 8$.
} \label{fig:12-3-1}
\end{figure}
\begin{figure}[htbp]
\centerline{\hbox{ \psfig{figure=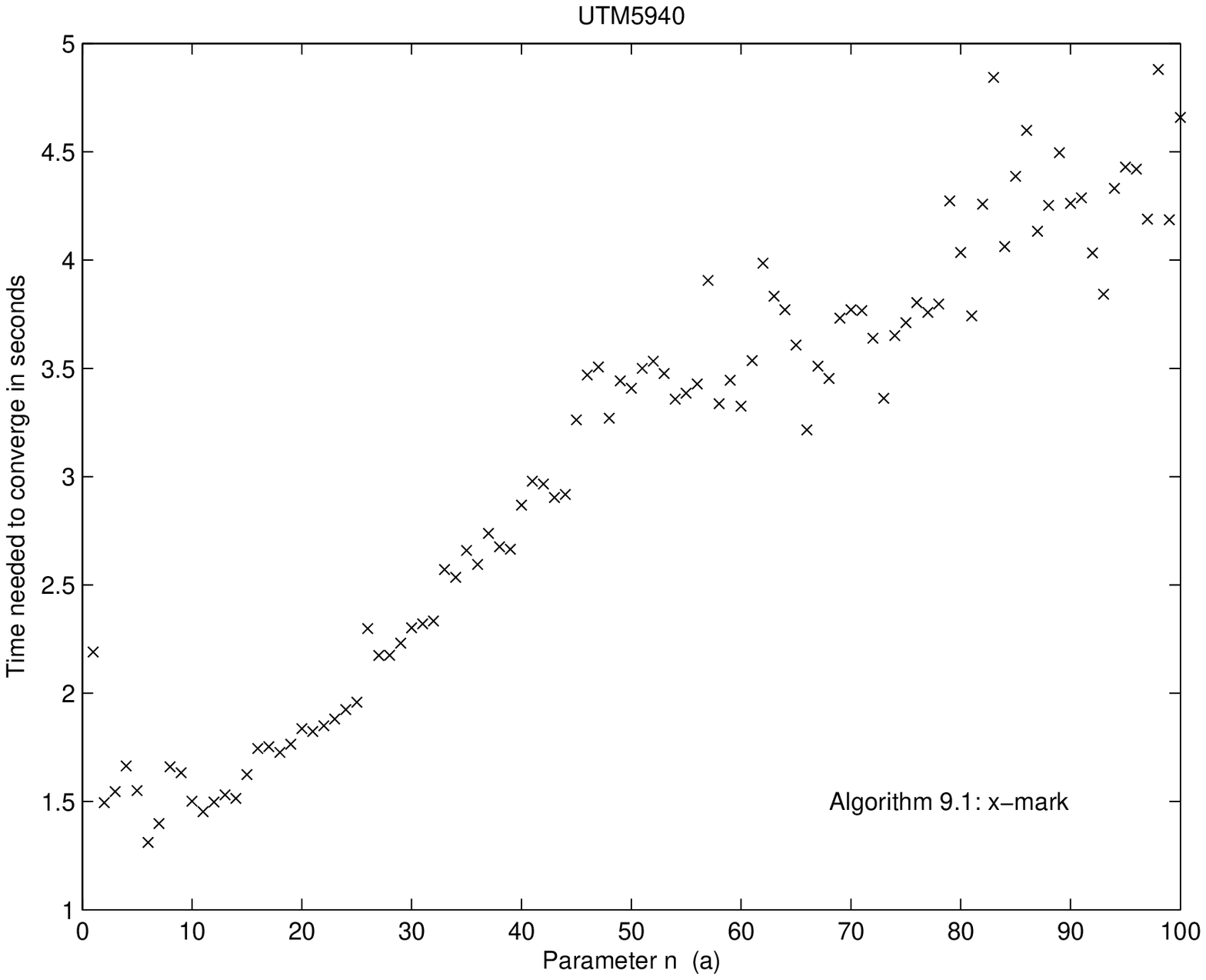,height=6cm}
\psfig{figure=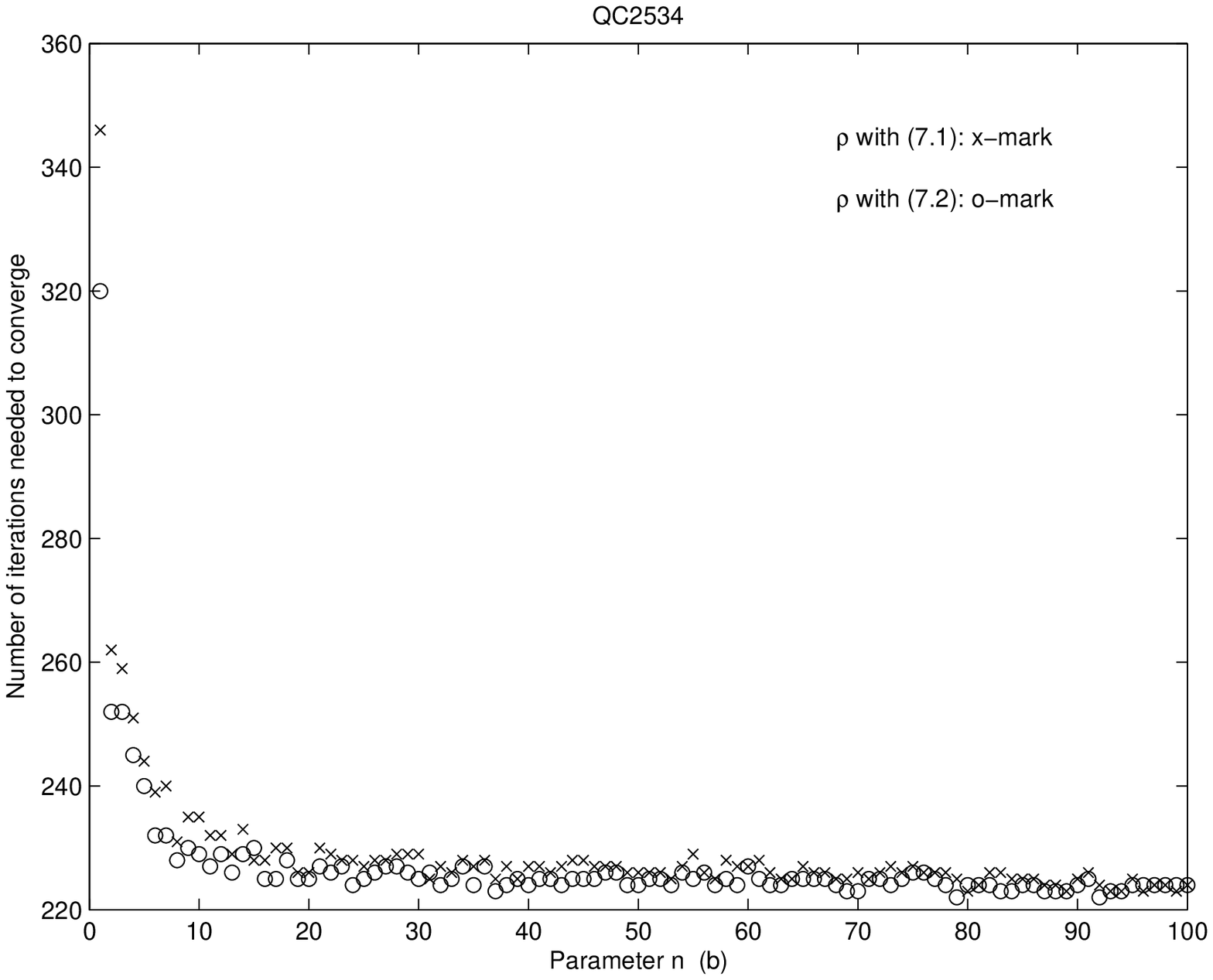,height=6cm}
 }} \caption{
 (a) utm5940:
Graph of $T_{conv}(n)$ of Algorithm \ref{alg:10-22} against $n$.
$T_{conv}(n)$ reaches its
 minimum at $n = 6$. (b) qc2534: Graphs of $I_{conv}(n)$ of Algorithm \ref{alg:10-22}
 against $n$ with choices (\ref{equ:12-6-1}) and (\ref{equ:12-6-2}) for $\rho$
 respectively. In this experiment, we picked $\kappa = 0.7$.
} \label{fig:12-3-2}
\end{figure}

We also plot the graphs of $T_{conv}(n)$ in Figures \ref{fig:12-3-1} and
\ref{fig:12-3-2}(a)
to provide more
information on how $n$ affects the performance of ML($n$)BiCGStab.

\subsection{Choice of $\rho$} 
\begin{figure}[htbp]
\centerline{\hbox{
\psfig{figure=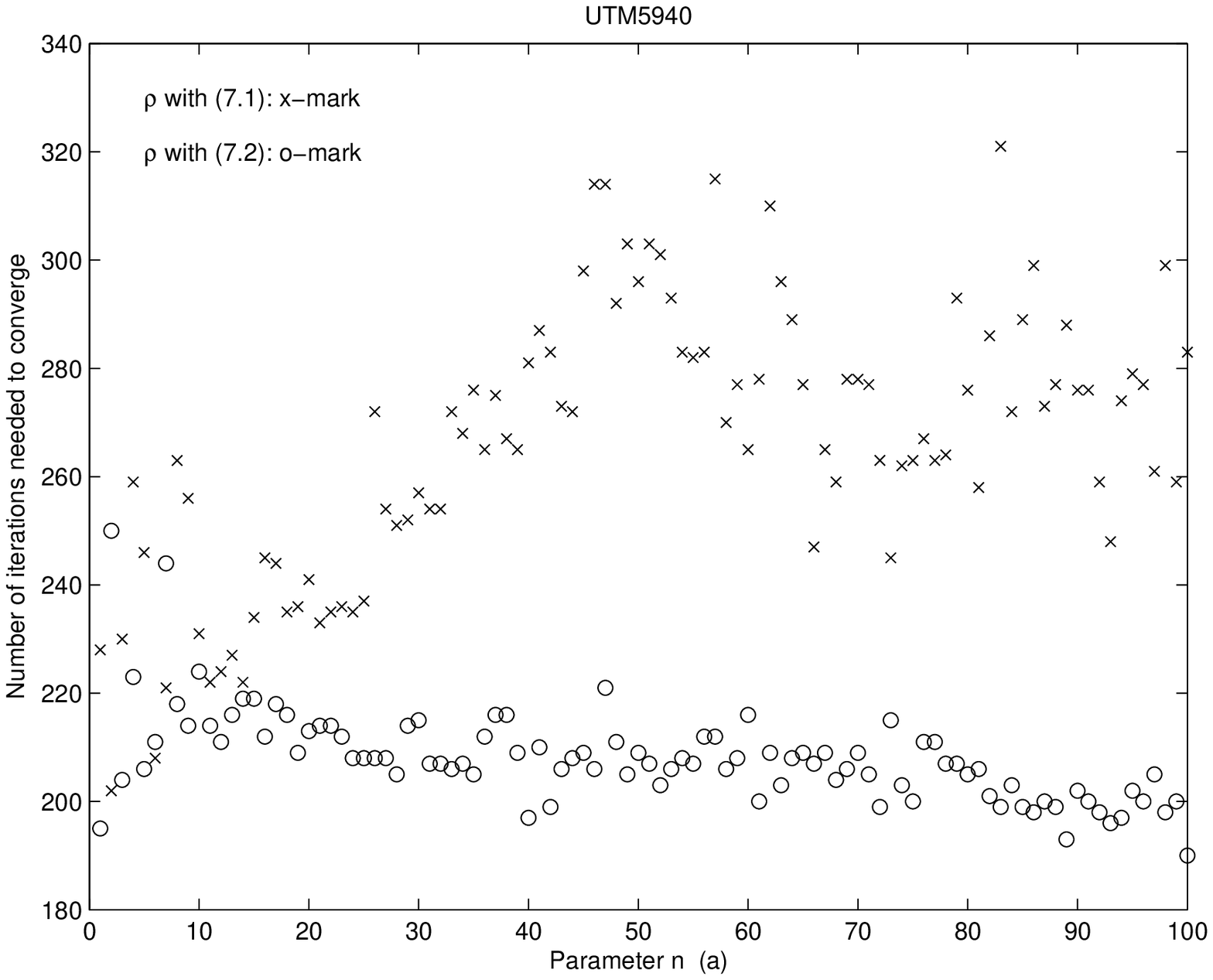,height=6cm }
\psfig{figure=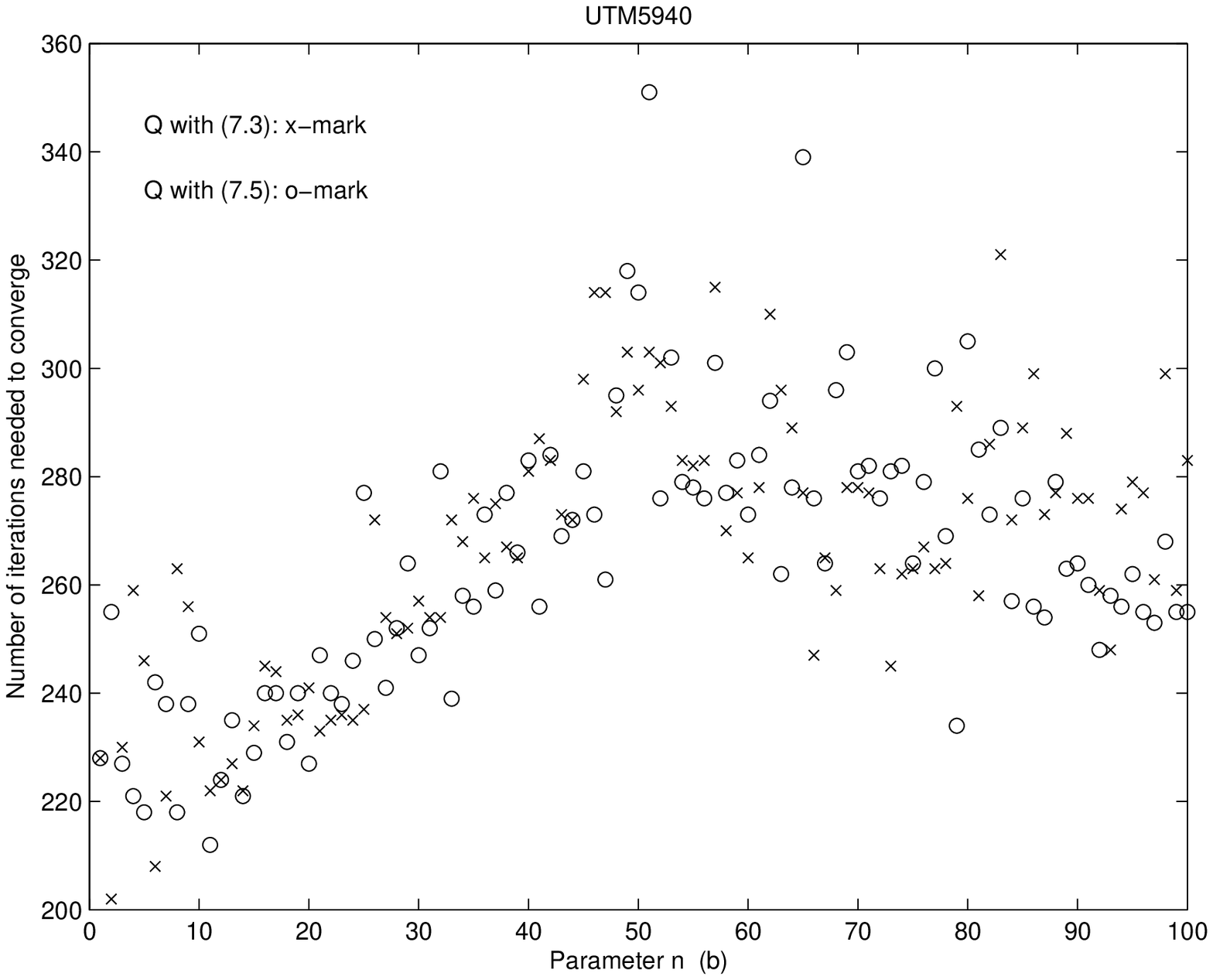,height=6cm}
 }} \caption{
utm5940: (a) Graphs of $I_{conv}(n)$ of Algorithm \ref{alg:10-22}
 against $n$ with choices (\ref{equ:12-6-1}) and (\ref{equ:12-6-2}) for $\rho$
 respectively. In this experiment, we picked $\kappa = 0.7$.
 (b) Graphs of $I_{conv}(n)$ of Algorithm \ref{alg:10-22}
 against $n$ with choices (\ref{equ:12-4-1}) and (\ref{equ:12-4-2}) for ${\bf Q}$ respectively.
 } \label{fig:12-3-3}
\end{figure}

The standard choice for the
$\rho_{j+1}$ in Algorithm \ref{alg:10-22} is
\begin{equation}\label{equ:12-6-1}
\rho_{j+1} = -({\bf A} \widetilde{\bf u}_{jn+1})^H {\bf u}_{jn+1} /
\| {\bf A} \widetilde{\bf u}_{jn+1} \|_2^2.
\end{equation}
This choice of $\rho_{j+1}$ minimizes the $2$-norm of ${\bf
r}_{jn+1} = \rho_{j+1} {\bf A} \widetilde{\bf u}_{jn+1} + {\bf
u}_{jn+1}$, but sometimes can cause instability due to that it can
be very small during an implementation. A remedy as follows has been
suggested in \cite{Svan1}:
\begin{equation}\label{equ:12-6-2}
\begin{array}{l}
 \rho_{j+1}
= -({\bf A} \widetilde{\bf u}_{jn+1})^H {\bf u}_{jn+1}
/ \| {\bf A} \widetilde{\bf u}_{jn+1} \|_2^2;\\
\omega = ({\bf A} \widetilde{\bf u}_{jn+1})^H {\bf u}_{jn+1} / (\|
{\bf A} \widetilde{\bf u}_{jn+1} \|_2\,\, \| {\bf u}_{jn+1}
\|_2);\\
\mbox{if  } |\omega| < \kappa,\,\,\, \rho_{j+1} = \kappa \rho_{j+1}
/ |\omega|; \,\,\, \mbox{end}
\end{array}
\end{equation}
where $\kappa$ is a user-defined parameter. In Figures
\ref{fig:12-3-2}(b) and \ref{fig:12-3-3}(a), we compare the performances of Algorithm
\ref{alg:10-22} with (\ref{equ:12-6-1}) and (\ref{equ:12-6-2})
respectively (we only plot the results of {\it qc2534} and {\it
utm5940}. The result of {\it e20r0100} with $\kappa = 0.1$ is
analogous to Figure \ref{fig:12-3-2}(b)). Also, see the numerical
experiments in \cite{sonn1} for more information about these $\rho$
choices.

\subsection{Choice of $\bf q$'s} We usually pick ${\bf Q} =
[{\bf q}_1, {\bf q}_2, \cdots, {\bf q}_n]$ as
\begin{equation}\label{equ:12-4-1}
{\bf Q} = [{\bf r}_0, randn(N, n-1)]
\end{equation}
for a real problem and
\begin{equation}\label{equ:12-4-11}
{\bf Q} = [{\bf r}_0, randn(N, n-1) + sqrt(-1)*randn(N,n-1)]
\end{equation}
for a complex problem. In our experiments, however, we observed a
comparable performance when we chose
\begin{equation}\label{equ:12-4-2}
{\bf Q} = [{\bf r}_0, sign(randn(N, n-1))].
\end{equation}
or
\begin{equation}\label{equ:12-4-10}
{\bf Q} = [{\bf r}_0, sign(randn(N, n-1)) +  sqrt(-1)* sign(randn(N,
n-1))].
\end{equation}
See Figure \ref{fig:12-3-3}(b) (we only plot the result
 of
 {\it utm5940} for saving space).

 The advantages of
(\ref{equ:12-4-2}) and (\ref{equ:12-4-10}) over
 (\ref{equ:12-4-1}) and (\ref{equ:12-4-11}) are that (i) the storage of ${\bf Q}$
is substantially reduced. In fact, we just need to store the random
signs (except its first column); (ii) an inner product with ${\bf
q}_i$, $2 \leq i \leq n$, is now reduced to a sum without involving
scalar multiplications.

For other choices for ${\bf Q}$, one is referred to \cite{sonn1}.


\section{Concluding Remarks}\label{con} With the help of index
functions, we re-derived the ML($n$)BiCGStab algorithm in
\cite{yeungchan} in a more systematic way.
This time, we have been able to find out and remove some
redundant operations so that the algorithm becomes more efficient. We also
realized that there are $n$ ways to define the
ML($n$)BiCGStab residuals ${\bf r}_k$. Each of the
definitions will lead to a different algorithm. We
presented
two definitions together with
their
associated algorithms, namely, (i) definition (\ref{equ:7-29}),
increasing the degree of $\phi$ at the beginning of an iteration
cycle,
and the associated Algorithm \ref{alg:2}; (ii) definition
(\ref{equ:9-24-1}),
increasing the degree of $\phi$ at the end of an iteration cycle,
and the associated Algorithm \ref{alg:3}.
By comparison, Algorithm \ref{alg:3} is cheaper in storage
and in computational cost, faster to converge, but less stable.
For other definitions of ${\bf r}_k$ that increase the degree of $\phi$
somewhere within a cycle, we expect that
the associated algorithms would lie between
Algorithms \ref{alg:2} and \ref{alg:3} in terms of computational cost, storage and performance.

We proved that the Lanczos-based BiCG/BiCGStab and the Arnoldi-based FOM are essentially methods of the same type.
Both are the extreme cases of ML($n$)BiCG /ML($n$)BiCGStab.

In this paper, we did not assume that ${\bf A}$ is a nonsingular matrix. When a singular system (\ref{equ:7-9-3})
is solved, selecting an appropriate initial guess
$\widehat{\bf x}_0$ is a crucial step. If $\widehat{\bf x}_0$ is selected
such that the affine space $\widehat{\bf x}_0 + span \{ {\bf A}^t
\widehat{\bf r}_0 | t \in {\mathcal N}_0\}$ contains a solution to the
system (\ref{equ:7-9-3})\footnote{For an example where the affine space contains no
solution, consider ${\bf A} = [ 0, 1; 0, 0], {\bf b} = [1, 0]^T$ and
select $\widehat{\bf x}_0 = [ 1, 0]^T$. Note that this linear system
is consistent.}, ML($n$)BiCG will almost surely converge (see
Corollary \ref{cor:4-1}). Otherwise, we shall have $p_{min}(0, {\bf
A}, \widehat{\bf r}_0) = 0$ (see the remark before Corollary
\ref{cor:4-1}) which yields $\det(\widehat{\bf S}_\nu) = 0$ (see the remark after
Proposition \ref{prop:7-11-1}). In this case,
$\prod_{l = 1}^\nu \det(\widehat{\bf S}_l) = 0$ and therefore there is no guarantee that the $LU$-factorizations in the
construction of ML($n$)BiCG exist (see the remark after
Proposition \ref{prop:7-11-1}). As a result, it is likely that
$\|\widehat{\bf r}_k\|_2$ blows up to $\infty$. A similar remark
also applies to ML($n$)BiCGStab.



\section{Appendix} \label{sec:apen} In this section, we present the preconditioned
ML($n$)BiCGStab algorithms together with their Matlab codes.

\subsection{ML($n$)BiCGStab with Definition (\ref{equ:7-29})} \label{sec:11-19-1} The following algorithm is a
preconditioned version of Algorithm \ref{alg:2}. \\

\begin{algorithm}\label{alg:10-22}
{\bf ML($n$)BiCGStab with preconditioning associated with 
(\ref{equ:7-29}). }
 \vspace{.2cm}
\begin{tabbing}
x\=xxx\= xxx\=xxx\=xxx\=xxx\=xxx\=xxx\=xxx\kill \>1. \> Choose an
initial guess ${\bf x}_0$ and $n$ vectors ${\bf q}_1, {\bf q}_2,
\cdots,
{\bf q}_n$. \\
\>2. \>  Compute ${\bf r}_0 = {\bf b} - {\bf A} {\bf x}_0$ and set
${\bf g}_0 = {\bf r}_0$. \\
\>\> Compute $\widetilde{\bf g}_0 = {\bf M}^{-1} {\bf g}_0, \,{\bf
w}_0 = {\bf A}
\widetilde{\bf g}_0,\,\, c_0 = {\bf q}_1^H {\bf w}_0$ and $e_0 = {\bf q}_1^H {\bf r}_0$. \\
\>3. \>  For $j = 0, 1, 2, \cdots$ \\
\>4. \>\> $\alpha_{jn+1} = e_{(j-1)n+n} / c_{(j-1)n+n}$;
\\
\>5. \>\>$
 {\bf u}_{jn+1} =
{\bf r}_{(j-1)n+n} - \alpha_{jn+1} {\bf w}_{(j-1)n+n}$; \\
\>6. \>\>${\bf x}_{jn+1} = {\bf x}_{(j-1)n+n} +
\alpha_{jn+1} \widetilde{\bf g}_{(j-1)n+n}$; \\
\>7. \>\>$\widetilde{\bf u}_{jn+1} = {\bf M}^{-1} {\bf u}_{jn+1}$;\\
\>8. \>\>$ \rho_{j+1} = -({\bf A} \widetilde{\bf u}_{jn+1})^H {\bf
u}_{jn+1} / \| {\bf A} \widetilde{\bf u}_{jn+1} \|_2^2
$;\\
\>9. \>\>${\bf x}_{jn+1} = {\bf x}_{jn+1} - \rho_{j+1} \widetilde{\bf u}_{jn+1}$; \\
\>10. \>\>${\bf r}_{jn+1} = \rho_{j+1} {\bf A} \widetilde{\bf
u}_{jn+1}
+ {\bf u}_{jn+1}$;\\
\>11. \>\>  For $i = 1, 2, \cdots, n-1$ \\
\>12. \>\>\> $f_{jn+i} = {\bf q}^H_{i+1} {\bf u}_{jn+i}$; \\
\>13. \>\>\> If $j \geq 1$\\
\>14. \>\>\>\> $\beta^{(jn+i)}_{(j-1)n+i} = - f_{jn+i} \big/ c_{(j-1)n+i}$;\\
\>15. \>\>\>\>If $i \leq n-2$\\
\>16. \>\>\>\>\>${\bf z}_d = {\bf u}_{jn+i} + \beta^{(jn+i)}_{(j-1)n+i} {\bf d}_{(j-1)n+i}$;\\
\>17. \>\>\>\>\>${\bf g}_{jn+i} = \beta^{(jn+i)}_{(j-1)n+i} {\bf g}_{(j-1)n+i}$;\\
\>18. \>\>\>\>\>${\bf z}_w = \beta^{(jn+i)}_{(j-1)n+i} {\bf w}_{(j-1)n+i}$; \\
\>19. \>\>\>\>\> $\beta^{(jn+i)}_{(j-1)n+i+1} = - {\bf q}^H_{i+2}
{\bf z}_d \big/ c_{(j-1)n+i+1}$; \\
\>20. \>\>\>\>\>   For $s = i+1, \cdots, n - 2$\\
\>21. \>\>\>\>\>\>${\bf z}_d = {\bf z}_d + \beta^{(jn+i)}_{(j-1)n+s} {\bf d}_{(j-1)n+s}$;\\
\>22. \>\>\>\>\>\>${\bf g}_{jn+i} = {\bf g}_{jn+i} + \beta^{(jn+i)}_{(j-1)n+s} {\bf g}_{(j-1)n+s}$;\\
\>23. \>\>\>\>\>\>${\bf z}_w = {\bf z}_w + \beta^{(jn+i)}_{(j-1)n+s} {\bf w}_{(j-1)n+s}$;\\
\>24. \>\>\>\>\>\> $\beta^{(jn+i)}_{(j-1)n+s+1} = - {\bf q}^H_{s+2}
{\bf
z}_d \big/ c_{(j-1)n+s+1}$;\\
\>25. \>\>\>\>\>End\\
\>26. \>\>\>\>\>${\bf g}_{jn+i} = {\bf g}_{jn+i} + \beta^{(jn+i)}_{(j-1) n+n - 1} {\bf g}_{(j-1) n+n - 1}$;\\
\>27. \>\>\>\>\>${\bf z}_w = {\bf z}_w + \beta^{(jn+i)}_{(j-1) n +n-
1}
{\bf w}_{(j-1) n+n - 1}$;\\
\>28. \>\>\>\>\>${\bf z}_w = {\bf r}_{jn+i} + \rho_{j+1}
{\bf z}_w$;\\
\>29. \>\>\>\> Else\\
\>30. \>\>\>\>\>${\bf g}_{jn+i} = \beta^{(jn+i)}_{(j-1) n+n - 1} {\bf g}_{(j-1) n+n - 1}$;\\
\>31. \>\>\>\>\>${\bf z}_w = {\bf r}_{jn+i} + \rho_{j+1}
\beta^{(jn+i)}_{(j-1) n+n - 1} {\bf w}_{(j-1) n+n - 1}
$; \\
\>32. \>\>\>\>           End \\
\>33. \>\>\>\> $\tilde{\beta}^{(jn+i)}_{(j-1) n+n} = - {\bf q}^H_{1} {\bf
z}_w
 \big/c_{(j-1) n+n}
$; \,\,\,\,\,\,\,\, \% $\tilde{\beta}^{(jn+i)}_{(j-1) n+n} = \rho_{j+1} \beta^{(jn+i)}_{(j-1) n+n}$\\
\>34. \>\>\>\> ${\bf z}_w = {\bf z}_w +
\tilde{\beta}^{(jn+i)}_{(j -1)n+n} {\bf w}_{(j-1) n+n}$; \\
\>35. \>\>\>\> ${\bf g}_{jn+i} = {\bf g}_{jn+i} + {\bf z}_w  +
(\tilde{\beta}^{(jn+i)}_{(j-1) n+n}/\rho_{j+1}) {\bf g}_{(j-1) n+n}$; \\
\>36. \>\>\>            Else \\
\>37. \>\>\>\> $\tilde{\beta}^{(jn+i)}_{(j-1) n+n} = - {\bf q}^H_{1} {\bf
r}_{jn+i}
 \big/ c_{(j -1)n+n}
$; \,\,\,\,\,\,\,\, \% $\tilde{\beta}^{(jn+i)}_{(j-1) n+n} = \rho_{j+1} \beta^{(jn+i)}_{(j-1) n+n}$\\
\>38. \>\>\>\> ${\bf z}_w = {\bf r}_{jn+i} +
\tilde{\beta}^{(jn+i)}_{(j-1) n+n} {\bf w}_{(j-1) n+n}$; \\
\>39. \>\>\>\> ${\bf g}_{jn+i} = {\bf z}_w  +
(\tilde{\beta}^{(jn+i)}_{(j-1) n+n}/\rho_{j+1}) {\bf g}_{(j-1) n+n}$;\\
\>40. \>\>\>            End \\
\>41. \>\>\>   For $s = 1, \cdots, i - 1$ \\
\>42. \>\>\>\> $\beta^{(jn+i)}_{jn+s} = - {\bf q}^H_{s+1} {\bf z}_w
 \big/
c_{jn+s}
$; \\
\>43. \>\>\>\>${\bf g}_{jn+i} = {\bf g}_{jn+i} + \beta_{jn+s}^{(jn+i)} {\bf g}_{jn+s}$;\\
\>44. \>\>\>\>${\bf z}_w = {\bf z}_w + \beta_{jn+s}^{(jn+i)} {\bf d}_{jn+s}$;\\
\>45. \>\>\> End \\
\>46. \>\>\>If $i < n-1$\\
\>47.\>\>\>\> $ {\bf d}_{jn+i} =
{\bf z}_w - {\bf u}_{jn+i}$; \\
\>48.\>\>\>\> $c_{jn+i} = {\bf q}_{i+1}^H {\bf d}_{jn+i}$;\\
\>49. \>\>\>\> $\widetilde{\alpha}_{jn+i+1} = f_{jn+i} / c_{jn+i}$;
\,\,\,\,\,\,\,\,\,\,\,\,\,\,\,\,\,\,  \%
$\widetilde{\alpha}_{jn+i+1} = \alpha_{jn+i+1} / \rho_{j+1}$
\\
\>50. \>\>\>\>${\bf u}_{jn+i+1} = {\bf u}_{jn+i} -
\widetilde{\alpha}_{jn+i+1}
{\bf d}_{jn+i}$; \\
\>51. \>\>\>Else\\
\>52.\>\>\>\> $c_{jn+i} = {\bf q}_{i+1}^H ({\bf z}_w - {\bf
u}_{jn+i})
$;\\
\>53. \>\>\>\> $\widetilde{\alpha}_{jn+i+1} = f_{jn+i} / c_{jn+i}$;
\,\,\,\,\,\,\,\,\,\,\,\,\,\,\,\,\,\,  \%
$\widetilde{\alpha}_{jn+i+1} = \alpha_{jn+i+1} / \rho_{j+1}$
\\
\>54. \>\>\>End\\
\>55.\>\>\> $\widetilde{\bf g}_{jn+i} = {\bf M}^{-1}
{\bf g}_{jn+i}$;\\
\>56.\>\>\> ${\bf w}_{jn+i} = {\bf A}
\widetilde{\bf g}_{jn+i}$;\\
\>57. \>\>\>$ {\bf x}_{jn+i+1} =
 {\bf x}_{jn+i} + \rho_{j + 1}
\widetilde{\alpha}_{jn+i+1} \widetilde{\bf g}_{jn+i}$; \\
\>58. \>\>\>${\bf r}_{jn+i+1} = {\bf r}_{jn+i} - \rho_{j + 1}
\widetilde{\alpha}_{jn+i+1}
{\bf w}_{jn+i}$; \\
\>59. \>\> End\\
\>60. \>\> $e_{jn+n} = {\bf q}^H_{1} {\bf r}_{jn+n}$; \\
\>61. \>\> $\tilde{\beta}^{(jn+n)}_{(j-1) n+n} = - e_{jn+n} \big/
c_{(j-1) n+n}$; \,\,\,\,\,\,\,\, \% $\tilde{\beta}^{(jn+n)}_{(j-1) n+n} = \rho_{j+1} \beta^{(jn+n)}_{(j-1) n+n}$\\
\>62. \>\> ${\bf z}_w = {\bf r}_{jn+n} +
\tilde{\beta}^{(jn+n)}_{(j-1)n+n}
{\bf w}_{(j-1)n+n}$; \\
\>63. \>\> ${\bf g}_{jn+n} = {\bf z}_w + (\tilde{\beta}^{(jn+n)}_{(j-1)n+n}/\rho_{j+1})
{\bf g}_{(j-1)n+n}$; \\
\>64. \>\>If $n \geq 2$ \\
\>65. \>\>\> $\beta^{(jn+n)}_{j n+1} = - {\bf q}^H_{2} {\bf z}_w \big/ c_{j n+1}$;\\
\>66. \>\>\>   For $s = 1, \cdots, n - 2$ \\
\>67. \>\>\>\> ${\bf g}_{jn+n} = {\bf g}_{jn+n} + \beta_{jn+s}^{(jn+n)} {\bf g}_{jn+s}$;\\
\>68. \>\>\>\> ${\bf z}_w = {\bf z}_w + \beta_{jn+s}^{(jn+n)} {\bf d}_{jn+s}$; \\
\>69. \>\>\>\> $\beta^{(jn+n)}_{jn+s+1} = - {\bf q}^H_{s+2} {\bf
z}_w \big/ c_{jn+s+1}$;\\
\>70. \>\>\>            End \\
\>71. \>\>\> ${\bf g}_{jn+n} = {\bf g}_{jn+n} + \beta_{jn+n - 1}^{(jn+n)} {\bf g}_{jn+n - 1}$;\\
\>72. \>\>End \\
\>73.\>\> $\widetilde{\bf g}_{jn+n} = {\bf M}^{-1}
{\bf g}_{jn+n}$;\\
\>74.\>\> ${\bf w}_{jn+n} = {\bf A}
\widetilde{\bf g}_{jn+n}$;\\
\>75.\>\> $c_{jn+n} = {\bf q}^H_{1} {\bf w}_{jn+n}$;\\
\>76. \> End
\end{tabbing}
\end{algorithm}

\vspace{.2cm}

{\bf
Matlab code of Algorithm \ref{alg:10-22}}
\begin{tabbing}
x\=xxx\=
xxx\=xxx\=xxx\=xxx\=xxx\=xxx\=xxx\=xxx\=xxx\=xxx\=xxx\=xxx\kill
\>1.\> function
$[x,err,iter,flag] = mlbicgstab(A,x,b,Q,M,max\_it,tol,kappa)$\\
\>2.\\
\>3.\>\% input:\>\>\>$A$:\,\,\, N-by-N matrix. $M$: \,\,\,N-by-N preconditioner matrix. \\
\>4.\>\% \>\>\>$Q$: N-by-n auxiliary matrix $[{\bf
q}_1,\cdots,{\bf q}_n]$. $x$: initial guess.\\
\>5.\>\% \>\>\> $b$:\,\,\, right hand side vector. $max\_it$:\,\,\, maximum number
of iterations.\\
\>6.\>\%\>\>\>$tol$:\,\,\, error tolerance.\\
\>7.\>\%\>\>\>$kappa$:\>\>\> (real number) minimization step controller:\\
\>8.\>\%\>\>\>\>\>\> $kappa = 0$, standard minimization\\
\>9.\>\%\>\>\>\>\>\> $kappa > 0$, Sleijpen-van der Vorst minimization \\
\>10.\>\%output:\>\>\>$x$: solution computed. $err$: error norm. $iter$: number of iterations
performed.\\
\>11.\>\% \>\>\>$flag$:\>\>\> $= 0$, solution found to
tolerance\\
\>12.\>\% \>\>\>\>\>\>  $= 1$, no convergence given $max\_it$ iterations\\
\>13.\>\% \>\>\>\>\>\>$= -1$, breakdown. \\
\>14.\>\% storage:\>\>\> $D$: $N \times (n-2)$ matrix defined only when $n > 2$.\\
\>15.\>\% \>\>\> $G, Q, W$: $N \times n$ matrices. $A, M$: $N \times N$ matrices.\\
\>16.\>\% \>\>\>$x, r, g\_t, u, z, b$: $N \times 1$ matrices. $c$: $1 \times n$ matrix.\\
\>17.\\
\>18.\>\>$N = size(A,2);\,\, n = size(Q,2)$;\\
\>19.\>\>$G = zeros(N,n);\,\, W = zeros(N,n)$; \,\,\,\,\,\%
initialize workspace for ${\bf d}$, ${\bf g}$, ${\bf w}$ and $c$
\\
\>20.\>\>if $n > 2$, \, $D = zeros(N,n-2)$; \, end\\
\>21.\>\>$c = zeros(1,n)$;\>\>\>\>\>\>\% end initialization \\
\>22.\>\>\\
\>23.\>\>$iter = 0;\,\, flag = 1;\,\, bnrm2 = norm(b)$;\\
\>24.\>\> if $bnrm2 == 0.0$,\, $ bnrm2 = 1.0$;\, end \\
\>25.\\
\>26.\>\>$r = b - A*x$;\, $err = norm( r ) / bnrm2$;\\
\>27.\>\> if $err < tol$,\, $flag = 0$;\,\, return, \, end \\
\>28.\\
\>29.\>\> $G(:,n) = r;\,\, g\_t = M \backslash r;\,\, W(:,n) = A*g\_t$;\\
\>30.\>\> $c(n) = Q(:,1)'*W(:,n)$;\\
\>31.\>\> if $c(n) == 0$, \, $flag = -1$; \,\,return,\, end\\
\>32.\>\> $e = Q(:,1)'*r$; \\
\>33.\>\> \\
\>34.\>\>   for $j = 0:max\_it$ \\
\>35.\>\>\> $alpha = e/c(n)$;\\
\>36.\>\>\> $x = x + alpha*g\_t$;\\
\>37.\>\>\> $u = r - alpha*W(:,n)$; \\
\>38.\>\>\>$err = norm( u ) / bnrm2$;\\
\>39.\>\>\> if $err < tol$,\, $flag = 0;\,\, iter = iter + 1$;\,\, return,\, end \\
\>40.\>\>\\
\>41.\>\>\> $g\_t = M \backslash u;\,\, z = A*g\_t;\,\, rho = z'*z$; \\
\>42.\>\>\> if $rho == 0$,\, $flag = -1$;\,\, return,\, end \\
\>43.\>\>\> $omega = z'*u$;\\
\>44.\>\>\> if $omega == 0$, \, $flag = -1$; \,\,return,\, end\\
\>45.\>\>\> $rho = - omega / rho$;\\
\>46.\>\>\> if $kappa > 0$\\
\>47.\>\>\>\>$omega = omega /(norm(z)*norm(u))$;\\
\>48.\>\>\>\> $abs\_om = abs(omega)$;\\
\>49.\>\>\>\>if $abs\_om < kappa$,\, $rho = rho*kappa/abs\_om$; \,end\\
\>50.\>\>\>end\\
\>51.\>\>\> $x = x - rho*g\_t$;\\
\>52.\>\>\> $r = rho*z + u$;\\
\>53.\>\>\> $ err = norm(r)/bnrm2$;\\
\>54.\>\>\> $iter = iter + 1$;\\
\>55.\>\>\> if $err < tol$,\, $flag = 0$;\,\, return,\, end\\
\>56.\>\>\>  if $iter >= max\_it$,\, return,\, end \\
\>57.\>\>\>\\
\>58.\>\>\>  for $i = 1:n-1$\\
\>59.\>\>\>\>$f = Q(:,i+1)'*u$;\\
\>60.\>\>\>\>if $j >= 1$ \\
\>61.\>\>\>\>\>$beta = - f / c(i)$;\\
\>62.\>\>\>\>\>if $i <= n-2$ \\
\>63.\>\>\>\>\>\>$D(:,i) = u + beta*D(:,i)$; \\
\>64.\>\>\>\>\>\>$G(:,i) = beta*G(:,i)$;\\
\>65.\>\>\>\>\>\>$W(:,i) = beta*W(:,i)$;\\
\>66.\>\>\>\>\>\> $beta = - Q(:, i+2)'*D(:,i) / c(i+1)$;\\
\>67.\>\>\>\>\>\>for $s = i+1:n-2$ \\
\>68.\>\>\>\>\>\>\>$D(:,i) = D(:,i) + beta*D(:,s)$;\\
\>69.\>\>\>\>\>\>\> $G(:,i) = G(:,i) + beta*G(:,s)$; \\
\>70.\>\>\>\>\>\>\>$W(:,i) = W(:,i) + beta*W(:,s)$;\\
\>71.\>\>\>\>\>\>\>$beta = - Q(:,s+2)'* D(:,i) / c(s+1)$;\\
\>72.\>\>\>\>\>\> end\\
\>73.\>\>\>\>\>\>$G(:,i) = G(:,i) + beta * G(:, n-1)$;\\
\>74.\>\>\>\>\>\>$W(:,i) = W(:,i) + beta * W(:,n-1)$;\\
\>75.\>\>\>\>\>\>$W(:,i) = r + rho*W(:,i)$;\\
\>76.\>\>\>\>\>else \\
\>77.\>\>\>\>\>\>$G(:,i) = beta*G(:,n-1)$;\\
\>78.\>\>\>\>\>\> $W(:,i) = r + (rho*beta)*W(:,n-1)$; \\
\>79.\>\>\>\>\>end \\
\>80.\>\>\>\>\>$beta = -Q(:,1)'*W(:,i)/c(n)$;\\
\>81.\>\>\>\>\> $W(:,i) = W(:,i) + beta*W(:,n)$; \\
\>82.\>\>\>\>\>$G(:,i) = G(:,i) + W(:,i) + (beta/rho)*G(:,n)$; \\
\>83.\>\>\>\>else \\
\>84.\>\>\>\>\>$beta = -Q(:,1)'*r/c(n)$;\\
\>85.\>\>\>\>\> $W(:, i) = r + beta*W(:,n)$;\\
\>86.\>\>\>\>\>$G(:,i) = W(:,i) + (beta/rho)*G(:,n)$; \\
\>87.\>\>\>\>end \\
\>88.\>\>\>\>for $s = 1:i-1$\\
\>89.\>\>\>\>\>$beta = -Q(:,s+1)'*W(:,i) / c(s)$;\\
\>90.\>\>\>\>\>$G(:,i) = G(:,i) + beta*G(:,s)$;\\
\>91.\>\>\>\>\>$W(:,i) = W(:,i) + beta*D(:,s)$;\\
\>92.\>\>\>\>end \\
\>93.\>\>\>\>if $i < n-1$\\
\>94.\>\>\>\>\>$D(:, i) = W(:,i) - u$;\\
\>95.\>\>\>\>\> $c(i) = Q(:,i+1)'*D(:,i)$;\\
\>96.\>\>\>\>\> if $c(i) == 0$,\, $flag = -1$;\,\, return,\, end \\
\>97.\>\>\>\>\>$alpha = f/c(i)$;\\
\>98.\>\>\>\>\> $u = u - alpha*D(:,i)$;\\
\>99.\>\>\>\>else \\
\>100.\>\>\>\>\>$c(i) = Q(:,i+1)'*(W(:,i) - u)$;\\
\>101.\>\>\>\>\> if $c(i) == 0$,\, $flag = -1$;\,\, return,\, end \\
\>102.\>\>\>\>\>$alpha = f/c(i)$;\\
\>103.\>\>\>\>end \\
\>104.\>\>\>\>$g\_t = M \backslash \, G(:,i);\,\, W(:,i) = A*g\_t$;\\
\>105.\>\>\>\>$alpha = rho*alpha$;\\
\>106.\>\>\>\> $x = x + alpha*g\_t$;\\
\>107.\>\>\>\> $r = r - alpha*W(:,i)$; \\
\>108.\>\>\>\>$err = norm(r)/bnrm2$;\\
\>109.\>\>\>\> $iter = iter + 1$;\\
\>110.\>\>\>\>if $err < tol$,\, $flag = 0$;\,\, return,\, end\\
\>111.\>\>\>\>if $iter >= max\_it$,\, return,\, end \\
\>112.\>\>\>    end \\
\>113.\>\>\>$e = Q(:,1)'*r;\,\, beta = - e/c(n)$; \\
\>114.\>\>\>$W(:,n) = r + beta*W(:,n)$;\\
\>115.\>\>\> $G(:,n) = W(:,n) + (beta/rho)*G(:, n)$; \\
\>116.\>\>\>if $n >= 2$ \\
\>117.\>\>\>\>$beta = - Q(:,2)'*W(:,n) / c(1)$;\\
\>118.\>\>\>\> for $s = 1 : n-2$ \\
\>119.\>\>\>\>\>$G(:,n) = G(:,n) + beta*G(:,s)$;\\
\>120.\>\>\>\>\>$W(:,n) = W(:,n) + beta*D(:, s)$;\\
\>121.\>\>\>\>\>$beta = - Q(:,s+2)'*W(:,n) / c(s+1)$;\\
\>122.\>\>\>\> end \\
\>123.\>\>\>\>$G(:,n) = G(:,n) + beta*G(:,n-1)$;\\
\>124.\>\>\> end\\
\>125.\>\>\>$g\_t = M \backslash \, G(:,n);\,\, W(:, n) = A*g\_t$;\\
\>126.\>\>\> $c(n) = Q(:,1)'*W(:,n)$;\\
\>127.\>\>\>  if $c(n) == 0$,\, $ flag = -1$;\,\, return,\, end \\
\>128.\>\>end
\end{tabbing}


\subsection{ML($n$)BiCGStab with Definition (\ref{equ:9-24-1})} \label{sec:11-19-2} The following algorithm is a
preconditioned version of
Algorithm \ref{alg:3}.\\

\begin{algorithm} \label{alg:10-27} {\bf ML($n$)BiCGStab with preconditioning associated
with (\ref{equ:9-24-1}).}
 \vspace{.1cm}
\begin{tabbing}
x\=xxx\= xxx\= xxx\= xxx\= xxx\= xxx\= xxx\=
xxx\=xxx\=xxx\=xxx\=xxx\=xxx\kill \>1. \> Choose an initial guess
${\bf x}_0$ and $n$ vectors ${\bf q}_1, {\bf q}_2, \cdots,
{\bf q}_n$. \\
\>2. \>  Compute ${\bf r}_0 = {\bf b} - {\bf A} {\bf x}_0, \,\,
\widetilde{\bf g}_0 = {\bf M}^{-1} {\bf r}_0,\,\, {\bf w}_0 = {\bf
A}
\widetilde{\bf g}_0,\,\, c_0 = {\bf q}^H_{1} {\bf w}_0$ and $e_0 = {\bf q}_1^H {\bf r}_0$.\\
\>3. \>  For $j = 0, 1, 2, \cdots$ \\
\>4. \>\>For $i = 1, 2, \cdots, n-1$ \\
\>5. \>\>\> $  \alpha_{jn+i} = e_{jn+i-1} /
c_{jn+i-1}$;\\
\>6. \>\>\>  $ {\bf x}_{jn+i} = {\bf x}_{jn+i-1} + \alpha_{jn+i}
\widetilde{\bf g}_{jn+i-1}$; \,\,\,\,\,\, \% $\widetilde{\bf g} = {\bf M}^{-1} {\bf g}$\\
\>7. \>\>\>  $ {\bf r}_{jn+i} = {\bf r}_{jn+i-1} - \alpha_{jn+i}
{\bf w}_{jn+i-1}$; \\
\>8.\>\>\>$e_{jn+i} = {\bf q}^H_{i+1}
 {\bf r}_{jn+i}$;\\
\>9. \>\>\>If $j \geq 1$ \\
\>10. \>\>\>\> $\tilde{\beta}^{(jn+i)}_{(j-1)n+i} = - e_{jn+i} \big/
c_{(j-1)n+i}$; \,\,\,\,\,\,\, \% $\tilde{\beta}^{(jn+i)}_{(j-1)n+i} = \rho_j \beta^{(jn+i)}_{(j-1)n+i}$ \\
\>11. \>\>\>\>${\bf z}_w = {\bf r}_{jn+i} +
\tilde{\beta}^{(jn+i)}_{(j-1)n+i}
{\bf w}_{(j-1)n+i}$; \\
\>12. \>\>\>\>$\widetilde{\bf g}_{jn+i} = \tilde{\beta}^{(jn+i)}_{(j-1)n+i}
\widetilde{\bf
g}_{(j-1)n+i}$;\\
\>13. \>\>\>\>   For $s = i +1, \cdots, n - 1$ \\
\>14. \>\>\>\>\> $\tilde{\beta}^{(jn+i)}_{(j-1)n+s} = - {\bf
q}^H_{s+1}
 {\bf z}_w \big/
c_{(j-1)n+s}$; \,\,\,\,\,\,\, \% $\tilde{\beta}^{(jn+i)}_{(j-1)n+s} = \rho_j \beta^{(jn+i)}_{(j-1)n+s}$\\
\>15. \>\>\>\>\>${\bf z}_w = {\bf z}_w +
\tilde{\beta}^{(jn+i)}_{(j-1)n+s}
{\bf w}_{(j-1)n+s}$; \\
\>16. \>\>\>\>\>$\widetilde{\bf g}_{jn+i} = \widetilde{\bf g}_{jn+i}
+ \tilde{\beta}^{(jn+i)}_{(j-1)n+s}
\widetilde{\bf g}_{(j-1)n+s}$; \\
\>17. \>\>\>\>End\\
\>18. \>\>\>\>$\displaystyle{\widetilde{\bf g}_{jn+i} =
{\bf M}^{-1}
{\bf z}_w +
\frac{1}{\rho_j} \widetilde{\bf g}_{jn+i}
}$;\\
\>19. \>\>\> Else \\
\>20. \>\>\>\>$\widetilde{\bf g}_{jn+i} = {\bf M}^{-1}
{\bf r}_{jn+i}$;\\
\>21. \>\>\> End\\
\>22. \>\>\>${\bf w}_{jn+i} = {\bf A} \widetilde{\bf g}_{jn+i}$;\\
\>23. \>\>\>   For $s = 0, \cdots, i - 1$\\
\>24. \>\>\>\>        $\beta^{(jn+i)}_{jn+s} = - {\bf q}^H_{s+1}
{\bf w}_{jn+i}
 \big/
c_{jn+s}$;\\
\>25. \>\>\>\>${\bf w}_{jn+i} = {\bf w}_{jn+i} + \beta^{(jn+i)}_{jn+s} {\bf w}_{jn+s}$;\\
\>26. \>\>\>\>$\widetilde{\bf g}_{jn+i} = \widetilde{\bf g}_{jn+i} +
\beta^{(jn+i)}_{jn+s}
\widetilde{\bf g}_{jn+s}$;\\
\>27. \>\>\>End \\
\>28. \>\>\>$c_{jn+i} = {\bf q}^H_{i+1} {\bf w}_{jn+i}$;\\
\>29. \>\> End \\
\>30. \>\> $  \alpha_{jn+n} = e_{jn+n-1} /
c_{jn+n-1}$;\\
\>31. \>\>  ${\bf x}_{jn+n} = {\bf x}_{jn+n-1} + \alpha_{jn+n}
\widetilde{\bf g}_{jn+n-1}$;\\
\>32. \>\>  $ {\bf u}_{jn+n} = {\bf r}_{jn+n-1} - \alpha_{jn+n}
{\bf w}_{jn+n-1}$; \\
\>33. \>\>  $ \widetilde{\bf u}_{jn+n} = {\bf M}^{-1} {\bf u}_{jn+n}$;\\
\>34. \>\> $\rho_{j+1} = - ({\bf A} \widetilde{\bf u}_{jn+n})^H {\bf
u}_{jn+n} /
\| {\bf A} \widetilde{\bf u}_{jn+n} \|_2^2$;\\
\>35. \>\> ${\bf x}_{jn+n} = {\bf x}_{jn+n} - \rho_{j+1}
\widetilde{\bf u}_{jn+n}$;\\
\>36. \>\> ${\bf r}_{jn+n} = \rho_{j+1} {\bf A} \widetilde{\bf
u}_{jn+n} +
{\bf u}_{jn+n}$;\\
\>37. \>\>$e_{jn+n} = {\bf q}^H_{1} {\bf r}_{jn+n}$; \\
\>38. \>\>  $\tilde{\beta}^{(jn+n)}_{(j-1) n+n} = - e_{jn+n} \big/
c_{(j-1) n+n}$; \,\,\,\,\,\, \% $\tilde{\beta}^{(jn+n)}_{(j-1) n+n} = \rho_{j+1} \beta^{(jn+n)}_{(j-1) n+n}$ \\
\>39. \>\>${\bf z}_w = {\bf r}_{jn+n} +
\tilde{\beta}^{(jn+n)}_{(j-1) n+n} {\bf w}_{(j-1) n+n}$;\\
\>40. \>\>$\widetilde{\bf g}_{jn+n} = \tilde{\beta}^{(jn+n)}_{(j-1) n+n}
\widetilde{\bf
g}_{(j-1) n+n}$; \\
\>41. \>\>   For $s = 1, \cdots, n - 1$ \\
\>42. \>\>\>        $\tilde{\beta}^{(jn+n)}_{jn+s} = - {\bf q}^H_{s+1} {\bf
z}_w \big/
c_{jn+s}$; \,\,\,\,\,\,\,\,\,\,\, \% $\tilde{\beta}^{(jn+n)}_{jn+s} = \rho_{j+1} \beta^{(jn+n)}_{jn+s}$ \\
\>43. \>\>\>${\bf z}_w = {\bf z}_w +
\tilde{\beta}^{(jn+n)}_{jn+s} {\bf w}_{jn+s}$; \\
\>44. \>\>\>$\widetilde{\bf g}_{jn+n} = \widetilde{\bf g}_{jn+n} +
\tilde{\beta}^{(jn+n)}_{jn+s}
\widetilde{\bf g}_{jn+s}$;\\
\>45. \>\>  End\\
\>46. \>\>$\displaystyle{\widetilde{\bf g}_{jn+n} = {\bf M}^{-1}
 {\bf z}_w + \frac{1}{\rho_{j+1}} \widetilde{\bf g}_{jn+n}}$;\\
\>47. \>\>${\bf w}_{jn+n} = {\bf A} \widetilde{\bf g}_{jn+n}
$; \\
\>48. \>\>$c_{jn+n} = {\bf q}^H_{1} {\bf w}_{jn+n}$;\\
\>49. \> End
\end{tabbing}
\end{algorithm}

\vspace{.2cm}

{\bf
Matlab code of Algorithm \ref{alg:10-27}}
\begin{tabbing}
x\=xxx\=
xxx\=xxx\=xxx\=xxx\=xxx\=xxx\=xxx\=xxx\=xxx\=xxx\=xxx\=xxx\kill \>1.
\>function $[x,err,iter,flag] = mlbicgstab(A,x,b,Q,M,max\_it,tol,kappa)$\\
\>2.\\
\>3.\>\% input: \>\>\>$A$: N-by-N matrix. $M$: N-by-N preconditioner matrix \\
\>4.\>\%\>\>\>$Q$: N-by-n auxiliary matrix with
columns ${\bf q}_1, \cdots, {\bf
q}_n$.\\
\>5.\>\%\>\>\>$x$: initial guess. $b$: right hand side vector.\\
\>6.\>\%\>\>\>$max\_it$: maximum
number of iterations. $tol$: error tolerance. \\
\>7.\>\%\>\>\>$kappa$: \,\,\,\,\,\,\,\,\,(real number) minimization step controller:\\
\>8.\>\%\>\>\>\>\>\>zero, standard minimization\\
\>9.\>\%\>\>\>\>\>\>positive, Sleijpen-van der Vorst minimization\\
\>10.\>\% output:\>\>\>\> $x$: solution computed. $err$: error norm.\\
\>11.\>\%\>\>\>\>$iter$: number of
iterations performed.\\
\>12.\>\%\>\>\>\>$flag$: \,\, 0 = solution
found to tolerance\\
\>13.\>\%\>\>\>\>\>\>\,\,1 = no convergence given $max\_it$\\
\>14.\>\%\>\>\>\>\>\> -1 = breakdown\\
\>15.\>\% storage: \>\>\>\>$c$: $1 \times n$ matrix. $x, r, b, u\_t, z$: N-by-1 matrices.\\
\>16.\>\% \>\>\>\>$A, M$: N-by-N matrices. $Q, G, W$: N-by-n matrices. \\
\>17.\\
\>18.\>\>$N = size(A,2); \,\,n = size(Q,2)$;\\
\>19.\>\>$G = zeros(N,n);\,\, W = zeros(N,n)$;\,\,\,\,\, \%
initialize
workspace for $\widetilde{\bf g}$, ${\bf w}$ and $c$\\
\>20.\>\>$c = zeros(1,n)$;
\,\,\,\,\,\,\,\,\,\,\,\,\, \% end initialization \\
\>21.\>\>\\
\>22.\>\>$iter = 0;\,\, flag = 1;\,\, bnrm2 = norm(b)$; \\
\>23.\>\>if $bnrm2 == 0.0$,\, $bnrm2 = 1.0$;\, end \\
\>24.\>\>$r = b - A*x; \,\, err = norm( r ) / bnrm2$;\\
\>25.\>\>if $err < tol$,\, $flag = 0$;\,\, return,\, end \\
\>26.\>\>\\
\>27.\>\>$G(:,1) = M \backslash r; \,\, W(:,1) = A*G(:,1); \,\, c(1) = Q(:,1)'*W(:,1)$; \\
\>28.\>\>if $c(1) == 0$,\, $flag = -1$;\,\, return,\, end \\
\>29.\>\>$e = Q(:,1)'*r$; \\
\>30.\>\>\\
\>31.\>\>  for $j = 0:max\_it$\\
\>32.\>\>\>for $i = 1:n-1$ \\
\>33.\>\>\>\> $alpha = e / c(i)$;\\
\>34.\>\>\>\> $x = x + alpha*G(:, i)$; \\
\>35.\>\>\>\> $r = r - alpha*W(:, i)$;\\
\>36.\>\>\>\> $err = norm(r)/bnrm2$;\\
\>37.\>\>\>\> $iter = iter + 1$;\\
\>38.\>\>\>\>if $err < tol$,\, $flag = 0$;\,\, return,\, end\\
\>39.\>\>\>\>if $iter >= max\_it$,\, return,\, end \\
\>40.\>\>\>\>\\
\>41.\>\>\>\>$e = Q(:,i+1)'*r$;\\
\>42.\>\>\>\>if $j >= 1$\\
\>43.\>\>\>\>\>$beta = -e / c(i+1)$; \\
\>44.\>\>\>\>\>$W(:,i+1) = r + beta*W(:,i+1)$;\\
\>45.\>\>\>\>\>$G(:,i+1) = beta*G(:,i+1)$;\\
\>46.\>\>\>\>\> for $s = i+1 : n-1$ \\
\>47.\>\>\>\>\>\>$beta = -Q(:,s+1)'*W(:,i+1)/c(s+1)$; \\
\>48.\>\>\>\>\>\>$W(:,i+1) = W(:,i+1) + beta*W(:,s+1)$; \\
\>49.\>\>\>\>\>\>$G(:,i+1) = G(:,i+1) + beta*G(:,s+1)$; \\
\>50.\>\>\>\>\>end \\
\>51.\>\>\>\>\>$G(:,i+1) = (M \backslash W(:,i+1)) + (1/rho)*G(:,i+1)$; \\
\>52.\>\>\>\>else\\
\>53.\>\>\>\>\>$G(:,i+1) = M \backslash r$;\\
\>54.\>\>\>\>end \\
\>55.\>\>\>\>$W(:,i+1) = A*G(:,i+1)$;\\
\>56.\>\>\>\>for $s = 0:i-1$ \\
\>57.\>\>\>\>\>$beta = -Q(:,s+1)'*W(:,i+1) / c(s+1)$; \\
\>58.\>\>\>\>\>$W(:,i+1) = W(:,i+1) + beta*W(:,s+1)$;\\
\>59.\>\>\>\>\>$G(:,i+1) = G(:,i+1) + beta*G(:,s+1)$;\\
\>60.\>\>\>\>end \\
\>61.\>\>\>\>$c(i+1) = Q(:,i+1)'*W(:,i+1)$; \\
\>62.\>\>\>\>if $c(i+1) == 0$,\, $flag = -1$;\,\, return,\, end \\
\>63.\>\>\>end\\
\>64.\>\>\> $alpha = e / c(n)$; \\
\>65.\>\>\> $x = x + alpha*G(:, n)$;\\
\>66.\>\>\> $r = r - alpha*W(:,n)$;\\
\>67.\>\>\>$err = norm(r)/bnrm2$;\\
\>68.\>\>\>if $err < tol$,\, $flag = 0; \,\, iter = iter + 1$;\, return,\, end \\
\>69.\>\>\>$u\_t = M \backslash r; \,\, z = A*u\_t$;\, $rho = z'*z$;\\
\>70.\>\>\>if $rho == 0$,\, $flag = -1$;\, return,\, end\\
\>71.\>\>\>$omega = z'*r$;\\
\>72.\>\>\>if $omega == 0$,\, $flag = -1$;\, return,\, end\\
\>73.\>\>\>$rho = - omega / rho$;\\
\>74.\>\>\>if $kappa > 0$\\
\>75.\>\>\>\>$omega = omega / (norm(z)*norm(r))$;\\
\>76.\>\>\>\>$abs\_om = abs(omega)$;\\
\>77.\>\>\>\>if $abs\_om < kappa$\\
\>78.\>\>\>\>\>$rho = rho*kappa/abs\_om$;\\
\>79.\>\>\>\>end\\
\>80.\>\>\>end\\
\>81.\>\>\>$x = x - rho*u\_t$; \\
\>82.\>\>\>$r = r + rho*z$;\\
\>83.\>\>\>$err = norm(r)/bnrm2$;\\
\>84.\>\>\>$iter = iter + 1$;\\
\>85.\>\>\>if $err < tol$,\, $ flag = 0$;\,\, return,\, end\\
\>86.\>\>\>if $iter >= max\_it$,\, return,\, end \\
\>87.\>\>\>\\
\>88.\>\>\>$e = Q(:,1)'*r; \,\, beta = - e/c(1)$; \\
\>89.\>\>\>$W(:,1) = r + beta*W(:,1)$;\\
\>90.\>\>\>$G(:,1) = beta*G(:,1)$; \\
\>91.\>\>\> for $s = 1:n-1$ \\
\>92.\>\>\>\>$beta = -Q(:,s+1)'*W(:,1)/c(s+1)$; \\
\>93.\>\>\>\>$W(:,1) = W(:,1) + beta*W(:,s+1)$; \\
\>94.\>\>\>\>$G(:,1) = G(:,1) + beta*G(:,s+1)$;\\
\>95.\>\>\> end \\
\>96.\>\>\>$G(:,1) = (M \backslash W(:,1)) + (1/rho)*G(:,1);\,\, W(:,1) = A*G(:,1)$;\\
\>97.\>\>\>$c(1) = Q(:,1)'*W(:,1)$; \\
\>98.\>\>\>if $c(1) == 0$,\, $flag = -1$;\, return,\, end\\
\>99.\>\>  end
\end{tabbing}

\vspace{.2cm}

{\bf
A sample run of ML($n$)BiCGstab}
\begin{tabbing}
x\=xxx\=
xxx\=xxx\=xxx\=xxx\=xxx\=xxx\=xxx\=xxx\=xxx\=xxx\=xxx\=xxx\kill \>1.
\>$N = 100$; $A = randn(N)$; $M = randn(N)$; $b = randn(N, 1)$;\\
\>2.\>$n = 10$; $kappa = 0.7$; $tol = 10^{-7}$; $max\_it = 3*N$;\\
\>3.\>$Q = sign(randn(N, n))$; $x = zeros(N, 1)$; $Q(:,1) = b-A*x$;\\
\>4.\>$[x,err,iter,flag] = mlbicgstab(A,x,b,Q,M,max\_it,tol, kappa)$;
\end{tabbing}

\section*{Acknowledgements} Thanks will be added.
%
%
%


\end{document}